%% file: OMINroot.tex
\documentclass[12pt,a4paper]{book}
\usepackage{amssymb}     % \mathbb, \mathfrak, \square (replaces the MSBM/EUFM font loading)
\DeclareFontShape{OMX}{cmex}{m}{n}{<-> sfixed * cmex10}{}
\usepackage{graphicx}
\usepackage{makeidx}
\usepackage[hidelinks]{hyperref}

\newenvironment{proof}%
    {\par\noindent{\it Proof. }\nopagebreak\normalsize}%
    {\linebreak[2]\hspace*{\fill}$\square$\ifdim\lastskip<12pt
    \removelastskip \penalty-200 \vskip12pt \fi}

\newtheorem{theorem}{Theorem}[chapter]
\newtheorem{proposition}[theorem]{Proposition}
\newtheorem{corollary}[theorem]{Corollary}
\newtheorem{lemma}[theorem]{Lemma}
\newtheorem{definition}[theorem]{Definition}
\newlength{\retraitexo}
\newlength{\largeurexo}
\newtheorem{exo}[theorem]{Exercise}%
\newcommand{\debexo}{\begin{exo} \rm\parshape=2 0cm \textwidth
\retraitexo
\largeurexo }
\newcommand{\finexo}{\par \end{exo}}

\renewcommand{\Bbb}{\mathbb}

\newlength{\dessinw}\newlength{\dessinh}
\def\dessin #1 by #2 (#3){%
  \setlength{\dessinw}{#1}\setlength{\dessinh}{#2}\dessinimg{#3}}
\newlength{\dessinmargin}
\newcommand{\dessinimg}[1]{%
  \vbox to \dessinh{\vfill
    \hbox to \dessinw{\kern-\dessinmargin
      \raisebox{-\dessinmargin}[0pt][0pt]{\includegraphics{#1}}\hss}}}

\newcommand{\illustration}[1]{\raisebox{-\dessinmargin}[0pt][0pt]{\kern-\dessinmargin\includegraphics{#1}}}

\def\clos{\mathop{\mbox{\rm clos}}\nolimits}

\def\et{\mathrel{\mbox{\rm\ and\ }}}
\def\ou{\mathrel{\mbox{\rm\ or\ }}}
\def\shitem{\setlength{\itemsep}{0pt plus1pt minus1pt}}

\let\origtheindex\theindex
\renewcommand{\theindex}{\origtheindex\addcontentsline{toc}{chapter}{Index}}
\makeindex

\begin{document}
\title{AN INTRODUCTION\\TO O-MINIMAL GEOMETRY}
\author{Michel COSTE\\ Institut de Recherche Math\'ematique de Rennes\\
michel.coste@univ-rennes.fr}
\date{November 1999}
\hypersetup{pageanchor=false}
\maketitle
\hypersetup{pageanchor=true}
\setcounter{page}{1}
\include{OMINintro}
\tableofcontents
\include{OMINun}

\include{OMINdeux}
\include{OMINtrois}
\include{OMINquatre}
\include{OMINcinq}
\include{OMINsix}
\include{OMINbib}
\printindex
\end{document}

%% file: OMINintro.tex
\chapter*{Preface}

These notes have served as a basis for a course in Pisa in Spring 1999. 
A parallel course on the construction of o-minimal structures was given
by A.~Macintyre.\par
The content of these notes owes a great deal to the excellent book by
L.~van den Dries \cite{vdD}. Some interesting topics contained in this
book are not included here, such as the Vapnik-Chervonenkis property. 
\par
 Part of the material which does not come from \cite{vdD} is taken from
the paper \cite{CoTT}. This includes the sections on the choice of good
coordinates and the triangulation of functions in Chapter 4 and Chapter 5.
The latter chapter contains the results on triviality in families of sets
or functions which were the main aim of this course.\par
The last chapter on smoothness was intended to establish property ``${\rm
D}{\cal C}^k$ - all $k$'' which played a crucial role in the course of
Macintyre (it can be easily deduced from the results in \cite{vDMi}). It
is also the occasion to give a few results on tubular neighborhoods.\par
 I am pleased to thank Francesca Acquistapace, Fabrizio
Broglia and all colleagues of the Dipartimento di Matematica for the
invitation to give this course in Pisa and their friendly hospitality.
Also many thanks to Antonio Diaz-Cano, Pietro Di Martino, Jesus
Escribano and Federico Ponchio for reading these notes and correcting
mistakes.

%% file: OMINun.tex
\chapter{O-minimal Structures}

\section{Introduction}

The main feature of o-minimal structures is that there are no ``monsters'' in such structures. Let us
take an example of the pathological behaviour that is ruled out. Let $\Gamma\subset {\Bbb
R}^2$ be the graph of the function $x\mapsto \sin(1/x)$ for $x>0$, and let $\clos (\Gamma)$ be
the closure of $\Gamma$ in ${\Bbb R}^2$. This set $\clos(\Gamma)$ is the union of $\Gamma$ and
the closed segment joining the two points $(0,-1)$ and $(0,1)$. We have 
$$\dim(\clos( \Gamma) \setminus \Gamma)= \dim \Gamma= 1.$$
Observe also that $\clos(\Gamma)$ is connected, but not arcwise connected: there is no
continuous path inside $\clos(\Gamma)$ joining the origin with a point of $\Gamma$. The
o-minimal structures will allow to develop a ``tame topology'' in which such bad things cannot
happen. The model for o-minimal structures is the class of semialgebraic sets.\par
A {\em semialgebraic subset}\index{semialgebraic!-- set} of ${\Bbb R}^n$ is a subset defined by a
boolean combination of polynomial equations and inequalities with real coefficients. In other
words, the semialgebraic subsets of ${\Bbb R}^n$ form the smallest class
${\cal SA}_n$ of subsets of ${\Bbb R}^n$ such that:
\begin{enumerate}
\item If  $P \in {\Bbb R}[X_1,\ldots, X_n]$, then  $\{x \in {\Bbb R}^n\ ;\ P(x) = 0\} \in {\cal
SA}_n$  and  $\{x \in {\Bbb R}^n\ ;\ P(x) > 0\} \in {\cal SA}_n$.
\item If $A \in {\cal SA}_n$ and $B \in {\cal SA}_n$, then $A \cup B$, $A \cap B$ and ${\Bbb R}^n
\setminus A$ are in ${\cal SA}_n$.
\end{enumerate}
On the one hand the class of semialgebraic sets
 is stable under many constructions (such as taking projections, closure, connected
components\dots), and on the other hand the topology of semialgebraic sets is very simple, without
pathological behaviour. The o-minimal structures may be seen as an axiomatic treatment of
semialgebraic geometry. An o-minimal structure (expanding the field of reals) is the data, for
every positive integer $n$, of a subset
${\cal S}_n$ of the powerset of ${\Bbb R}^n$, satisfying certain axioms.  There are axioms which
allow to perform many constructions  inside the structure, and an ``o-minimality axiom'' which
guarantee the tameness of the topology.\par
One can distinguish two kinds of activities in the study of o-minimal structures. The first one is to
develop the geometry of o-minimal structures from the axioms. In this activity one tries to follow
the semialgebraic model. The second activity is to discover new interesting classes statisfying the
axioms of o-minimal structures. This activity is more innovative, and the progress made in this
direction (starting with the proof by Wilkie that the field of reals with the exponential function
defines an o-minimal structure) justifies the study of o-minimal structures. The subject of this
course is the first activity (geometry of o-minimal structures), while the course of
A.~Macintyre is concerned with the construction of o-minimal structures.

\section{Semialgebraic Sets}

The semialgebraic subsets of ${\Bbb R}^n$ were defined above. We denote by ${\cal SA}_n$ the set of
all semialgebraic subsets of ${\Bbb R}^n$\index{semialgebraic!-- set}. Some stability properties
of the class of semialgebraic sets follow immediately from the definition.
\begin{enumerate}
\item All algebraic subsets of ${\Bbb R}^n$ are in ${\cal SA}_n$. Recall that an algebraic subset is
a subset defined by a finite number of polynomial equations
$$P_1(x_1,\ldots,x_n)=\ldots=P_k(x_1,\ldots,x_n)=0\;.$$
\item ${\cal SA}_n$ is stable under the boolean operations,
i.e. finite unions and intersections and taking complement. In other words, ${\cal SA}_n$ is a
Boolean subalgebra of the powerset ${\cal P}({\Bbb R}^n)$.
\item The cartesian product of semialgebraic sets is semialgebraic. If $A\in{\cal SA}_n$ and
$B\in{\cal SA}_p$, then $A\times B\in{\cal SA}_{n+p}$.
\end{enumerate}

Sets are not sufficient, we need also maps. Let $A\subset {\Bbb R}^n$ be a semialgebraic set. A map
$f: A\to {\Bbb R}^p$ is said to be semialgebraic\index{semialgebraic!-- map} if its graph
$\Gamma(f)\subset {\Bbb R}^n\times{\Bbb R}^p={\Bbb R}^{n+p}$ is semialgebraic. For instance, the
polynomial maps are semialgebraic. The function $x\mapsto \sqrt{1-x^2}$ for $|x|\leq 1$ is
semialgebraic.\par

The most important stability property of semialgebraic sets is known as ``Tarski-Seidenberg
theorem''\index{Tarski-Seidenberg theorem}. This central result in semialgebraic geometry is not
obvious from the definition.

\begin{enumerate}
\setcounter{enumi}{3}
\item Denote by $p: {\Bbb R}^{n+1}\to {\Bbb R}^n$ the projection on the first $n$ coordinates. Let
$A$ be a semialgebraic subset of ${\Bbb R}^{n+1}$. Then the projection $p(A)$ is semialgebraic.
\end{enumerate}

The Tarski-Seidenberg theorem has many consequences. For instance, it implies that the composition
of two semialgebraic maps is semialgebraic. We shall say more on the consequences of the
stability under projection in the context of o-minimal structures.\par

The semialgebraic subsets of the line are very simple to describe. Indeed, a semialgebraic subset
$A$ of ${\Bbb R}$ is described by a boolean combination of sign conditions ($<0$, $=0$ or $>0$) on
polynomials in one variable. Consider the finitely many roots of all polynomials appearing in the
definition of $A$. The signs of the polynomials are
constant on the intervals delimited by these roots. Hence, such an interval is either disjoint
from or contained in $A$. We obtain the following description.

\begin{enumerate}
\setcounter{enumi}{4}
\item The elements of ${\cal SA}_1$ are the finite unions of points
and open intervals.
\end{enumerate}

We cannot hope for such a simple description of semialgebraic subsets of ${\Bbb R}^n$, $n>1$.
However, we have that every semialgebraic set has a finite partition into semialgebraic subsets
homeomorphic to open boxes (i.e.\ cartesian product of open intervals). This is a consequence of the
so-called ``cylindrical algebraic decomposition''\index{cylindrical algebraic decomposition}
(cad), which is the main tool in the study of semialgebraic sets. Actually, the Tarki-Seidenberg
theorem can be proved by using cad. A
 cad of ${\Bbb R}^n$ is a partition of ${\Bbb R}^n$ into finitely many semialgebraic subsets (which
are called the cells\index{cell} of the cad), satisfying certain properties. We define the cad of
${\Bbb R}^n$ by induction on $n$.
\begin{itemize}
\item A cad of ${\Bbb R}$ is a subdivision by finitely many points $a_1<\ldots<a_\ell$. The cells are
the singletons $\{a_i\}$ and the open intervals delimited by these points.
\item For $n>1$, a cad of ${\Bbb R}^n$ is given by a cad of ${\Bbb R}^{n-1}$ and, for each cell $C$
of ${\Bbb R}^{n-1}$, continuous semialgebraic functions
$$\zeta_{C,1}<\ldots<\zeta_{C,\ell_C}: C\to {\Bbb R}\;.$$
The cells of ${\Bbb R}^n$ are the graphs of the functions $\zeta_{C,i}$ and the bands in the cylinder
$C\times{\Bbb R}\subset{\Bbb R}^n$ delimited by these graphs. For $i=0,\ldots,\ell_C$, the
band
$(\zeta_{C,i},\zeta_{C,i+1})$ is
$$(\zeta_{C,i},\zeta_{C,i+1})=\{ (x',x_n)\in {\Bbb R}^n\;;\ x'\in C \et
\zeta_{C,i}(x')<x_n<\zeta_{C,i+1}(x')\}\;,$$
where we set $\zeta_{C,0}=-\infty$ and $\zeta_{C,\ell_C+1}=+\infty$.
\end{itemize}

Observe that every cell of a cad is homeomorphic to an open box. This is proved by induction on
$n$, since a graph of $\zeta_{C,i}$ is homeomorphic to
$C$ and a band $(\zeta_{C,i},\zeta_{C,i+1})$ is homeomorphic to $C\times(0,1)$. \par

Given a finite list $(P_1,\ldots,P_k)$ of polynomials in ${\Bbb R}[X_1,\ldots,X_n]$, a subset $A$
of ${\Bbb R}^n$ is said to be $(P_1,\ldots,P_k)$-invariant if the sign ($<0$, $=0$ or $>0$)  of each
$P_i$ is constant on $A$. The main result concerning cad is the following.

\begin{theorem}\label{cad} 
Given a finite list $(P_1,\ldots,P_k)$ of polynomials in ${\Bbb
R}[X_1,\ldots,X_n]$, there is a cad of ${\Bbb R}^n$ such that each cell is 
$(P_1,\ldots,P_k)$-invariant.
\end{theorem}

It follows from Theorem \ref{cad} that, for every semialgebraic subset $A$ of ${\Bbb R}^n$, there
is a cad such that $A$ is a union of cells (such a cad is called adapted to $A$). Indeed, $A$ is
defined by a boolean combination of sign conditions on a finite  list of polynomials
$(P_1,\ldots,P_k)$, and it suffices to take a cad such that each cell is
$(P_1,\ldots,P_k)$-invariant.\par

\begin{figure}[ht]
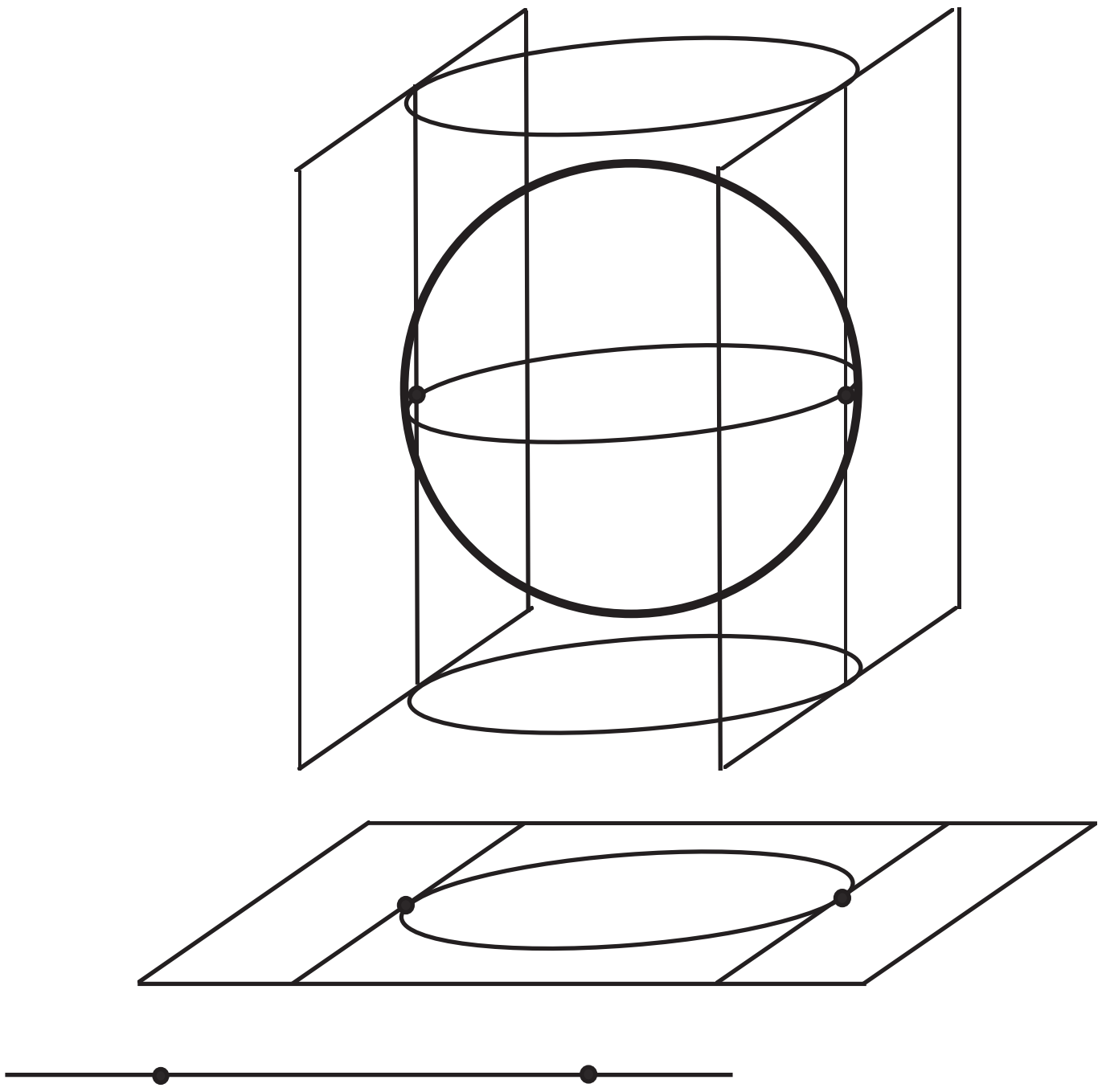

\hspace*{20mm}\dessin 120mm by 90mm (sphereeps)
\caption{A cad adapted to the sphere\label{sphere}}
\end{figure}

We illustrate Theorem \ref{cad} with the example of a cad of ${\Bbb R}^3$ such that each cell is
$(X_1^2+X_2^2+X_3^2-1)$-invariant (a cad adapted to the unit sphere). Such a cad is shown on Figure
\ref{sphere}.

\debexo How many cells of ${\Bbb R}^3$ are there in this cad? Is it possible to have a cad of
${\Bbb R}^3$, adapted to the unit sphere,  with less cells?
\finexo

We refer to \cite{Co} for an introduction to semialgebraic geometry.

Semialgebraic geometry can also be developed over an arbitrary real closed field, instead of the
field of reals. A real closed field $R$ is an ordered field satisfying one of the equivalent
conditions:
\begin{itemize}
\item Every positive element is a square and every polynomial in
$R[X]$ with odd degree has a root in $R$.
\item For every polynomial $F\in R[X]$ and all $a,b$ in $R$ such that $a<b$ and $F(a)F(b)<0$, there
exists $c\in R$, $a<c<b$, such that
$F(c)=0$.
\item $R[\sqrt{-1}]=R[X]/(X^2+1)$ is an algebraically closed field.
\end{itemize}
A semialgebraic subset of $R^n$ is defined as for ${\Bbb R}^n$. The properties of semialgebraic
sets that we have recalled in this section also hold for semialgebraic subsets of $R^n$. We
refer to \cite{BCR} for the study of semialgebraic geometry over an arbitrary real closed field.
 
\section{Definition of an O-minimal Structure}

We shall now define the o-minimal structures expanding a real closed fied $R$. The fact that we
consider a situation more general than the field of reals will be important only in Chapter 5.
Otherwise, the reader may take $R={\Bbb R}$.\par
 An {\em interval in $R$}\index{interval} will always be an open
interval $(a,b)$ (for $a<b$) or
$(a,+\infty)$ or $(-\infty,b)$. We insist that an interval always has endpoints in
$R\cup\{-\infty,+\infty\}$. For instance, if $R$ is the field of real algebraic numbers (which is
the smallest real closed field), the set of $x$ in $R$ such that
$0<x<\pi$ ($\pi=3.14\ldots$) is not an interval, because there is no right endpoint in $R$.
Occasionally, we shall also use the notation $[a,b]$ for closed segments in $R$. \par   
The field
$R$ has a topology for which the intervals form a basis. Affine spaces $R^n$ are endowed with the
product topology. The {\em open boxes}\index{open box}, i.e. the cartesian products of open
intervals
$(a_1,b_1)\times\cdots\times(a_n,b_n)$ form a basis for this topology. The polynomials are
continuous for this topology.
\debexo Prove that the polynomials are continuous. \finexo
 When dealing with the topology, one
should take into account that a real closed field is generally not locally connected nor locally
compact (think of real algebraic numbers).\par

\begin{definition} A {\em structure expanding the real closed field $R$}\index{structure} is a
collection
${\cal S}=({\cal S}^n)_{n\in{\Bbb N}}$, where each ${\cal S}^n$ is a set of subsets of the
affine space
$R^n$, satisfying the following axioms:
\begin{enumerate}
\item All algebraic subsets of $R^n$ are in $S_n$.
\item For every $n$, $S_n$ is a Boolean subalgebra of the powerset of $R^n$.
\item If $A\in S_m$ and $B\in S_n$, then $A\times B\in S_{m+n}$.
\item If $p:R^{n+1}\to R^n$ is the projection on the first $n$ coordinates and $A\in S_{n+1}$, then
$p(A)\in S_n$.
\end{enumerate} 
The elements of $S_n$ are called the {\em definable subsets of
$R^n$}\index{definable!-- set}. The structure ${\cal S}$ is said to be {\em
o-minimal}\index{structure!o-minimal --} if, moreover, it satisfies:
\begin{enumerate}\setcounter{enumi}{4}
\item The elements of $S_1$ are precisely the finite unions of points and intervals.
\end{enumerate}
\end{definition}

{\em In the following, we shall always work in a o-minimal structure expanding a real closed field
$R$.}

\begin{definition} A map $f: A\to R^p$ (where $A\subset R^n$) is called {\em
definable}\index{definable!-- map} if its graph is a definable subset of $R^n\times R^p$.
(Applying $p$ times property 4, we deduce that
$A$ is definable).
\end{definition}

\begin{proposition} The image of a definable set by a definable map is definable.
\end{proposition}

\begin{proof}
 Let $f:A\to R^p$ be
definable, where $A\subset R^n$, and let $B$ be a definable subset of $A$. Denote by
$\Gamma_f=\{(x,f(x))\;;\ x\in A\}\subset R^{n+p}$ the graph of $f$. Let $\Delta$ be the algebraic
(in fact, linear) subset of $R^{p+n+p}$ consisting of those $(z,x,y)\in R^p\times R^n\times R^p$
such that $z=y$. Then $C=\Delta\cap (R^p\times \Gamma)\cap (R^p\times B\times R^p)$ is a definable
subset of $R^{p+n+p}$. Let $p_{p+n+p,p}:R^{p+n+p}\to R^p$ be the projection on the first $p$
coordinates. We have $p_{p+n+p,p}(C)=f(B)$, and, applying $n+p$ times property 4, we deduce that
$f(B)$ is definable.
\end{proof}

Observe that every polynomial map is a definable map, since its graph is an algebraic set.

\debexo  Every semialgebraic subset of $R^n$ is definable (Hint: the set defined by
$P(x_1,\ldots,x_n)>0$ is the projection of the algebraic set with equation
$x_{n+1}^2 P(x_1,\ldots,x_n)-1=0$). Hence, the collection of ${\cal SA}_n$ is the smallest
o-minimal structure expanding $R$.
\finexo 

\debexo \label{liminf} Show that every nonempty definable subset of
$R$ has a least upper bound in $R\cup\{+\infty\}$.
\finexo 

\debexo  Assume that ${\cal S}$ is an o-minimal structure expanding an ordered field $R$ (same
definition as above). Show that
$R$ is real closed. (Hint: one can use the second equivalent condition for a field to be real closed,
the continuity of polynomials and property 5 of o-minimal structures.)
\finexo 

\debexo Let $f=(f_1,\ldots,f_p)$ be a map from $A\subset R^n$ into $R^p$. Show that $f$ is
definable if and only if each of its coordinate functions $f_1,\ldots,f_p$ is definable
\finexo

\debexo Show that the composition of two definable maps is definable. Show that the definable
functions $A\to R$ form an $R$-algebra.
\finexo

\begin{proposition} The closure and the interior of a definable subset of $R^n$ are definable.
\end{proposition}

\begin{proof}
It is sufficient to prove the assertion concerning the closure. The case of the interior follows
by taking complement. Let
$A$ be a definable subset of $R^n$. The closure of
$A$ is
$$\displaylines{
\quad\clos(A) =\hspace*{\fill}\cr
\hspace*{\fill} \left\{ x
\in R^n\ ;\ \forall \varepsilon \in R,\;  \varepsilon > 0 \Rightarrow
\exists y \in R^n,\;  y \in A  \et  \sum_{i=1}^n (x_i-y_i)^2 < \varepsilon^2
\right\}\;.\quad\cr}$$
where $x=(x_1,\ldots,x_n)$ and $y=(y_1,\ldots,y_n)$.  
The closure of $A$ can also be described as
$$\clos(A) = R^n \setminus \left(p_{n+1,n}\left(R^{n+1} \setminus
p_{2n+1,n+1}(B)\right)\right)\;,$$
where
$$B =\left(R^n \times R \times A\right)\cap\left\{ (x, \varepsilon, y) \in R^n
\times R \times
R^n\ ;\  \sum_{i=1}^n (x_i-y_i)^2 < \varepsilon^2\right\}\;,$$ 
 $p_{n+1,n}(x, \varepsilon) = x$ and
 		$p_{2n+1,n+1}(x, \varepsilon, y) =  (x, \varepsilon)$. Then observe that $B$ is definable.
\end{proof}

The example above shows that it is usually boring to write down projections in order to show that a
subset is definable. We are more used to write down formulas. Let us make precise what is meant
by a {\em first-order formula}\index{first-order formula} (of the language of the o-minimal
structure). A first-order formula is constructed according to the following rules.
\begin{enumerate}\shitem
\item If $P \in R[X_1 ,\ldots, X_n]$, then $P(x_1,\ldots,x_n) = 0$ and $P(x_1,\ldots,x_n) > 0$ are
first-order formulas.
\item If $A$ is a definable subset of $R^n$, then $x\in A$ (where $x=(x_1,\ldots,x_n))$ is a
first-order formula.
\item If $\Phi(x_1,\ldots,x_n)$ and $\Psi(x_1,\ldots,x_n)$ are first-order formulas, then ``$\Phi \et
\Psi$'', ``$\Phi
\ou
\Psi$'', ``$\mbox{\rm not }\Phi$'' , $\Phi\Rightarrow\Psi$ are first order formulas.
\item If $\Phi(y,x)$ is a first-order formula (where $y=(y_1,\ldots,y_p)$ and
$x=\linebreak (x_1,\ldots,x_p)$) and
$A$ is a definable subset of $R^n$, then
$\exists x\in A\
\Phi(y,x)$ and
$\forall x\in A\ \Phi(y,x)$ are first-order formulas.
\end{enumerate}

\begin{theorem} If $\Phi(x_1,\ldots,x_n)$ is a first-order formula,
the set of 
$(x_1 ,\ldots, x_n)$ in $R^n$ which satisfy $\Phi(x_1 ,\ldots, x_n)$ is definable.
\end{theorem}

\begin{proof}
By induction on the construction of formulas. Rule 1 produces semialgebraic sets, which are
definable. Rule 2 obviously produces definable sets. Rule 3 works because ${\cal SA}_n$ is closed
under boolean operation. Rule 4 reflects the fact that the projection of a definable set is
definable. Indeed, if $B=\{(y,x)\in R^{p+n}\;;\ \Phi(y,x)\}$ is definable and
$p_{p+n,p}:R^{p+n}\to R^p$ denotes the projection on the first $p$ coordinates, we have
\begin{eqnarray*}
\{y\in R^p\;;\ \exists x\in A\ \Phi(y,x)\} &=&  p_{p+n,p}\big((R^p\times A)\cap B\big)\;,\\
\{y\in R^p\;;\ \forall x\in A\ \Phi(y,x)\} &=& R^p\setminus p_{p+n,p}\big((R^p\times
A)\cap(R^{p+n}\setminus B)\big)\;,
\end{eqnarray*}
which shows that both sets are definable. 
\end{proof}

 One should pay attention to the fact that the quantified variables  have to
range over definable sets. For instance,
$$\{ (x,y) \in {\Bbb R}^2\ ;\ \exists n \in {\Bbb N}\  
y=nx \}$$
is not definable (why ?).

We refer the reader to \cite{Pr} for notions of model theory.

\debexo Let $f:A\to R$ be a definable function. Show that $\{ x\in A\;;\ f(x)>0\}$ is definable.
Hence, we can accept inequalities involving definable functions in formulas defining definable
sets.
\finexo

\debexo Let $A$ be a non empty definable subset of $R^n$. For $x\in R^n$, define ${\rm
dist}(x,A)$ as the greatest lower bound of the set of $\Vert
y-x\Vert=\sqrt{\sum_{i=1}^n (y_i-x_i)^2}$ for $y\in A$. Show that ${\rm dist}(x,A)$ is
well-defined and that $x\mapsto{\rm dist}(x,A)$ is a continuous definable function on $R^n$.
\finexo

\debexo
Let $f:A\to R$ be a definable function and
$a\in\clos(A)$. For $\varepsilon>0$, define
$$m(\varepsilon)= \inf\{ f(x)\;;\ x\in A\et \Vert x-a\Vert<\varepsilon\}\in R\cup\{-\infty\}\;.$$
Show that $m$ is a definable function (this means that $m^{-1}(-\infty)$ is definable and
$m|_{m^{-1}(R)}$ is definable).  Show that
$\liminf_{x\to a} f(x)$ exists in
$R\cup\{-\infty,+\infty\}$.
\finexo

%% file: OMINdeux.tex
\chapter{Cell Decomposition}

In this chapter we prove the cell decomposition for definable sets,
which generalizes the cylindrical algebraic decomposition of
semialgebraic sets. This result is the most important for the study
of o-minimal geometry. The proofs are rather technical. The main
results come from \cite{PS,KPS}, and we follow \cite{vdD} rather
closely.

\section{Monotonicity Theorem}

\begin{theorem}[Monotonicity Theorem]\label{monotonicity}
Let $f:(a,b)\to R$ be a definable function. There exists a finite
subdivision $a=a_0<a_1<\ldots<a_k=b$ such that, on each interval
$(a_i,a_{i+1})$, $f$ is continuous and either constant or strictly
monotone.
\end{theorem}

The key of the proof of the Monotonicity Theorem is the following
Lemma.

\begin{lemma}\label{subinterval} Let $f:(a,b)\to R$ be a definable
function. There exists a subinterval of $(a,b)$ on which $f$ is
constant, or there exists a subinterval on which $f$ is strictly
monotone and continuous.
\end{lemma}

\begin{proof}
Suppose that there is no subinterval of $(a,b)$ on which $f$ is
constant.\par\smallskip
{\em First step: there exists a subinterval on which $f$ is
injective.} It follows from the assumption that, for all
$y$ in
$R$, the definable set
$f^{-1}(y)$ is finite: otherwise, $f^{-1}(y)$ would contain an
interval, on which
$f=y$. Hence the definable set $f((a,b))$ is infinite and contains an
interval $J$. The function $g:J\to (a,b)$ defined by
$g(y)=\min(f^{-1}(y))$ is definable and satisfies $f\circ g={\rm
Id}_J$. Since $g$ is injective, $g(J)$ is infinite and contains a
subinterval $I$ of $(a,b)$. We have $g\circ f|_I={\rm Id}_I$, and $f$
is injective on $I$.\par\smallskip
{\em Second step: there exists a subinterval on which $f$ is
strictly monotone.} We know that $f$ is injective on a subinterval
$I$ of $(a,b)$. Take $x\in I$. The sets
\begin{eqnarray*}
I_+&= &\{ y\in I\;;\ f(y) > f(x)\}\\
I_-&= &\{ y\in I\;;\ f(y) < f(x)\}
\end{eqnarray*}
form a definable partition of $I\setminus\{x\}$. Therefore, there is
$\varepsilon>0$ such that $(x-\varepsilon,x)$ (resp.\
$(x,x+\varepsilon)$) is contained in $I_+$ or in $I_-$. There are
four possibilities $\Phi_{+,+}(x)$, $\Phi_{+,-}(x)$, $\Phi_{-,+}(x)$
and $\Phi_{-,-}(x)$ which give a definable partition of $I$. For
instance, $\Phi_{+,-}(x)$ is
$$\displaylines{
\qquad \exists\varepsilon\ \ \forall y \in I\ \
\Big( (x-\varepsilon<y<x \Rightarrow f(y)>f(x)) \et\hspace*{\fill}\cr
\hspace*{\fill}(x<y<x+\varepsilon\Rightarrow
f(y)<f(x))\Big)\;.\qquad\cr}$$
\par We claim that $\Phi_{+,+}$ and
$\Phi_{-,-}$ are finite. It is sufficient to prove the claim for
$\Phi_{+,+}$ (for
$\Phi_{-,-}$, replace $f$ with $-f$). If $\Phi_{+,+}$ is not finite,
it holds on a subinterval of $I$, which we still call $I$. Set
$$B=\{x\in I\;;\ \forall y \in I\quad y>x\Rightarrow f(y)> f(x)\}\;.$$
If $B$ contains an interval, then $f$ is strictly increasing on this
interval, which contradicts $\Phi_{+,+}$. Hence, $B$ is finite.
Replacing $I$ with $(\max(B),+\infty)\cap I$ if $B\not=\emptyset$, we
can assume
$$ \forall x\in I\ \exists y \in I\ (y>x \et f(y)<f(x))\;.
\leqno{(*)}$$
Take $x\in I$. The definable set $C_x=\{y\in I\;;\ y>x \et
f(y)<f(x)\}$ is nonempty. If $C_x$ were finite, its maximum would be
an element of $I$ contradicting the property (*). Therefore, $C_x$
contains an interval. Let $z$ be the greatest lower bound of the
interior of
$C_x$. We have $z>x$ because $f>f(x)$ on some interval
$(x,x+\varepsilon)$. We have also $f>f(x)$ on some interval
$(z-\varepsilon,z)$ and
$f<f(x)$ on some interval $(z,z+\varepsilon)$. Hence, the following
formula
$\Psi_{+,-}(z)$ holds:
$$\exists\varepsilon>0\ \forall t\in I\ \forall u\in I \ 
(z-\varepsilon<t<z<u<z+\varepsilon \Rightarrow f(t)>f(u))\;.$$
We have shown
$$\forall x\in I\ \exists z\in I\ (x<z \et \Psi_{+,-}(z))\;.
\leqno{(**)}$$
The definable set of elements of $I$ satisfying $\Psi_{+,-}$ is not
empty and it is not finite: otherwise, its maximum would contradict
(**). It follows that, replacing $I$ with a smaller subinterval, we can
assume that $\Psi_{+,-}$ holds on $I$. Consider the definable
function $h: -I\to R$ defined by $h(x)=f(-x)$. The property
$\Phi_{+,+}$ for $h$ holds on $-I$. Hence, by the preceding argument,
there is a subinterval of $-I$ on which $\Psi_{+,-}$ for $h$ holds.
This means that  there is a subinterval of $I$ on which $\Psi_{-,+}$
for $f$ holds (exchange left and right). On this subinterval we have
both $\Psi_{+,-}$ and $\Psi_{-,+}$, which are contradictory. This
contradiction proves the claim.\par
Since $\Phi_{+,+}$ and $\Phi_{-,-}$ are finite, replacing $I$ with a
smaller subinterval, we can assume that $\Phi_{-,+}$ holds
on $I$ or $\Phi_{+,-}$ holds on $I$. Say $\Phi_{-,+}$ holds on $I$.
Then $f$ is strictly increasing on $I$. Indeed, for all $x\in I$, the
definable set $\{y\in I\;;\ y>x\ \mbox{and}\ f>f(x)\ \mbox{on}\ 
(x,y)\}$ is nonempty and its least upper bound is necessarily the
right endpoint of $I$ (Otherwise, this l.u.b.\ would not satisfy
$\Phi_{-,+}$). Similarly, if 
$\Phi_{+,-}$ holds on
$I$,
$f$ is strictly decreasing on $I$.\par\smallskip
{\em Last step: there is a subinterval on which $f$ is strictly
monotone and continuous}. Recall that $f$ is strictly monotone on $I$.
The definable set $f(I)$ is infinite and contains an interval $J$.
The inverse image $f^{-1}(J)$ is the interval $(\inf(f^{-1}(J)),
\sup(f^{-1}(J)))$. Replacing $I$ with this interval, we can assume that
$f$ is a strictly monotone bijection from the interval $I$ onto the
interval $J$. Since the inverse image of a subinterval of $J$ is a
subinterval of $I$, $f$ is continuous.
\end{proof}

{\par\noindent{\it Proof of the Monotonicity Theorem.
}\nopagebreak\normalsize}
Let $X_=$ (resp.\ $X_\nearrow$, resp. $X_\searrow$) be the
definable set of $x\in(a,b)$ such that $f$ is constant
(resp.\ continuous and strictly increasing, resp.\ continuous and
strictly decreasing) on an interval containing $x$. Then 
$$(a,b)\setminus (X_= \cup X_\nearrow \cup X_\searrow)$$
is finite: otherwise, it would contain an interval, and we get a
contradiction by applying Lemma \ref{subinterval} to this interval.
Hence, there is a subdivision $a=a_0<a_1<\ldots<a_k=b$ such that each
$(a_i,a_{i+1})$ is contained in $X_=$ or in $X_\nearrow$ or in
$X_\searrow$. Clearly, $f$ is continuous on each $(a_i,a_{i+1})$.
Take $x$ in $(a_i,a_{i+1})$. If $(a_i,a_{i+1})$ is contained in 
$X_=$ (resp.\ $X_\nearrow$, resp.\ $X_\searrow$), let $D_x$ be the
set of $y\in (a_i,a_{i+1})$ such that $x<y$ and $f=f(x)$ (resp.\
$f>f(x)$, resp.\ $f<f(x)$) on $(x,y)$. The definable set $D_x$ is
nonempty and its least upper bound is necessarily $a_{i+1}$. It
follows that $f$ is constant or strictly monotone on $(a_i,a_{i+1})$.
    {\linebreak[2]\hspace*{\fill}$\square$\ifdim\lastskip<12pt
    \removelastskip \penalty-200 \vskip12pt \fi}

\debexo \label{leftrightlimit} Let $f:(a,b)\to R$ be a
definable function. Then
$\lim_{x\to b_-} f(x)$ and $\lim_{x\to a_+} f(x)$ exist in
$R\cup\{-\infty,+\infty\}$.
\finexo 

\section[Cell Decomposition]{Uniform Finiteness, Cell Decomposition
and Piecewise Continuity}\label{sectioncdcd}

The notion of a {\em cylindrical definable cell decomposition (cdcd) of
$R^n$}\index{cylindrical definable cell decomposition (cdcd)}  is a
generalization of cylindrical algebraic cell decomposition. We
define the cdcd by induction on
$n$

\begin{definition}\label{defcdcd} A cdcd of $R^n$ is a finite partition
of
$R^n$ into definable sets $(C_i)_i\in I$ satisfying certain properties
explained below. The $C_i$ are called the {\em cells}\index{cell} of
the cdcd.
\begin{description}
\item[$n=1$:] A cdcd of
$R$ is given by a finite subdivision $a_1<\ldots<a_\ell$ of $R$. The
cells of $R$ are the singletons
$\{a_i\}$,
$0<i\leq\ell$, and the intervals $(a_i,a_{i+1})$, $0\leq i\leq\ell$,
where $a_0=-\infty$ and
$a_{\ell+1}=+\infty$.
\item[$n>1$:] A cdcd of $R^n$ is given by a cdcd of $R^{n-1}$ and,
for each cell $D$ of $R^{n-1}$, continuous definable functions
$$\zeta_{D,1}<\ldots<\zeta_{D,\ell(D)}:D\longrightarrow R\;.$$ 
The cells of $R^n$ are the {\em graphs}\index{graph}
$$\{ (x,\zeta_{D,i}(x))\;;\ x\in D\}\;,\qquad
0<i\leq\ell(D)\;,$$
and the {\em bands}\index{band}
$$(\zeta_{D,i},\zeta_{D,i+1})=\{(x,y)\;;\ x\in D \et
\zeta_{D,i}(x)<y<\zeta_{D,i+1}(x)\}$$ for $0\leq i\leq\ell(D)$,
where $\zeta_{D,0}=-\infty$ and $\zeta_{D,\ell(D)+1}=+\infty$.
\end{description}
\end{definition}

Note that the fact that a cell of a cdcd is definable follows
immediately from the definition and the axioms of o-minimal structure.
Note also that if $p_{n,m}:R^n\to R^m$, $m<n$, is the projection on the
first $m$ coordinates, the images by $p$ of the cells of a cdcd of
$R^n$ are the cells of a cdcd of $R^m$.\par
We define by induction the {\em dimension of a
cell}\index{dimension!-- of a cell}. The dimension of a singleton is
$0$ and the dimension of an interval is
$1$. If $C$ is a cell of $R^n$, its dimension is $\dim(p_{n,n-1}(C))$
if $C$ is a graph and $\dim(p_{n,n-1}(C))+1$ if $C$ is  a band.

\begin{proposition}\label{cell=affine} For each
cell
$C$ of a cdcd of $R^n$, there is a definable homeomorphism 
$\theta_C: C\to R^{\dim(C)}$. 
\end{proposition}

\begin{proof} Let $D=p_{n,n-1}(C)$ and assume $\theta_D: D\to
R^{\dim(D)}$ is already defined. Let $(x,y)\in C$, where $x\in D$. We
define $\theta_C(x,y)$ as $\theta_D(x)$ if $C$ is a graph, and as
$$\left( \theta_D(x)\;,\ {1\over \zeta_{D,i}(x)-y} + y + {1\over
\zeta_{D,i+1}(x)-y}\right)$$
if $C$ is the band $(\zeta_{D,i},\zeta_{D,i+1})$ (fractions where
infinite functions appear in the denominator have to be omitted).
\end{proof}

\debexo\label{opencells} Prove  (by induction on $n$) that a cell is
open in
$R^n$ if and only if
 its dimension is $n$. Prove that the union of all cells of dimension
$n$ is dense in $R^n$ (hint: show by induction on $n$ that every open
box in $R^n$ meets a cell of dimension $n$).
\finexo

We have already said that the notion of connectedness does not behave
well if $R\not={\Bbb R}$. It has to be replaced by the notion of
definable connectedness.

\begin{definition}\label{defconnected} A definable set $A$ is said to
be {\em definably connected}\index{connected!definably --} if, for all
disjoint definable open subsets
$U$ and $V$ of
$A$ such that $A=U\cup V$, one has $A=U$ or $A=V$. A definable set $A$
is said to be {\em definably arcwise
connected}\index{connected!definably arcwise --} if, for all points
$a$ and $b$ in $A$, there is a definable continuous map
$\gamma:[0,1]\to A$ such that $\gamma(0)=a$ and $\gamma(1)=b$.
\end{definition}

\debexo \label{exoconnected} Prove the following facts: The segment
$[0,1]$ is definably connected. Definably arcwise connected implies
definably connected. Every box $(0,1)^d$ is definably connected. Every
cell of a cdcd is definably connected.
\finexo 

We now state the three main results of this section. In the
following, we denote by $\sharp S$ the number of elements of a finite
set $S$.

\begin{theorem}[Uniform Finiteness ${\rm UF}_n$]\label{UF} Let $A$ be
a definable subset of $R^n$ such that, for every $x\in R^{n-1}$, the
set
$$A_{x}=\{ y\in R\;;\ (x,y)\in A\}$$ is finite. Then there exists
$k\in {\bf N}$ such that $\sharp A_x \leq k$ for every $x\in R^{n-1}$.
\end{theorem}

\begin{theorem}[Cell Decomposition ${\rm CDCD}_n$]\label{CDCD} Let
$A_1,\ldots,A_k$ be definable subsets of $R^n$. There is a cdcd of
$R^n$ such that each
$A_i$ is a union of cells.
\end{theorem}

A cdcd of $R^n$ satisfying the property of the theorem will be called
{\em adapted to $A_1,\ldots,A_k$}\index{adapted cdcd}.

\debexo\label{nonemptyint} 
Prove the following consequence of ${\rm CDCD}_n$: Let $X$ be
a definable subset of $R^n$ such that $\clos(X)$ has nonempty
interior.
Then $X$ has nonempty interior. (Hint: use Exercise \ref{opencells}
for a convenient cdcd of $R^n$.)
\finexo

\begin{theorem}[Piecewise Continuity ${\rm PC}_n$]\label{PC} Let $A$
be a definable subset of $R^n$ and $f:A\to R$ a definable function.
There is a cdcd of
$R^n$ adapted to $A$ such that, for every cell $C$ contained in $A$,
$f|_C$ is continuous. 
\end{theorem} 

\debexo\label{boxcontinuity} Prove the following consequence of ${\rm
PC}_n$: Let $X$ be an open definable subset of $R^n$ and $f: X\to R$
a definable function. Then there is an open box $B\subset X$ such
that $f|_B$ is continuous (hint: use Exercise \ref{opencells}).
\finexo

We prove the three theorems simultaneously by induction on $n$. For
$n=1$, ${\rm UF}_1$ is trivial ($R^0$ is just one point), ${\rm
CDCD}_1$ follows immediately from the o-minimality, and ${\rm PC}_1$
is a consequence of the Monotonicity Theorem \ref{monotonicity}. In
the following we assume $n>1$ and ${\rm UF}_m$, ${\rm CDCD}_m$
and ${\rm PC}_m$ hold for all integers $m$ such that
$0<m<n$.\par
\medskip
{\par\noindent{\it Proof of ${\rm UF}_n$.
}\nopagebreak\normalsize} 
We can assume that, for every $x\in R^{n-1}$, the set $A_x$ is
contained in $(-1,1)$. Let $\mu: R\to (-1,1)$ be the semialgebraic
homeomorphism defined by $\mu(y)=y/\sqrt{1+y^2}$. We can replace $A$
with its image by $(x,y)\mapsto (x,\mu(y))$, which satisfies the
assumption.\par
 For
$x\in R^{n-1}$ and
$i=1,2,\ldots$, define
$f_i(x)$ to be the $i$-th element of $A_x$, if it exists. Note that
the function
$f_i$ is definable. Call $a\in R^{n-1}$ {\em good} if
$f_1,\ldots,f_{\sharp(A_a)}$ are defined and continuous on an open box
containing $a$ and $a$ does not belong to the closure of the domain
of $f_{\sharp(A_a)+1}$. In other words, $a$ is good if and only if
there is an open box $B\ni a$ in $R^{n-1}$, such that $(B\times R)\cap
A$ is the union of graphs of continuous definable functions 
$$\zeta_1<\ldots<\zeta_{\sharp(A_a)} : B\to (-1,1)\qquad
(\zeta_i=f_i|B)\;.$$
Call $a\in R^{n-1}$ {\em bad} if it is not good.\par\smallskip

{\it First step: the set of good points is definable.}
 Let $a$
be a point of $R^{n-1}$. If $b$ belongs to $[-1,1]$, we say that {\em
$(a,b)\in R^n$ is normal} if there is an open box
$C=B\times(c,d)\subset R^n$ containing
$(a,b)$ such that $A\cap C$ is either empty or the graph of a
continuous definable function $B\to (c,d)$.\par
 Clearly, if $a$ is
good,
$(a,b)$ is normal for every $b\in [-1,1]$. Now assume that
$a$ is bad. Let $f_\ell$ be the first function $f_i$ such
that $a$ is in the closure of the domain of $f_i$ and there is no
open box containing $a$ on which $f_i$ is defined and continuous. We
set
$\beta(a)=\liminf_{x\to a} f_\ell(x)\in[-1,1]$. We claim that
$(a,\beta(a))$ is not normal. Suppose $(a,\beta(a))$ is normal. There
is an open box
$B\times(c,d)$ containing $(a,\beta(a))$ whose intersection with $A$ is
the graph of a continuous function $g:B\to (c,d)$. We can assume that
$f_\ell(x)>c$ for all $x\in B$ such that $f_\ell(x)$ is defined. If
$\ell>1$ and
$\beta(a)=f_{\ell-1}(a)$, we would deduce $g=f_{\ell-1}|_B$ since $B$
is definably connected. We would have $f_\ell(x)\geq d$ for all $x\in
B$ such that $f_\ell(x)$ is defined, which contradicts $\beta(a)<d$.
Hence, we can assume $\ell=1$ or $f_{\ell-1}<c$ on
$B$. It follows that $g=f_\ell|_B$, which contradicts the definition
of
$\ell$. Therefore the claim is proved.\par
We have shown that $a\in R^{n-1}$ is good if and only if, for all
$b\in [-1,1]$, $(a,b)$ is normal. From this we
deduce easily that the set of good points is definable.\par\smallskip

{\it Second step: the set of good points is dense.}
Otherwise, there is an open box $B\subset R^{n-1}$ contained in the
set of bad points. Consider the definable function $\beta:B\to
[-1,1]$ defined as in Step 1. By ${\rm PC}_{n-1}$, we
can assume that  $\beta$ is
continuous (see Exercise \ref{boxcontinuity}).
For $x\in B$,  we define $\beta_-(x)$ (resp.\ $\beta_+(x)$) to be the
maximum (resp.\ minimum) of the $y\in A_x$ such that  
$y<\beta(x)$ (resp. $y>\beta(x)$), if such $y$ exist. Using ${\rm
PC}_{n-1}$ and shrinking the box
$B$, we can assume that $\beta_-$ (resp.\ $\beta_+$) either is
nowhere defined on $B$ or is continuous on $B$. Then the set of
$(x,y)\in A\cap (B\times R)$ such that $y\not=\beta(x)$ is open
and closed in $A\cap (B\times R)$. Shrinking further the box $B$, we
can assume that the graph of $\beta|_B$ is either disjoint from $A$ or
contained in $A$. The first case contradicts the
definition of $\beta$. In the second case, $(x,\beta(x))$ would be
normal for every $x\in B$, which contradicts what was
proved in the first step. We have proved that the set of good points
is dense.\par\smallskip

{\it Third step.} By ${\rm CDCD}_{n-1}$, there is a cdcd of $R^{n-1}$
adapted to the set of good points. Let $C$  be a cell of dimension
$n-1$ of $R^{n-1}$. Since good points are dense, every $x\in C$ is
good. Take $a\in C$. The set of $x\in C$ such that
$\sharp(A_x)=\sharp(A_a)$ is definable, open and closed in $C$. By
definable connectedness of $C$, it is equal to $C$. If $D$ is a cell
of $R^{n-1}$ of smaller dimension, we can use the definable
homeomorphism of $\theta_D: D\to R^{\dim(D)}$ and the
assumption that
${\rm UF}_{\dim(D)}$ holds to prove that $\sharp(A_x)$ is uniformly
bounded for $x\in D$. Since there are finitely many cells, the proof of
${\rm UF}_n$ is completed.
 {\linebreak[2]\hspace*{\fill}$\square$\ifdim\lastskip<12pt
    \removelastskip \penalty-200 \vskip12pt \fi}

{\par\noindent{\it Proof of ${\rm CDCD}_n$. }\nopagebreak\normalsize}
Let $A$ be the set of $(x,y)\in R^{n-1}\times R$ such that $y$
belongs to the boundary of one of $A_{1,x},\ldots,A_{k,x}$ (the
boundary of $S$ in $R$ is $\clos(S)\setminus {\rm int}(S)$). Clearly
$A$ is definable and satisfies the assumptions of ${\rm UF}_n$. Hence
$\sharp(A_x)$ has a maximum $\ell$ for $x\in R^{n-1}$, and $A$ is the
union of the graphs of functions $f_1,\ldots,f_\ell$ defined at the
beginning of the proof of ${\rm UF}_n$.  We define the type of 
$x$  in $R^{n-1}$ to consist of the following data:
\begin{itemize}
\item $\sharp(A_x)$, 
\item the sets of $j\in \{1,\ldots,k\}$ such that $f_i(x)\in A_{j,x}$
for $i=1,\ldots,\sharp(A_x)$,
\item the sets of $j\in \{1,\ldots,k\}$ such that
$(f_i(x),f_{i+1}(x))\subset A_{j,x}$, for
$i=0,\ldots, \sharp(A_x)$ (where
$f_0(x)=-\infty$ and
$f_{\sharp(A_x)+1}(x)=+\infty$).
\end{itemize}
 Since there are finitely
many possible types and the set of points in $R^{n-1}$ with a given
type is definable, we deduce from
${\rm CDCD}_{n-1}$ that there is a cdcd of $R^{n-1}$ such that two
points in the same cell of
$R^{n-1}$ have the same type. Moreover, using
${\rm PC}_{n-1}$, we can assume that the cdcd of $R^{n-1}$ is such
that, for each cell $C$ and
$i=1,\ldots,\ell$, either $f_i$ is defined nowhere on $C$ or $f_i$ is
defined and continuous on $C$. The cdcd of $R^{n-1}$ we have
obtained, together with the restrictions of the functions $f_i$ to
the cells of this cdcd, define a cdcd of $R^n$ adapted to
$A_1,\ldots,A_k$.
 {\linebreak[2]\hspace*{\fill}$\square$\ifdim\lastskip<12pt
    \removelastskip \penalty-200 \vskip12pt \fi}

{\par\noindent{\it Proof of ${\rm PC}_n$.}\par\smallskip
{\it First step}. Assume that $A$ is
an open box  $B\times (a,b) \subset R^{n-1}\times R$. We claim that
there is an open box $A'\subset A$ such that $f|_{A'}$ is continuous.
\par
For every $x\in B$, let
$\lambda(x)$ be the least upper bound of the set of $y\in (a,b)$ such
that
$f(x,\cdot)$ is continuous and monotone on $(a,y)$. The
Monotonicity Theorem \ref{monotonicity} implies $\lambda(x) >a$ for
all $x\in B$. The function $\lambda$ is definable. Applying
${\rm PC}_{n-1}$ to $\lambda$ and replacing $B$ with a smaller
open subbox, we can assume that $\lambda$ is continuous (see Exercise
\ref{boxcontinuity}). Replacing again $B$ with a smaller subbox, we can
assume that there is
$c>a$ such that $\lambda>c$. Replacing $b$ with $c$, we can
assume that, for every $x\in B$, $f(x,\cdot)$ is continuous and
monotone on $(a,b)$.
\par
Now consider the set $C$ of points $(x,y)\in B\times (a,b)$ such that
$f(\cdot,y)$ is continuous at $x$. The set $C$ is definable. It
follows from ${\rm PC}_{n-1}$ that, for every $y\in(a,b)$, the set of
$x$ such that $f(\cdot,y)$ is continuous at $x$ is dense in $B$.
Hence, $C$ is dense in $A$. Applying ${\rm CDCD}_n$, we deduce that $C$
contains an open subbox of $A$. Replacing $A$ with this smaller
subbox, we can assume that, for every $y\in (a,b)$, $f(\cdot,y)$ is
continuous.\par
So it suffices to consider the case where $f(x,\cdot)$ is continuous
and monotone on $(a,b)$ for every $x\in B$ and $f(\cdot,y)$ is
continuous on $B$ for every $y\in (a,b)$. In this situation, $f$ is
continuous on $B\times (a,b)$. Indeed, take $(x^0,y^0)\in B\times
(a,b)$ and $I$ an interval containing $f(x^0,y^0)$. By continuity of
$f(x^0,\cdot)$, we find $y^1<y^0<y^2$ such that $f(x^0,y^i)\in I$ for
$i=1,2$. By continuity of $f(\cdot,y^i)$, we find an open box $B'\ni
x^0$ in $B$ such that $f(B'\times\{y^i\})\subset I$ for $i=1,2$. It
follows from the monotonicity of $f(x,\cdot)$ that $f(B'\times
(y^1,y^2))$ is contained in $I$. This proves the continuity of $f$ and
completes the proof of the claim.\par\smallskip
{\it Second step.} Now we can prove ${\rm PC}_n$. Consider the set
$D$ of points of $A$ where $f$ is continuous. The set $D$ is
definable. By
${\rm CDCD}_n$, there is a cdcd of $R^n$ adapted to $A$ and $D$. If
$E$ is an open cell contained in $A$, the first step shows that
$E\cap D$ is nonempty, therefore $E\subset D$ and $f$ is continuous
on $E$. If $F$ is a cell of dimension $d<n$ contained in $A$, there
is a definable homeomorphism $\theta_F: F\to R^d$. Composing $f$
with $\theta_F^{-1}$ and applying ${\rm PC}_d$, we obtain a finite
partition of $F$ into definable subsets $F_i$ such that $f|_{F_i}$ is
continuous. Hence, there is a finite partition of $A$ into definable
subsets $A_1,\ldots,A_k$ such that $f|_{A_i}$ is continuous for
$i=1,\ldots,k$. Using ${\rm CDCD}_n$, we obtain a cdcd of $R^n$
adapted to $A_1,\ldots,A_k$. So we can assume that the $A_i$ are
cells of a cdcd of $R^n$.
\nopagebreak\normalsize}{\linebreak[2]\hspace*{\fill}$\square$\ifdim\lastskip<12pt
    \removelastskip \penalty-200 \vskip12pt \fi}

%% file: OMINtrois.tex
\chapter{Connected Components and Dimension}

\section{Curve Selection}

We begin with a useful result. It says that, if a formula ``$\forall
x\in X\ \exists y\in Y\ (x,y)\in Z$'' holds, where $X$, $Y$ and $Z$
are definable sets, then $y$ can be choosed as a {\em definable}
function of $x\in X$.

\begin{theorem}[Definable Choice]\label{defchoice}\index{definable!-- choice} Let $A$ be a
definable subset of
$R^m\times R^n$. Denote by $p : R^m\times R^n \to R^m$ the projection
on the first $m$ coordinates. There is a definable function $f:p(A)\to R^n$
such that, for every $x\in p(A)$, $(x,f(x))$ belongs to $A$.
\end{theorem}

\begin{proof} It is sufficient to consider the case $n=1$. The
general case follows by decomposing the projection $R^{m+n}\to R^m$
as $R^{m+n}\to R^{m+n-1}\to\ldots\to R^{m+1}\to R^m$.\par
Make a cdcd of $R^{m+1}$ adapted to $A$. The projection $p(A)$ is the
union of the images by $p$ of cells contained in $A$. Hence, we can
assume that $A$ is a cell of $R^{m+1}$, and consequently $p(A)$ is a
cell of $R^m$. If $A$ is the graph of $\zeta_i: p(A)\to R$, we take
$f=\zeta_i$. If $A$ is a band  $(\zeta_i,\zeta_{i+1})$, where both
functions are finite, we take
$f={1\over 2}(\zeta_i+\zeta_{i+1})$. If for instance $\zeta_i$ is
finite and $\zeta_{i+1}=+\infty$, we take $f=\zeta_i+1$.
\end{proof}

We introduce a notation. Given $x=(x_1,\ldots,x_n)$ and $y=(y_1,\ldots,y_n)$ in $R^n$,
we set $\Vert x-y\Vert= \sqrt{\sum_{i=1}^n (x_i-y_i)^2}$. Given $r>0$ in $R$, we define
the open and closed balls
\begin{eqnarray*}
B(x,r)&=&\{z\in R^n\;;\ \Vert z-x\Vert <r\}\;,\\
\overline{B}(x,r)&=&\{z\in R^n\;;\ \Vert z-x\Vert \leq r\}\;.
\end{eqnarray*}
The open balls $B(x,r)$ form a basis of open neighborhoods of $x$.

\begin{theorem}[Curve Selection Lemma]\label{curvesel}\index{curve selection lemma} Let $A$
be a definable subset of $R^n$, $b$ a point in $\clos(A)$. There is a
continuous definable map $\gamma :[0,1)\to R^n$ such that
$\gamma(0)=b$ and $\gamma((0,1))\subset A$.
\end{theorem}

\begin{proof} Set $X=\{ (t,x)\in R\times R^n\;;\ x\in A \et \Vert
x-b\Vert<t\}$. Let $p:R\times R^n\to R$ be the projection on the first coordinate. Since
$b\in \clos(A)$, we have $p(X)=\{ t\in R\;;\ t>0\}$. Applying
Definable Choice \ref{defchoice} and the Monotonicity Theorem
\ref{monotonicity}, we find $\varepsilon>0$ and a continuous
definable map $\delta : (0,\varepsilon)\to A$ such that $\Vert
\delta(t)-b\Vert < t$. Clearly, $\delta$ extends continuously to
$\overline{\delta}:[0,\varepsilon)\to R^n$ with
$\overline{\delta}(0)=b$, and we define
$\gamma:[0,1)\to R^n$ by $\gamma(t)=\overline\delta(t\varepsilon)$.
\end{proof}

The curve selection lemma replaces the use of sequences in many
situations.

\debexo
Show that a definable function $f:A\to R$ is continuous if and
only if, for every continuous definable $\gamma: [0,1)\to A$,
$\lim_{t\to 0_+}f(\gamma(t))=f(\gamma(0))$.\finexo

We have seen that compactness is a problem for a general real closed
field. However, the following result shows that we retain the good
properties of compactness if we deal with definable objects.

\begin{theorem}\label{compact} Let $A$ be a definable subset of
$R^n$. The following properties are equivalent.
\begin{enumerate}
\item $A$ is closed and bounded.
\item Every definable continuous map $(0,1)\to A$ extends by
continuity to a map $[0,1)\to A$.
\item For every definable continuous function $f:A\to R$, $f(A)$ is
closed and bounded.
\end{enumerate}\end{theorem}

\begin{proof} $1 \Rightarrow 2$. A definable continuous map $(0,1)\to
A$ extends by continuity to a map $[0,1)\to R^n$: every
coordinate of the map has a limit as $t\to 0_+$ (cf.\ Exercise
\ref{leftrightlimit}), and this limit is in $R$ since $A$ is bounded.
Since
$A$ is closed, the value of the extension at $0$ belongs to $A$. \par
$2\Rightarrow 3$. Suppose that $f(A)$ is not bounded. Set $$X=\{
(t,x)\in R\times R^n \;;\ x\in A \et t\,|f(x)|=1\}\;.$$ Then the
projection of $X$ on the first coordinate contains some interval
$(0,\varepsilon)$. Using Definable Choice
\ref{defchoice} and the Monotonicity Theorem
\ref{monotonicity}, and rescaling the interval of definition, we can
assume that there is a continuous map 
$\delta : (0,1)\to A$ such that $\lim_{x\to 0_+}
|f(\delta(t))|=+\infty$. This implies that $\delta$ cannot be
extended continuously to a map $[0,1)\to A$, which contradicts 2.
Hence,
$f(A)$ is bounded. Now let $b$ belong to $\clos(f(A))$. Set
$$Y=\{(t,x)\in R\times R^n\;;\ x\in A \et |b-f(x)|<t\}\;.$$
The same argument as above shows that there is a continuous map
$\gamma : (0,1)\to A$ such that $\lim_{t\to 0_+}
f(\gamma(t))=b$. By 2, $\gamma$ extends continuously to a map
$[0,1)\to A$ and its value at $0$ is an element $a$ of $A$ such that
$f(a)=b$. This shows that $f(A)$ is closed.
\par
$3\Rightarrow 1$.  Since the image of $A$ by each coordinate function
is bounded, $A$ is bounded. Let $b$ belong to $\clos(A)$. Since the
image of
$A$ by $x\mapsto \Vert x-b\Vert$ is closed, it contains $0$. Therefore
$b$ belongs to $A$, which shows that $A$ is closed.
\end{proof}

\debexo Show that a definable continuous function on a closed and
bounded definable set is uniformly continuous.
\finexo

\begin{corollary}\label{closedbddintrisic} Let $A$ be a closed and
bounded definable set. If $B$ is a definable set definably
homeomorphic to $A$, then $B$ is also closed and bounded.
\end{corollary}

The preceding corollary shows that the property of being closed and\linebreak
bounded, for a definable set, is intrisic (in the sense that it does
not depend on an imbedding in affine space). The property of being
locally closed is also intrisic, in the same sense. Recall that a
subset of $R^n$ is said to be locally closed if it is open in its
closure.

\begin{proposition}\label{locclosedintrisic} 1) A definable set
$A\subset R^n$ is locally closed if and only if every point $x\in A$
has a basis of closed and bounded definable neighborhoods in $A$.\par
2) If $A$ is a locally closed definable set, and $B$, a definable set
definably homeomorphic to $A$, then $B$ is locally closed.
\end{proposition}

\begin{proof}
1) First assume that $A$ is locally closed. Take $x\in A$. There is
$r_0>0$ such that the intersection of the ball $B(x,r_0)\subset R^n$
with
$\clos(A)$ is contained in $A$. Then the intersections of $\clos(A)$
with the closed balls $\overline{B}(x,r)$, for $0<r<r_0$, form a basis
of closed and bounded definable neighborhoods of $x$ in $A$. \par
Conversely, assume that every
point
$x\in A$ has a closed neighborhood
$N$ in $A$. Take $s>0$ such that $B(x,s)\cap A\subset N$. Since
$B(x,s)\cap (R^n\setminus N)$ is open and disjoint from $A$, we have
$B(x,s)\cap
\clos(A)\subset N\subset A$. This shows that $A$ is locally
closed.\par
2) This follows from 1) and Corollary \ref{closedbddintrisic}.
\end{proof}

\debexo\label{unionlocclosed} Show:\\
1) A cell of a cdcd is locally closed. (Use \ref{locclosedintrisic},
statement 2.)\\
2) Every definable set is the union of finitely many locally closed
definable sets.
\end{exo}

\section{Connected Components}

\begin{theorem}\label{cc} Let $A$ be a definable subset of $R^n$. There
is a partition of $A$ into finitely many definable subsets
$A_1,\ldots,A_k$ such that each $A_i$ is nonempty, open and closed in
$A$, and definably arcwise connected. Such a partition is unique. The
$A_1,\ldots,A_k$ are called the {\em definable connected components
of $A$}\index{connected! definable -- component}.
\end{theorem}

\begin{proof} Make a cdcd of $R^n$ adapted to $A$. We say that a cell
$C$ is adjacent to a cell $D$ if $C\cap\clos(D)\not=\emptyset$, and
we denote this fact by $C\prec D$. \par
We claim that, if $C\prec  D$, every $x\in C$ can be joined to every
$y\in D$ by a definable continuous path in $C\cup D$. Take $c\in
C\cap\clos(D)$. By the Curve Selection Lemma, there is a
continuous definable $\gamma: [0,1)\to C\cup D$ such that $\gamma(0)=c$
and $\gamma((0,1))\subset D$. Set $d=\gamma(1/2)$. The points $c$ and
$d$ are joined by a continuous definable path in $C\cup D$. Since
every cell is definably homeomorphic to an affine space, every point $x\in
C$ can be joined to $c$ by a definable continuous path in $C$ and 
every point $y\in D$ can be joined to $d$ by a definable continuous
path in $D$. This proves the claim.\par
Let $\sim$ be the smallest equivalence relation on the set of cells
contained in $A$ containing the adjacency relation: we have $C\sim D$
if and only if there is a chain $C=C_0\prec  C_1\succ C_2\prec \ldots
\succ C_{2k}=D$ where all cells $C_i$ are contained in $A$ (note that
we can have
$C_i=C_{i+1}$). Let ${\cal E}_1,\ldots,{\cal E}_k$ be the
equivalence classes for $\sim$, and let $A_i$ be the union of all cells
in
${\cal E}_i$. The $A_i$ form a finite partition of $A$ into definable
sets. The claim proved above and the definition of
$\sim$ imply that each $A_i$ is definably arcwise connected. If a
cell $C\subset A$ has a nonempty intersection with $\clos(A_i)$, it
is adjacent to a cell in ${\cal E}_i$, which implies
$C \subset A_i$. Hence, each $A_i$ is closed in $A$. Since the $A_i$
form a finite partition of $A$, they are also open in $A$. \par
Assume that $A=B_1\cup\ldots\cup B_\ell$ is another partition with
the same properties. We have $A_i=\bigcup_{j=1}^\ell (A_i\cap B_j)$,
and each $(A_i\cap B_j)$ is definable, open and closed in $A_i$.
Since $A_i$ is definably connected, there is exactly one $j$ such
that $A_i\subset B_j$. For the same reason, for every $j$ there is
exactly one $i$ such that $B_j\subset A_i$. This shows that the two
partitions coincide up to order.
\end{proof}

Remark that it follows from the proof that the number of definable
connected components of $A$ is not greater that the number of cells
contained in $A$ for any cdcd of $R^n$ adapted to $A$. 

\begin{corollary}\label{arcconn} A definably connected definable set is
definably arcwise connected.
\end{corollary}

Consider the case where $R={\Bbb R}$, i.e.\ the o-minimal structure
expands the field of real numbers. It is clear that a connected
definable set is definably connected. It is also clear that a
definably arcwise connected definable set is arcwise connected.
Hence, a definable set is connected if and only it is definably
connected, and the definable connected components of a definable set
are its usual connected components. It follows that a definable set
has finitely many connected components, which are definable.\par
\medskip

We shall prove now that there is a uniform bound for the number of
connected components in a definable family of subsets of $R^n$. First
we make precise this notion of definable family. Let $A$ be a
definable subset of $R^m\times R^n$. For $x\in R^m$, set $A_x=\{ y\in
R^n\;;\ (x,y)\in A\}$. We call the family $(A_x)_{x\in R^m}$ a {\em
definable family of subsets of $R^n$ parametrized by $R^m$}\index{definable family!-- of
sets}.\par

The next result allows to regard a cdcd of $R^m\times R^n$ as a
``definable family of cdcd of $R^n$''.  Assume that a cdcd of
$R^m\times R^n$ is given. Let $p_{m+n,m}: R^m\times R^n\to R^m$ be the
projection on the first $m$ coordinates. Recall that the
 $p_{m+n,m}(C)$, for all cells $C$ of $R^m\times R^n$, are the cells of a cdcd of $R^m$.

\begin{proposition}\label{cdcdfamily} Let $B$ be a cell of $R^m$, $a\in
B$. For all cells $C\subset R^m\times R^n$ such that $p_{m+n,m}(C)=B$,
we set
$C_a=\{y\in R^n\;;\ (a,y)\in C\}$. Then these $C_a$ are the cells of a
cdcd of
$R^n$. The dimension of the cell
$C_a$ is equal to $\dim(C)-\dim(B)$.
\end{proposition}

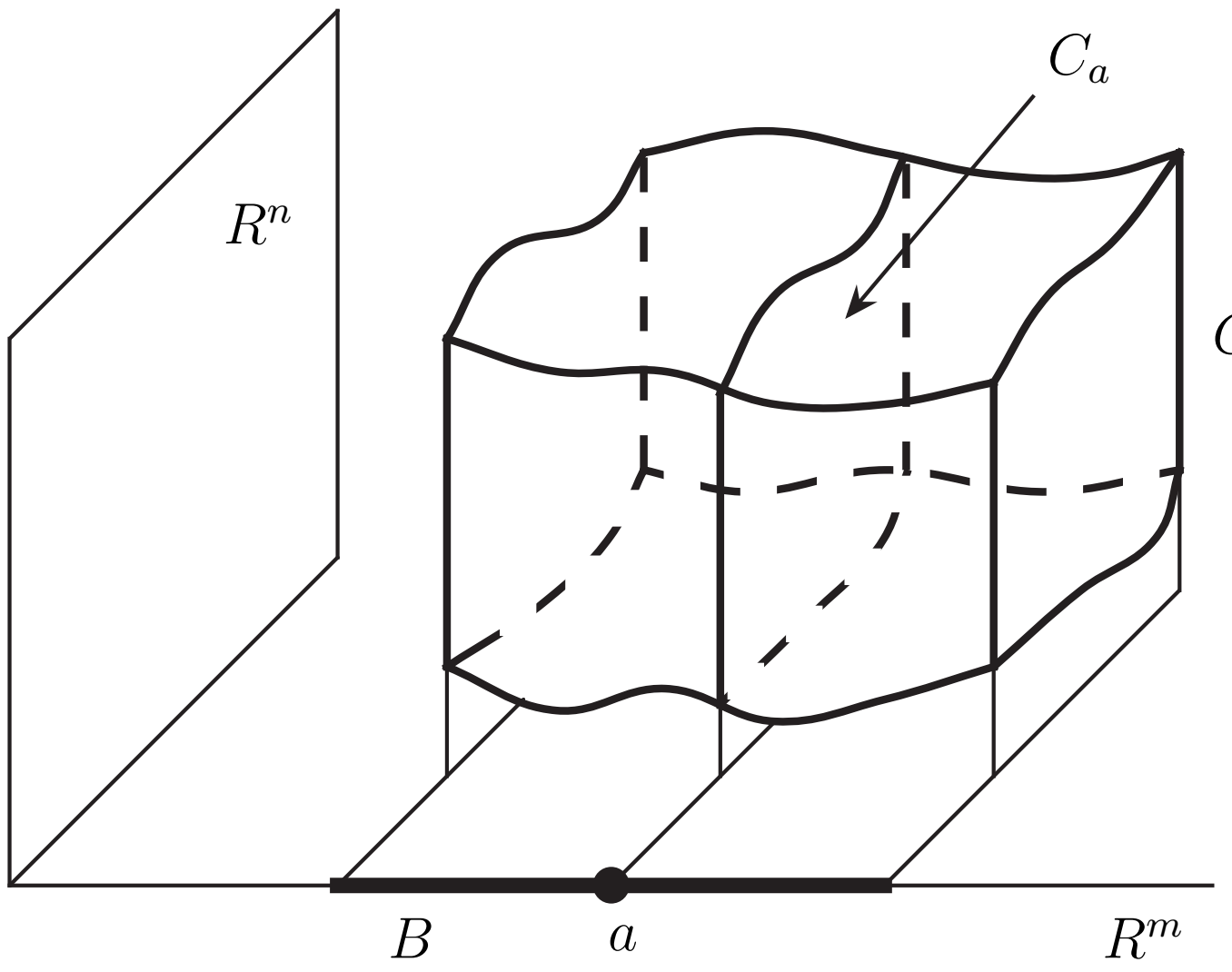
\begin{figure}
\setlength{\unitlength}{1mm}
\begin{picture}(112,72)
\put(20,6){\illustration{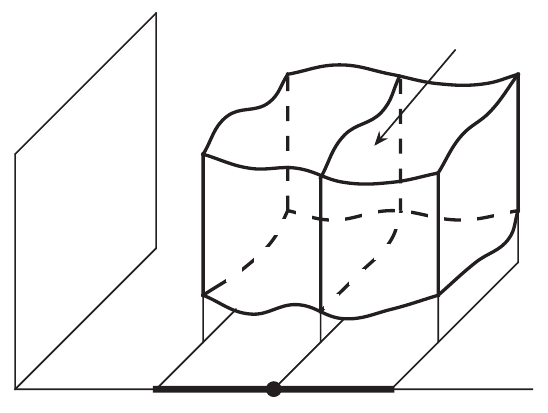}}
\put(48,4){\makebox(0,0)[l]{$B$}}
\put(64,4){\makebox(0,0)[l]{$a$}}
\put(100,4){\makebox(0,0)[l]{$R^m$}}
\put(108,48){\makebox(0,0)[l]{$C$}}
\put(96,68){\makebox(0,0)[l]{$C_a$}}
\put(36,56){\makebox(0,0)[l]{$R^n$}}
\end{picture}
\caption{The slice $C_a$ of a cell $C$.}
\end{figure}

\begin{proof} We proceed by induction on $n$. The case $n=0$ is
trivial. Now assume that $n>0$ and the proposition is proved for
$n-1$. The cdcd of $R^m\times R^n$ induces by the projection
$p_{m+n,m+n-1}:R^m\times R^n\to R^m\times R^{n-1}$ a cdcd of
$R^m\times R^{n-1}$. For each cell $D$ of $R^m\times R^{n-1}$, there
are continuous definable functions
$\zeta_{D,1}<\ldots<\zeta_{D,\ell(D)}:D\to R$, and the cells $C$ of
$R^m\times R^n$ such that $p_{m+n,m+n-1}(C)=D$ are graphs of
$\zeta_{D,i}$ or bands
$(\zeta_{D,i},\zeta_{D,i+1})$. By the inductive assumption, the
$D_a$, for all cells $D$ of $R^m\times R^{n-1}$ such that
$p_{m+n-1,m}(D)=B$, are the cells of a cdcd of $R^{n-1}$ and
$\dim(D_a)=\dim(D)-\dim(B)$. For such cells $D$, define
$(\zeta_{D,i})_a:D_a\to R$ by $(\zeta_{D,i})_a(y)=\zeta_{D,i}(a,y)$.
Now $C\subset R^m\times R^n$ is the graph of $\zeta_{D,i}$ (resp. the
band
$(\zeta_{D,i},\zeta_{D,i+1})$) if and only if $C_a\subset R^n$ is the
graph of $(\zeta_{D,i})_a$ (resp. the band
$((\zeta_{D,i})_a,(\zeta_{D,i+1})_a)$). Hence the $C_a$, for all
cells $C$ such that $B=p_{m+n,m}(C)$, form
a cdcd of
$R^n$. If $C$ is a graph over $D=p_{m+n,m+n-1}(C)$,
$\dim(C_a)=\dim(D_a)=\dim(D)-\dim(B)=\dim(C)-\dim(B)$. If $C$ is
a band, $\dim(C_a)=\dim(D_a)+1=\dim(D)-\dim(B)+1=\dim(C)-\dim(B)$.
\end{proof}

\begin{theorem}\label{uniboundcc} Let $A$ be a definable subset of
$R^m\times R^n$. There is $\beta\in {\Bbb N}$ such that, for every
$x\in R^m$ the number of definable connected components of $A_x$ is not greater
than $\beta$.
\end{theorem}

\begin{proof}
We choose a cdcd of $R^m\times R^n$ adapted to $A$. We adopt the
notation of Proposition \ref{cdcdfamily} and its proof. Take $x\in
R^m$ and let $B$ be the cell of $R^m$ containing $x$. The
collection of all $C_x$, for
$C$ cell of $R^m\times R^n$ such that $p_{m+n,m}(C)=B$, is a cdcd of
$R^n$ adapted to $A_x$.  Hence, the number of definable connected components of
$A_x$ is not greater that the number of cells of $R^m\times
R^n$ contained in $A$.
\end{proof}

\debexo The aim of this exercise is to prove the following fact: 
\begin{itemize}
\item[$(\sharp)$] Let $A$ be a definable subset of $R^n$.
There exists $\beta(A)\in {\Bbb N}$ such that, for every affine
subspace $L$ of $R^n$, the number of definable connected components of
$L\cap A$ is not greater than $\beta(A)$.
\end{itemize}
\hskip\retraitexo 1) (The Grassmannian of affine subspaces of $R^n$ as
an algebraic subset of $R^{n^2+n}$.) Let $M$ be a $n\times n$ matrix
with coefficients in
$R$, which we identify with an element of $R^{n^2}$, and let $v$ be a
vector in $R^n$. We associate to $g=(M,v)\in R^{n^2}\times R^n$
the affine subspace $L_g=\{ x\in R^n\;;\ Mx=v\}$ of $R^n$. Show that
the mapping $g\mapsto L_g$ restricted to the algebraic subset 
$$G=\{ g=(M,v)\in R^{n^2}\times R^n\;;\ M^2=M \et {}^tM=M \et
Mv=v\}$$
is a bijection onto the set of all affine subspaces of $R^n$.\\
2) Let $A$ be a definable subset of $R^n$.
Construct a definable set ${\cal A}\subset G\times R^n$ such
that, for all $g\in G$, ${\cal A}_g=A\cap L_g$ (where ${\cal A}_g$ is
defined as $\{ x\in R^n\;;\ (g,x)\in {\cal A}\}$).\\
3) Prove $(\sharp)$.
\finexo

\section{Dimension}

We have already defined the dimension of a cell of a cdcd. Now let
$A$ be a definable subset of $R^n$. Take a cdcd of $R^n$ adapted to
$A$. A ``naive'' definition of the dimension of $A$ is the maximum of
the dimension of the cells contained in $A$. But this definition is
not intrinsic. We have to prove that the dimension so defined does
not depend on the choice of a cdcd adapted to $A$. We shall rather
introduce an intrisic definition of dimension, and show that it
coincides with the ``naive'' one

\begin{definition}\label{dim} The dimension\index{dimension} of a definable set $A$ is
the sup of $d$ such that there exists a injective definable map from
$R^d$ to $A$. By convention, the dimension of the empty set is
$-\infty$. 
\end{definition}

Remark that it is not obvious for the moment that the dimension is
always $<+\infty$. It is also not clear that this definition of
dimension agrees with the one already given for cells. Both facts will
follow from the next lemma.

\begin{lemma}\label{invdomain}
Let $A$ be a definable subset of $R^n$ with nonempty interior. Let
$f: A\to R^n$ be an injective definable map. Then $f(A)$ has nonempty
interior.
\end{lemma}

\begin{proof} We prove the lemma by induction on $n$. If $n=1$, $A$
is infinite, hence $f(A)\subset R$ is infinite and contains an
interval. Assume that $n>1$ and the lemma is proved for all $m<n$.
Using Piecewise Continuity \ref{PC}, we can assume moreover that $f$ is
continuous. Take a cdcd of
$R^n$ adapted to
$f(A)$. If
$f(A)$ has empty interior, it contains no open cell. Hence $f(A)$ is
the union of nonopen cells $C_1,\ldots,C_k$ and, for $i=1,\ldots,k$,
there is a definable homeomorphism $C_i\to R^{m_i}$ with $m_i<n$.
Take a cdcd of $R^n$ adapted to the $f^{-1}(C_i)$. Since
$A=\bigcup_{i=1}^k f^{-1}(C_i)$ has nonempty interior, one of the
$f^{-1}(C_i)$, say $f^{-1}(C_1)$, must contain an open cell $B$. The
restriction of $f$ to $B$ gives an injective continuous definable map
$B\to C_1$. Since $B$ is definably homeomorphic to $R^n$ and $C_1$
definably homeomorphic to $R^m$ with $m<n$, we obtain an
injective continuous definable map $g: R^n\to R^m$. Set
$a=(0,\ldots,0)\in R^{n-m}$ and consider the mapping
$g_a:R^m\to R^m$ defined by $g_a(x)=g(a,x)$. We can apply the
inductive assumption to $g_a$. It implies that $g_a(R^m)$ has
nonempty interior in $R^m$. Take a point $c=g_a(b)$ in the interior
of $g_a(R^m)$. Since $g$ is  continuous we can find $x\in
R^{n-m}$, $x\not =a$ and close to $a$, such that
$g(x,b)\in g_a(R^m)$. There is $y\in R^m$ such that
$g(x,b)=g_a(y)=g(a,y)$, which contradicts the fact that $g$ is
injective. Hence, $f(A)$ has nonempty interior.
\end{proof}

\begin{corollary}\label{dimcell} The dimension of $R^d$ (according
to Definition
\ref{dim}) is $d$. The dimension of a cell, as defined in Section
\ref{sectioncdcd}, agrees with its dimension according to Definition
\ref{dim}.
\end{corollary}

\begin{proof} There is no injective definable map from $R^e$ to
$R^d$ if $e>d$. Otherwise, the composition of such a map with the
embedding of $R^d$ in $R^e=R^d\times R^{e-d}$ as
$R^d\times\{0\}$ would contradict Lemma \ref{invdomain}. This
shows the first part of the corollary. The second part follows, using
the fact that the dimension according to \ref{dim} is
invariant by definable bijection.
\end{proof}

\begin{proposition}\label{pptiesdim}\begin{enumerate}
\item If $A\subset B$ are definable sets, $\dim A \leq \dim B$.
\item If $A$ and $f:A\to R^n$ are definable, $\dim(f(A))\leq
\dim(A)$. If moreover $f$ is injective, $\dim(f(A))=\dim(A)$.
\item If $A$ and $B$ are definable subsets of $R^n$, $\dim (A\cup B) =
\max (\dim A, \dim B)$.
\item Let $A\subset R^n$ be definable and take a cdcd of $R^n$ adapted
to $A$. Then the dimension of $A$ is the maximum of the dimension of
the cells contained in $A$.
\item  If $A$ and $B$ are definable sets, $\dim(A\times B)=\dim
A+\dim B$.
\end{enumerate}\end{proposition}

\begin{proof} 1 is clear from the definition.\par
2. The second part is obvious since dimension is invariant by
definable bijection. If $f$ is definable, we get by Definable Choice
a definable mapping $g:f(A)\to A$ such that $f\circ g={\rm
Id}_{f(A)}$. Hence, $g$ is injective and $\dim (f(A))=\dim
(g(f(A)))\leq \dim(A)$.\par  
3. The inequality
$\dim (A\cup B)
\geq
\max (\dim A,
\dim B)$ follows from 1. Now let $f: R^d\to A\cup B$ be a definable
injective map. Taking a cdcd of $R^d$ adapted to $f^{-1}(A)$ and
$f^{-1}(B)$, we see that $f^{-1}(A)$ or $f^{-1}(B)$ contains a cell of
dimension $d$. Since $f$ is injective, we have $\dim A\geq d$ or
$\dim B\geq d$. This proves the reverse inequality $\dim (A\cup B)
\leq \max (\dim A,\dim B)$.\par
4 is an immediate consequence of 3.\par
5. By 3, it is sufficient to consider the case where $A$ and $B$ are
cells. Since $A\times B$ is definably homeomorphic to $R^{\dim
A}\times R^{\dim B}$, the assertion in this case follows from
Corollary \ref{dimcell}.
\end{proof}

Now we study the variation of dimension in a definable family. Recall
that a definable subset $A$ of $R^m\times R^n$ can be considered as a
definable family $(A_x)_{x\in R^m}$ of subsets of $R^n$, where $A_x=\{
y\in R^n\;;\ (x,y)\in A\}$.

\begin{theorem}\label{diminfamily}
Let $A$ be a definable subset of $R^m\times R^n$.  For
$d\in {\Bbb N}\cup\{-\infty\}$, set $X_d=\{x\in R^m\;;\ \dim(A_x)=d\}$.
Then $X_d$ is a definable subset of $R^m$, and $\dim(A\cap (X_d\times
R^n))= \dim(X_d) + d$.
\end{theorem}

\begin{proof}
We take a cdcd of $R^m\times R^n$ adapted to $A$ and use Proposition
\ref{cdcdfamily}. Let $B$ be a cell of $R^m$. For every $x$ in $B$,
$A_x$ is the union of the cells $C_x\subset R^n$, for all cells
$C\subset R^m\times R^n$ contained in $A$ such that
$p_{m+n,m}(C)=B$. Moreover, $\dim(C_x)= \dim C -\dim B$. It follows
that, for every $x\in B$, we have $\dim A_x= \dim(A\cap
(B\times R^n))-\dim B$. Hence, each $X_d$ is the union of some cells
$B$ of
$R^m$. This implies that $X_d$ is definable. Since $\dim(A\cap
(B\times R^n))= \dim B + d$ for a cell $B\subset X_d$, we have
$\dim (A\cap (X_d\times R^n))= \dim X_d + d$.
\end{proof}

\debexo\label{smallfibsmallfam} Let $A$ and $B$ be subsets of
$R^m\times R^n$, with $A$ nonempty. Assume that, for every $x\in
R^m$, $\dim (B_x) <\dim (A_x)$ or $B_x$ is empty. Prove that
$\dim B < \dim A$.
\finexo 

\debexo Let $f:A\to R^m$ be a definable map. Prove that the set
$Y_d$ of $x\in R^m$ such that $\dim(f^{-1}(x))=d$ is definable. Prove
that $\dim(Y_d)= \dim(f^{-1}(Y_d)) -d$.
\finexo

We finish this section on dimension by showing that the dimension
behaves well with respect to closure. The following lemma will be
useful.

\begin{lemma}\label{closfibandfibclos} Let $A$ be a definable subset of
$R^m\times R^p$. Let
$M$ be the set of $x\in R^m$ such that the closure of $A_x$ in $R^p$
is different from $(\clos(A))_x$. Then $M$ is definable and
$\dim M <m$. In particular, if $m=1$, $M$ is finite.
\end{lemma}

\begin{proof} Note that we have always
$\clos(A_x)\subset (\clos(A))_x$, since $(\clos(A))_x$ is closed and
contains $A_x$.\par
We leave the verification of the definability of
$M$ as an exercise. Suppose that
$\dim(M)=m$. Then $M$ contains an open cell $C$ of a cdcd of $R^m$. For
every
$x\in C$, we can find a box $B=\prod_{i=1}^p (a_i,b_i)$ such that
$B\cap (\clos(A))_x\not=\emptyset$ and $B\cap A_x=\emptyset$. By
Definable Choice \ref{defchoice} and Piecewise Continuity \ref{PC}, we
can assume that all
$a_i$ and
$b_i$ are definable continuous functions of $x\in C$. Let $U$ be the
set of $(x,y_1,\ldots,y_p)\in C\times R^p$ such that
$a_i(x)<y_i<b_i(x)$ for $i=1,\ldots,p$. The set $U$ is open in
$R^m\times R^p$, disjoint from $A$ and has nonempty intersection
with $\clos(A)$. This is impossible. Hence, $\dim M <m$.
\end{proof}

\begin{theorem}\label{dimboundary} Let $A$ be a nonempty definable
subset of
$R^n$. Then $$\dim(\clos(A)\setminus
A)<\dim(A)\;.$$ 
\end{theorem}

It follows from the theorem that $\dim(\clos(A))=\dim A$.

\begin{proof}
We proceed by induction on $n$. The theorem is obvious for $n=1$. We
assume that $n>1$ and the theorem is proved for $n-1$. We denote by
$\xi_1,\ldots,\xi_n$ the coordinate functions on $R^n$. For
$i=1,\ldots,n$, let $\clos_i(A)$ be the set of $x\in R^n$
such that
$x$ belongs to the closure of the intersection of $A$ with the
hyperplane $\xi_i^{-1}(\xi_i(x))$. \par
{\it First step.} We claim that $\clos(A)\setminus A$ has
dimension not greater than the maximum of $0$ and
$\dim(\clos_i(A)\setminus A)$ for
$i=1,\ldots,n$. We have $\clos_i (A)\subset
\clos(A)$. Applying Lemma \ref{closfibandfibclos} after a permutation
of the coordinates which puts the $i$-th coordinate in the first
place, we obtain that the difference
$\clos(A)\setminus\clos_i(A)$ is contained in finitely many
hyperplanes $\xi_i^{-1}(a_{i,j})$, for
$j=1,\ldots,\ell(i)$. Hence,
$\clos(A)\setminus\bigcup_{i=1}^n\clos_i A$ is contained in the finite
set consisting of the $\prod_{i=1}^n \ell(i)$ points $(a_{1,j_1},\ldots,a_{n,j_n})\in
R^n$. The claim follows.
\par
{\it Second step.} We claim that $\dim(\clos_i(A)\setminus A)<\dim A$
for $i=1,\ldots,n$. Take $a\in R$. Since the hyperplane
$\xi_1^{-1}(a)$ has dimension $n-1$, the inductive assumption implies
that 
$\clos(A\cap\xi_1^{-1}(a))\setminus (A\cap\xi_1^{-1}(a))$ is empty or
has dimension strictly smaller than the dimension of
$A\cap\xi_1^{-1}(a)$; note that the first set is $(\clos_1(A)\setminus
A)\cap\xi_1^{-1}(a)$. Hence, applying the result of Exercise
\ref{smallfibsmallfam}, we obtain
$\dim(\clos_1(A)\setminus A)<\dim A$.  Permuting the coordinates, 
we prove the claim.\par
Now we can complete the proof of the theorem. Steps 1 and 2 imply
that $\dim(\clos(A)\setminus A)\leq 0$ or $\dim(\clos(A)\setminus
A) <\dim A$. If $\dim A=0$, $A$ is closed. Hence, for every
nonempty $A$,  $\dim(\clos(A)\setminus A) <\dim A$.
\end{proof}

%% file: OMINquatre.tex
\chapter{Definable Triangulation}

In this chapter we show that the topology of definable sets can be entirely encoded
in finite terms. This will be done by means of a triangulation.

\section{Good Coordinates}

\begin{figure}[ht]
\setlength{\unitlength}{1mm}
\begin{picture}(112,65)
\put(20,0){\illustration{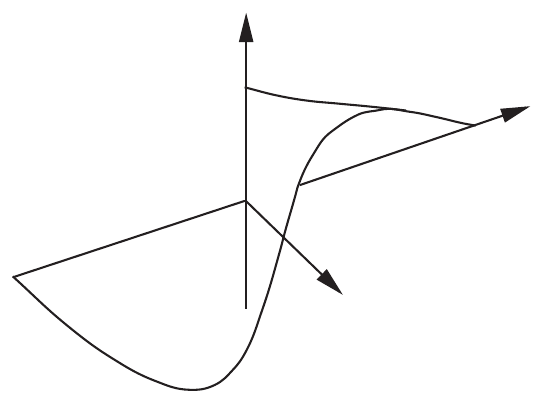}}
\put(78,15){\makebox(0,0)[l]{$x$}}
\put(110,49){\makebox(0,0)[l]{$y$}}
\put(63,65){\makebox(0,0)[l]{$z$}}
\end{picture}
\caption{$\Gamma\subset R^3$ near the origin\label{gamma}}
\end{figure}
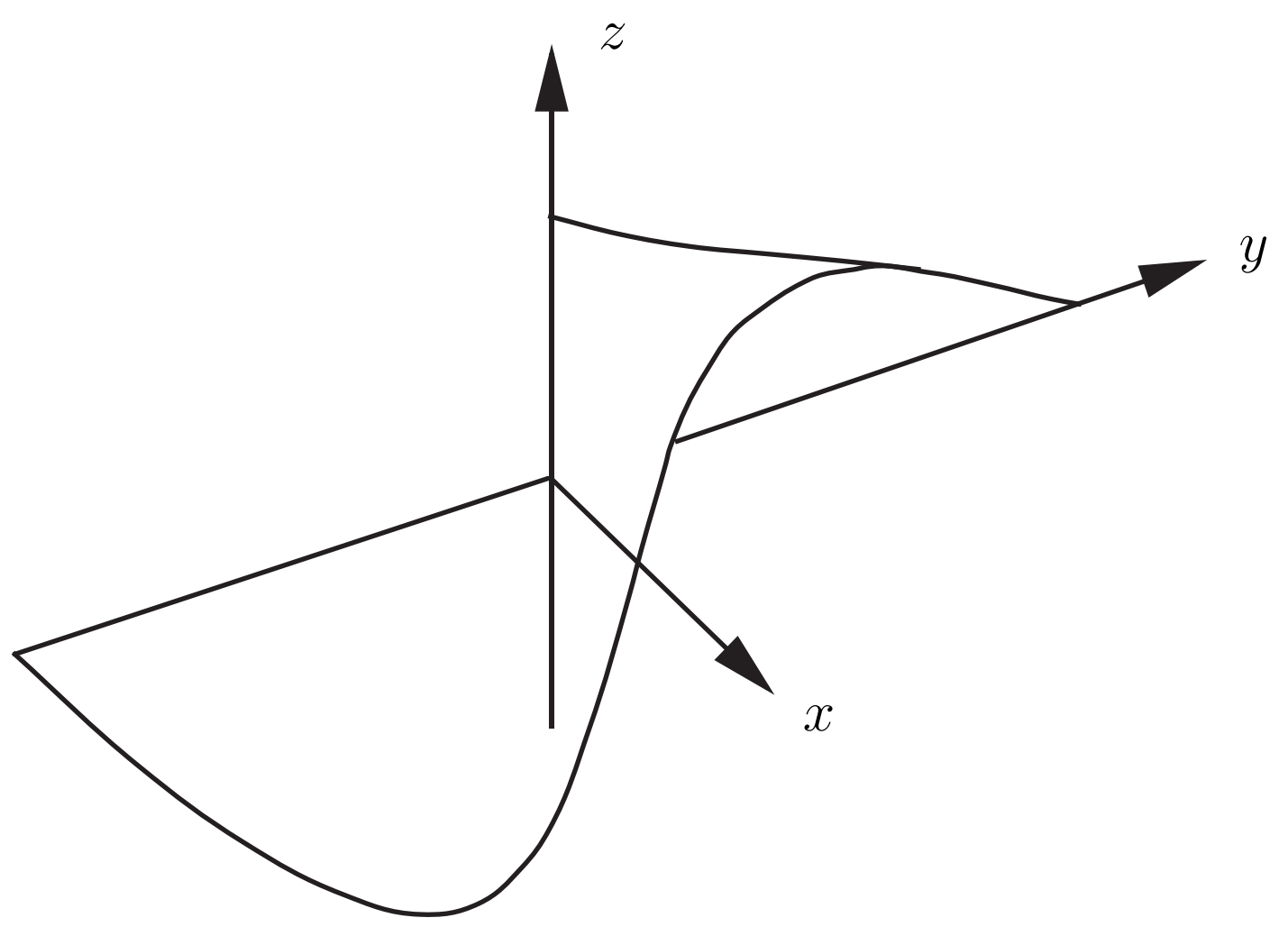

As we have seen in the last chapter, the cdcd is a very powerful tool. But it does
not give sufficient control on the relative disposition of the cells, when they are
not contained in the same cylinder. In particular, one cannot, in general,
reconstruct the topology of a definable set from its decomposition into cells of an
adapted cdcd. The main difficulty is that we have no control on how a definable
continuous function $\zeta:C\to R$ on a cell $C$ behaves as one approaches the
boundary of $C$. The function $\zeta$, even if it is bounded, may not extend to a
continuous function on $\clos(C)$. For instance, the definable continuous function
$\zeta$ defined on the set of $(x,y)\in R^2$ such that $x>0$ by
$\zeta(x,y)=2xy/(x^2+y^2)$ does not extend continuously to $(0,0)$. All points
$(0,0,z)$ with $-1\leq z\leq 1$ belong to the closure $\Gamma$ of the graph of
$\zeta$. The definable set $\Gamma\subset R^3$ has dimension $2$, but the
restriction to $\Gamma$ of the projection $R^3\to R^2$ on the first two
coordinates is not finite-to-one (see Figure \ref{gamma}).
\par
 In this section, we consider the following
problem: given a definable subset $G$ of $R^n$ of dimension $<n$, can we make a
polynomial change of coordinates in $R^n$ such that the restriction to $G$ of
the projection on the first $n-1$ new coordinates becomes finite-to-one? In other
words, we want a polynomial automorphism $u: R^n\to R^n$ such that
$p_{n,n-1}|_{u(G)}:u(G)\to R^{n-1}$ is finite-to-one. We look for $u$
 of the form $u(x',x_n)=(x'-v'(x_n),x_n)$, where $x'=(x_1,\ldots,x_{n-1})$ and
$v': R\to R^{n-1}$ is a polynomial map. The condition
that $p_{n,n-1}|_{u(G)}$ is finite-to-one reads: 
$$\qquad(*)\qquad\cases{ 
\mbox{for all }y'\in R^{n-1},\mbox{ the set of }x_n\in R&\cr
\mbox{ such that }(y'+v'(x_n),x_n)\in G\mbox{ is finite.}\cr}$$  
We introduce the definable set
$$W=\{ (y',x',x_n)\in R^{n-1}\times R^{n-1}\times R^n\;;\ (y'+x',x_n)\in G\}\;,$$
and we set, for $y'\in R^{n-1}$, $W_{y'}=\{ x\in R^n\;;\ (y',x)\in W\}$. Note that, 
for all $y'\in R^{n-1}$, $W_{y'}=G-(y',0)$ has dimension $<n$. We can translate the
condition $(*)$ as follows: 
$$\mbox{for all }y'\in R^{n-1},\mbox{ the set of }x_n\in R\mbox{ such that }
(v'(x_n),x_n)\in W_{y'}\mbox{ is finite.}$$
 Hence, our problem will be solved by the
following lemma. In the statement of the lemma, we exchange $x_1$ and $x_n$; it
will make the proof by induction easier to write. 
  
\begin{lemma}\label{evite} Let $W\subset R^m\times R^n$ ($n\geq 2$) be a
definable set. For
$s\in R^m$, define $W_s=\{ y\in R^n\;;\ (s,y)\in W\}$. Assume that, for all
$s$ in $R^m$, $\dim(W_s)<n$. Then there exist a polynomial map $v' : R\to
R^{n-1}$ of degree not greater than $m$ such that, for all $s$ in
$R^m$, the set of $x_1\in R$ such that $(x_1,v'(x_1))\in W_s$ is finite.  
\end{lemma}

\begin{proof}
We proceed by induction on  $n$. We begin with $n=2$. Let $V$
be a definable subset of dimension $1$ of $R^2$. A cdcd of $R^2$ adapted to $V$
decomposes $V$ as the disjoint union of finitely many points, vertical open intervals
and graphs of definable continuous functions
$\zeta_i:I_i\to R$, where $I_i$ is an interval. Consider such a function
$\zeta_i$. For $a=(a_0,a_1,\ldots,a_m)\in R^{m+1}$, define $f_a:R\to R$ by
$f_a(t)=\sum_{i=0}^m a_it^i$. The set
$$\displaylines{
\qquad\Big\{x_1\in I_i\;;\  \exists a\in R^{m+1}\ \exists \epsilon\in
R\hspace*{\fill}\cr
\hspace*{\fill}\Big(\epsilon>0\ {\rm and}\  
\forall y \in I_i\ \big(|y-x_1|<\epsilon \Rightarrow f_a(y)=\zeta_i(y)\big)\Big)\Big\}
\qquad\cr}$$ is
a definable open subset of $R$. It has finitely many definable connected
components which are open intervals. If $U$ is one of its definable connected
components, there is a unique $a$ such that
$f_a|_U=\zeta_i|_U$. Hence, there are finitely many $a\in R^{m+1}$ such that the set of
$x_1\in R$ such that
$(x_1,f_a(x_1))\in V$ is infinite. From this we deduce that, for all $s\in R^m$,
there are finitely many
$a\in R^{m+1}$ such that the set $B_{a,s}=\{x_1\in R\;;\ (x_1,f_a(x_1))\in W_s\}$ is
infinite. Therefore the definable set of $(a,s)$ such that $B_{a,s}$ is infinite has
dimension at most $m$. Then also the set of
$a\in R^{m+1}$, such that there is
$s\in R^m$ such that $B_{a,s}$ is infinite, has dimension at most $m$. Hence,
there exists a polynomial
$f_a$ of degree not greater than $m$ such that, for all $s\in R^m$, the set $B_{a,s}$ is
finite.\par  
Given $n>2$, assume the lemma proved for $n-1$. Let
$Z$ be the definable set of $(s,x_1,u) \in R^m\times R\times R^{n-2}$ such that
$\{x_n\in R\;;\  (s,x_1,u,x_n)\in W\}$ is infinite. For all $s\in R^m$,
the set $Z_s$ has dimension at most $n-2$ because $W_s$ has dimension at most $n-1$. 
Therefore, we can apply the inductive assumption to obtain a polynomial map $g:R\to
R^{n-2}$ of degree at most
$m$ such that, for all $s\in R^m$, the set of
$x_1\in R$ such that $(x_1,g(x_1))\in Z_s$ is finite. Consider the definable
subset of
$R^m\times R^2$
$$M=\{(s,x_1,x_n)\in R^m\times R^2\;;\  (s,x_1,g(x_1), x_n)\in W\}\;.$$ 
For all $s\in R^m$, the set $M_s$ has dimension at most $1$. Therefore,
by the argument above, there is a polynomial $f$ of degree at most $m$ such
that for all
$s\in R^m$ the set of $x_1\in R$ such that
$(x_1,f(x_1))\in M_s$ is finite. Set $v'=(g,f)$, and the proof is complete.
\end{proof}

The lemma actually solves more than our initial problem. 

\begin{proposition}\label{goodcoord} Let $G$ be a definable subset of $R^q\times
R^n$.
Assume that, for every
$t$ in $R^q$, the dimension of $G_t=\{x\in R^n\mid (t,x)\in G\}$ is $<n$.
Let 
$p:R^q\times R^n\to R^q\times R^{n-1}$ be the projection on the
first
$q+n-1$ coordinates. Then there is a polynomial automorphism $u$ of $R^n$ such
that
$$\forall (t,x') \in R^q\times R^{n-1}\quad 
p^{-1}(t,x')\cap \left(R^q\times u(G_t)\right)\ \mbox{ is
finite}.$$ 
\end{proposition}

\begin{proof} 
 Set
\begin{eqnarray*} 
W&=&\{(t,y',x',x_n)\in R^q\times R^{n-1}\times R^{n-1}\times R
\;;\ (t,y'+x',x_n)\in G\}\\ 
W_{(t,y')}&=&\{(x',x_n)\in R^{n-1}\times R \;;\ 
(t,y',x',x_n) \in W\}\;.
\end{eqnarray*}
For all $(t,y')$ in $R^q\times R^{n-1}$, the dimension of
$W_{(t,y')}= G_t-(y',0)$ is $< n$. By Lemma \ref{evite} (exchanging $x_1$ and $x_n$),
there is a polynomial mapping
$v':R\to R^{n-1}$ such that the set of $x_n$ in $R$ such that
$(v'(x_n),x_n)\in W_{(t,y')}$ is finite. Now set
$u(x',x_n)=(x'-v'(x_n),x_n)$. Then $u$ is a polynomial automorphism
of $R^n$ which satisfies the condition of the proposition.
\end{proof}

The case $q=0$ is our initial problem. In fact, this problem can be solved with $u$ a
linear automorphism (see \cite{vdD}, Theorem 4.2). The case $q=1$ will be useful for
the triangulation of functions. The point of having a finite-to-one projection is
explained by the following proposition.

\begin{proposition}\label{contextension} Let $F$ be a closed and bounded definable subset
of
$R^n$, such that $p_{n,n-1}|_F$ is finite-to-one. Let
$X\subset p_{n,n-1}(F)$ be a definable subset of
$R^{n-1}$ such that every
$x'\in\clos(X)$ has a basis of neighborhoods $U$ such that $U\cap X$ is definably
connected (this is the case, for instance, if $X$ is convex). Then every continuous
definable function
$\zeta:X\to R$ whose graph is contained in $F$ extends continuously to $\clos(X)$.
\end{proposition}

\begin{proof} Let $\Gamma\subset F$ be the graph of $\zeta$. Its closure is bounded,
and\linebreak
$p_{n,n-1}(\clos(\Gamma))=\clos(X)$. Take $x'\in\clos(X)$. The set $\clos(\Gamma)\cap
p_{n,n-1}^{-1}(x')$ is nonempty and finite. Take $a\in R$ and $\delta>0$ such that
$(\{x'\}\times (a-\delta,a+\delta))\cap
\clos(\Gamma)= (x',a)$. For every $\varepsilon$ such that $0<\varepsilon<\delta$, there
is neighborhood $U$ of
$x'$ such that $U\cap X$ is definably connected, $U\cap
\zeta^{-1}(a-\varepsilon,a+\varepsilon)\not=\emptyset$, and $\zeta(y')\not=
a\pm\varepsilon$ for every $y'\in X\cap U$. It follows that $\zeta(U\cap X)\subset
(a-\varepsilon,a+\varepsilon)$. This proves that $\zeta$ extends continuously to
$x'$, with value $a$.\end{proof}

\section{Simplicial Complex}

We recall some definitions concerning simplicial complexes that we shall need.
Let  $a_0 ,\ldots, a_d$ be points of $R^n$ which are affine independent (i.e.\
not contained in an affine subspace of dimension $d-1$). The $d$-{\em
simplex\/}\index{simplex} with vertices
$a_0 ,\ldots, a_d $ is
$$\displaylines{\quad[a_0 ,\ldots, a_d] =\hspace*{\fill}\cr
\hspace*{\fill} \{ x \in  R^n\  ;\ \exists
\lambda_0 ,\ldots,
\lambda_d
\in  [0,1]\quad  \sum_{i=0}^d \lambda_i = 1 
\et x = \lambda_0 a_0 +\ldots+ \lambda_d a_d \}\;.\quad\cr}$$
\begin{figure}[ht]
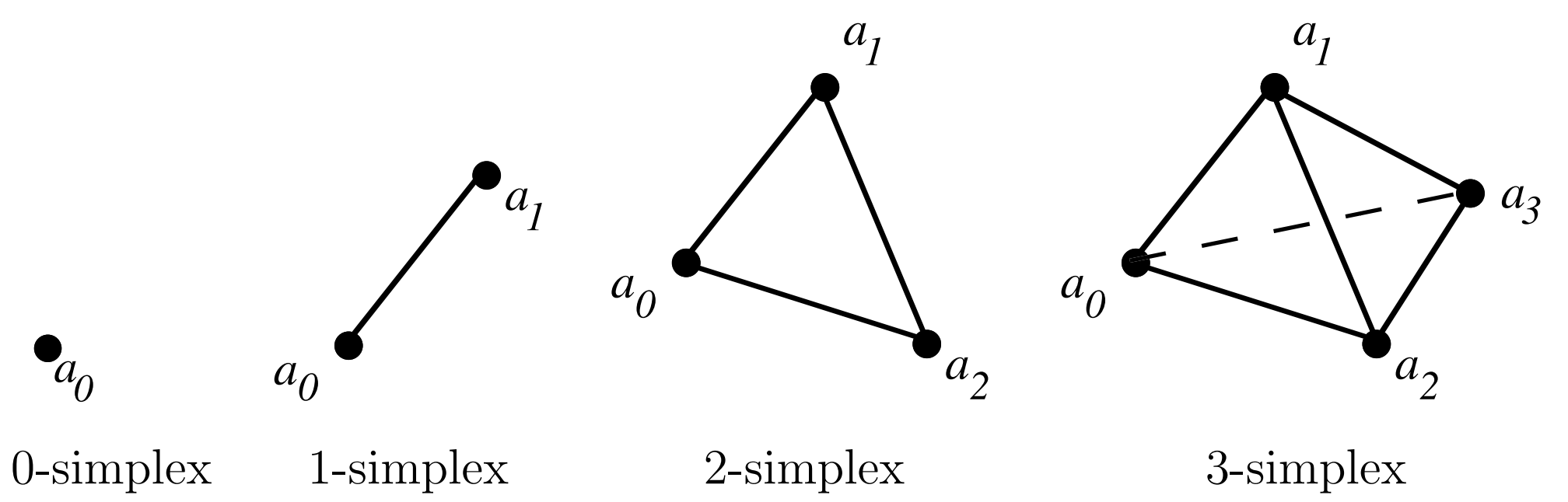

\hspace*{\fill}
\vbox{
\hbox{\hfil\dessin 7.1 mm by 8.2 mm (0-simp)\hfil}
\hbox{0-simplex}}
\hfill
\vbox{
\dessin 25 mm by 23 mm (1-simp)
\hbox to 25 mm{\hfil 1-simplex\hfil}}
\hfill
\vbox{
\dessin 34mm by 36 mm (2-simp)
\hbox to 34mm{\hfil 2-simplex\hfil}}
\hfill
\vbox{
\dessin 43 mm by 36 mm (3-simp)
\hbox to 43mm{\hfil 3-simplex\hfil}}
\hspace*{\fill}
\caption{\label{simplices}Simplices}
\end{figure} 

The corresponding {\em open simplex\/}\index{simplex!open --} is
$$\{ x \in  R^n\  ; \ \exists \lambda_0 ,\ldots,
\lambda_d \in  (0,1]\  \sum_{i=0}^d \lambda_i = 1  \et x=\lambda_0 a_0 +\cdots+
\lambda_d a_d \}$$  We shall denote by $\sigma$ the open simplex
corresponding to the simplex $\overline\sigma$.  A {\em face}\index{face} of the simplex
$\overline\sigma = [a_0 ,\ldots, a_d]$ is a simplex $\overline\tau  = [b_0 ,\ldots, b_e]$
such that
$$\{b_0 ,\ldots, b_e \} \subset  \{a_0 ,\ldots, a_d \}\;.$$
 
\begin{figure}[hb]
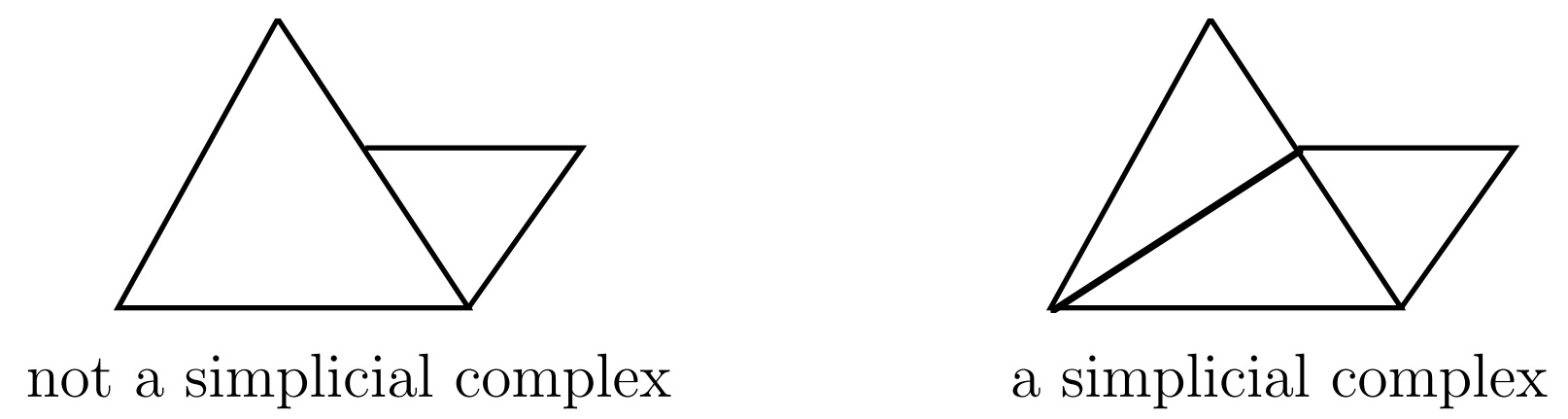

\noindent
\hspace*{\fill}
\vbox{
\hbox to 50mm{\hfil\dessin 33.1mm by 21.1 mm (pas_comp)\hfil}
\hbox to 50mm{\hfil			         not a simplicial complex\hfil}}
\hfill
\vbox{
\hbox to 50mm{\hfil\dessin 33.1mm by 21.1 mm (comp)\hfil}
\hbox to 50mm{\hfil			      a simplicial complex\hfil}}
\hspace*{\fill}
\caption{\label{simpcomp}Simplicial complex}
\end{figure}

A {\em finite simplicial complex\/} in $R^n$ is a finite collection $K =
\{\overline\sigma_1 ,\ldots, \overline\sigma_p \}$ of simplices $\overline\sigma_i
\subset  R^n$ such that, for every
$\overline\sigma_i,
\overline\sigma_j \in  K $, the intersection $\overline\sigma_i \cap \overline\sigma_j$
is a common face of
$\overline\sigma_i$ and $\overline\sigma_j$ (see Figure \ref{simpcomp}).

We set $\vert K\vert  = \bigcup_{\overline\sigma_i \in  K} \overline\sigma_i $; this is a
semialgebraic subset of
$R^n$. A {\em polyhedron\/} in $R^n$ is a subset $P$ of $R^n $,
such that there exists a finite simplicial complex $K$ in $R^n$ with $P = \vert
K\vert $. Such a $K$ will
be called a {\em simplicial decomposition} of
$P$. Note that a polyhedron is a closed and bounded definable set. In the following, it
will be convenient to agree that if a simplex
$\overline\sigma$ belongs to a finite simplicial complex $K$, then all faces of
$\overline\sigma$ also belong to $K$. With this convention, $\vert K\vert$  is the
disjoint union of all open simplices
$\sigma$ for $\overline\sigma \in  K$.
\par
We shall also use the notion of a {\em cone}\index{cone}. Let $B$ be a polyhedron in
$R^n$, and
$a\in(R^n\setminus B)$ such that every half-line from $a$ intersects $B$ in at most one
point (i.e., for every $x\in B$, $[a,x]\cap B=\{x\}$). The cone with base $B$ and
vertex $a$ is the polyhedron 
$$a*B=\{ tx+(1-t)a\;;\ x\in B \et t\in [0,1]\}\;.$$ 
Given a simplicial decomposition of
$B$, we obtain a simplicial decomposition of $a*B$ by taking all $a*\overline\sigma$,
for $\overline\sigma$ a simplex of the simplicial decomposition of $B$

\section{Triangulation of Definable Sets}

\begin{theorem}\label{triangsets} Let $A$ be a closed and bounded definable
subset of $R^n$ and let
$B_i$, $i=1,\ldots,k$, be definable subsets of $A$. Then there exist a finite
simplicial complex
$K$ with vertices in ${\Bbb Q}^n$ and a definable homeomorphism  $\Phi:
|K|\to A$ such that each $B_i$ is a union of images by
$\Phi$ of open simplices of $K$.
\end{theorem}

\begin{proof}
  We proceed by induction on $n$. The case of $n=1$ is obvious. We can subdivide
$R$ with finitely many points $x_1<\ldots<x_p$ such that $A$ and the $B_i$ are unions of
points $x_i$ and intervals $(x_j,x_{j+1})$. Then we choose a definable homeomorphism
$\tau : R\to R$ such that $\tau(x_i)=i$ for $i=1,\ldots,p$. We take for
$\Phi$ the restriction of $\tau^{-1}$ to $\tau(A)$.\par
  Now assume that
$n>1$ and the theorem is proved for $n-1$. Since every definable set is a
finite union of locally closed definable sets (cf. Exercise \ref{unionlocclosed}), we may
assume that the
$B_i$ are locally closed. Then we can replace $B_i$ with its closure
$\clos(B_i)$ and the difference $\clos(B_i)\setminus B_i$. Hence, we can assume
that all  $B_i$ are closed.
 Let
$F_0$ be the boundary $A\cap \clos(R^n\setminus A)$ of
$A$ and 
$F_i$ the boundary of $B_i$, for $i=1,\ldots,k$. Set $F=\bigcup_{i=0}^k F_i$. Then
$F$ is a closed and bounded definable set of dimension $<n$. Let 
$p=p_{n,n-1}:R^n\to R^{n-1}$ be the projection on the first $n-1$ coordinates. By
Proposition \ref{goodcoord}, we can assume that, for all $x'\in
R^{n-1}$, $p^{-1}(x')\cap F$ is finite.\par
  We make a cdcd of $R^n$ adapted to
$F_0,F_1,,\ldots,F_k$. We get a finite partition of
$p(A)=p(F_0)$ into definably connected definable subsets
$X_\lambda$ of $R^{n-1}$, and definable continuous functions
$$\zeta_{\lambda,1}<\ldots<\zeta_{\lambda,m(\lambda)} : X_\lambda\longrightarrow R\;,$$ such
that every graph of $\zeta_{\lambda,\mu}$ is contained in some
$F_i$ and every
$F_i$ is a union of graphs of  $\zeta_{\lambda,\mu}$. Since every band
$(\zeta_{\lambda,\mu},\zeta_{\lambda,\mu+1})$ is definably connected and disjoint from the
boundaries $F_i$, it is contained in or disjoint from each one of $A,B_1,\ldots,B_k$.\par
 Applying the inductive assumption, we obtain a simplicial complex
$L$ with vertices in ${\Bbb Q}^{n-1}$ and a definable homeomorphism $\Psi:|L|\to p(A)$
and such that all
$X_\lambda$ are the images by $\Psi$ of unions of open simplices of $L$. 
Replacing $A$ with
$\{(x',x_n)\in|L|\times R\;;\ (\Psi(x'),x_n)\in A\}$, we can assume that $\Psi$ is
the identity. Moreover, we can partition the $X_\lambda$ and assume that they are the
open simplices $\sigma_\lambda$ of 
$L$ (the open simplices are no longer cells of a cdcd). We denote by
$(C_\alpha)_{\alpha=1,\ldots,\ell}$ the collection of all graphs of
$\zeta_{\lambda,\mu}:\sigma_\lambda\to R$ and all bands
$(\zeta_{\lambda,\mu},\zeta_{\lambda,\mu+1})\subset \sigma_\lambda\times R$ which are
contained in $A$. By Proposition \ref{contextension}, since
$F$ is closed and bounded and $p|_{F}:F\to R^{n-1}$ is finite-to-one, every
$\zeta_{\lambda,\mu}:\sigma_\lambda\to R$ can be continuously extended to the closed
simplex
$\overline{\sigma}_\lambda$.  We denote the extension by
$\overline{\zeta}_{\lambda,\mu}$. The graph of this extension is contained in $F$. It
follows that, if
$\sigma_{\lambda'}$ is contained in
$\overline{\sigma}_\lambda$,  the restriction of $\overline{\zeta}_{\lambda,\mu}$ to
$\sigma_{\lambda'}$ coincides with some $\zeta_{\lambda',\mu'}$. We denote by
$\overline{C}_\alpha$ the closure of $C_\alpha$. It is either the graph of some
$\overline{\zeta}_{\lambda,\mu}$ or some closed band
$[\overline{\zeta}_{\lambda,\mu},\overline{\zeta}_{\lambda,\mu+1}]\subset
\overline{\sigma}_\lambda \times R$ (obvious notation). The set $\partial
C_\alpha=\overline{C}_\alpha\setminus C_\alpha$ is a union of $C_\beta$ with $\dim
(C_\beta)<\dim (C_\alpha)$.
\par For every simplex $\overline\sigma_\lambda$ of $L$, let $b(\sigma_\lambda)\in {\Bbb
Q}^{n-1}$ be its barycenter. If $C_\alpha$ is the graph of
$\zeta_{\lambda,\mu}:\sigma_\lambda\to R$, we set $b_\alpha=(b(\sigma_\lambda),\mu)\in
{\Bbb Q}^n$. If $C_\alpha$ is the band $(\zeta_{\lambda,\mu},\zeta_{\lambda,\mu+1})$, we
set
$b_\alpha=(b(\sigma_\lambda),\mu+{1\over 2})\in {\Bbb Q}^n$. Now we associate to
each $\overline{C}_\alpha$ a polyhedron $\overline{D}_\alpha$, by
induction on $\dim(\overline{C}_\alpha)$. Moreover, we give a simplicial decomposition of
$\overline{D}_\alpha$. If
$\overline{C}_\alpha$ is a point, we set $\overline{D}_\alpha=\{b_\alpha\}$. If
$\dim(\overline{C_\alpha})>0$, then $\overline{D_\alpha}$ is the cone with vertex
$b_\alpha$ and base the union $\partial D_\alpha$ of all $\overline{D}_\beta$ such
that
$C_\beta\subset \partial{C}_\alpha$. We decompose $\overline{D_\alpha}$ by
taking the cones with vertex $b_\alpha$ over all simplices of the simplicial
decompositions of
$\overline{D}_\beta\subset\partial D_\alpha$. We obtain in this way a finite simplicial
complex
$K$ such that
$|K|=\bigcup_{\alpha=1}^\ell \overline{D}_\alpha$. The simplices of $K$ are all
simplices $[b_{\alpha_0},b_{\alpha_1},\ldots,b_{\alpha_q}]$ such that
$C_{\alpha_{i-1}}\subset \partial{C}_{\alpha_i}$ for
$i=1,\ldots,q$. The complex $K$ has all its vertices in ${\Bbb Q}^n$. 
% The projection
% $p:R^n\to R^{n-1}$ induces a simplicial morphism from
% $K$ to the barycentric subdivision $L'$ of $L$.
Note that, by construction, $p(\overline{D}_\alpha)=p(\overline{C}_\alpha)$. We
construct a definable homeomorphism $\theta_\alpha: \overline{D}_\alpha\to
\overline{C}_\alpha$ such that $p\circ\theta_\alpha=p|_{\overline{D}_\alpha}$. If
$\overline{C}_\alpha$ is the graph of $\overline{\zeta}_{\lambda,\mu}$, we must have
$\theta_\alpha(x',x_n)=(x',\overline{\zeta}_{\lambda,\mu}(x'))$. If $\overline{C}_\alpha$
is the closed band $[\overline{\zeta}_{\lambda,\mu},\overline{\zeta}_{\lambda,\mu+1}]$, we
choose $\theta_\alpha$ so that, for every $x'\in\overline{\sigma}_\lambda$, it carries
the segment $(\{x'\}\times R)\cap \overline{D}_\alpha$ onto the segment $\{x'\}\times
[\overline{\zeta}_{\lambda,\mu}(x'),\overline{\zeta}_{\lambda,\mu+1}(x')]$ in an affine
way (note that $(\{x'\}\times R)\cap \overline{D}_\alpha$ is reduced to a point if and
only if $\overline{\zeta}_{\lambda,\mu}(x')=\overline{\zeta}_{\lambda,\mu+1}(x')$). The
homeomorphism $\theta_\alpha$ carries $\partial D_\alpha$ onto $\partial C_\alpha$.

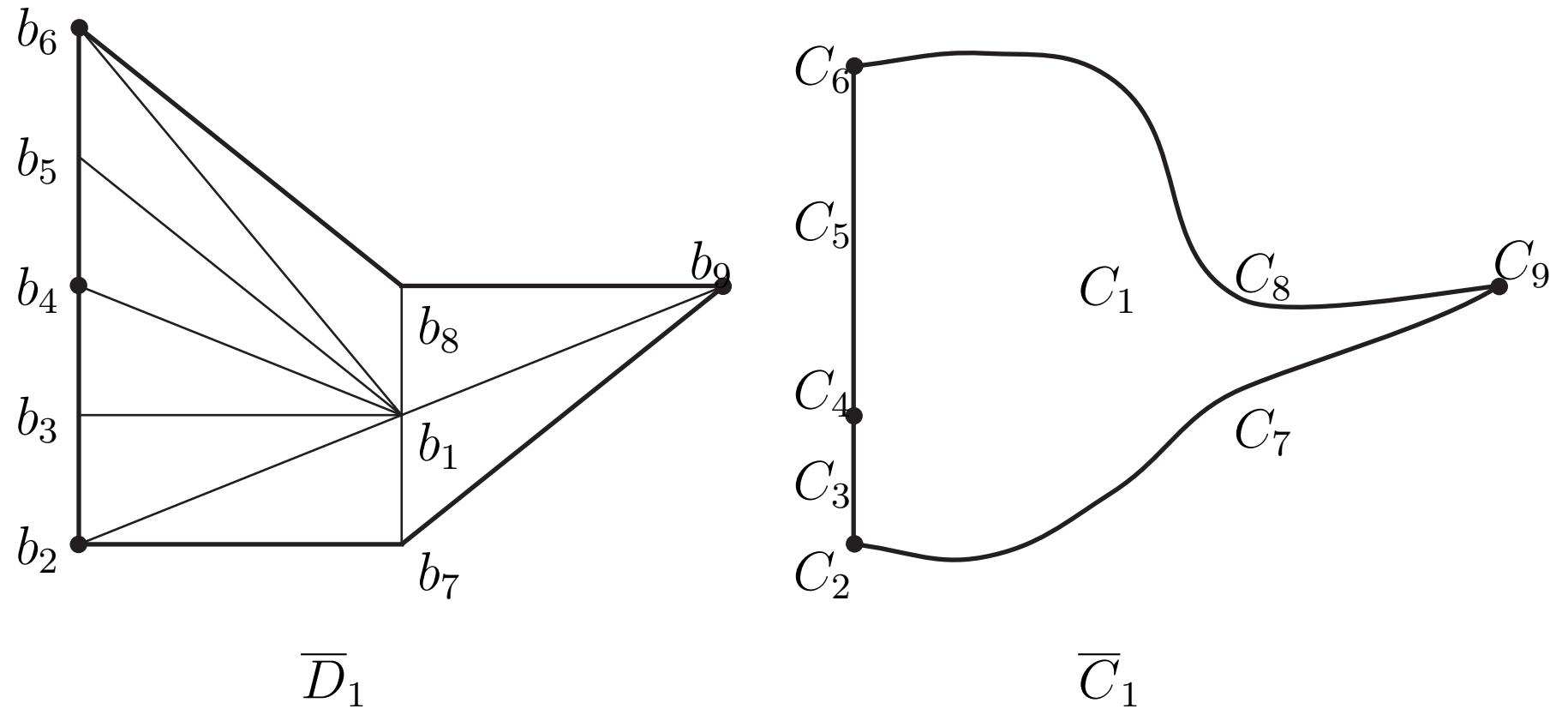
\begin{figure}[h]
\setlength{\unitlength}{1 mm}
\begin{picture}(130,60)
\put(12,10){\illustration{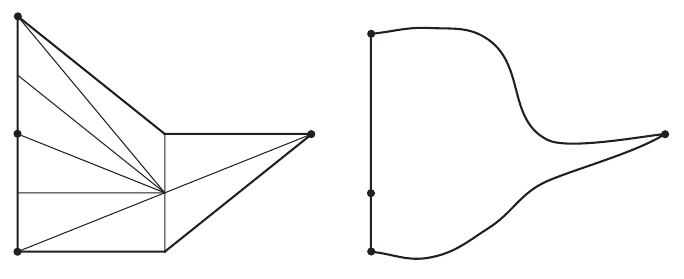}}
\put(8,12){\makebox(0,0)[l]{$b_2$}}
\put(8,22){\makebox(0,0)[l]{$b_3$}}
\put(8,32){\makebox(0,0)[l]{$b_4$}}
\put(8,42){\makebox(0,0)[l]{$b_5$}}
\put(8,52){\makebox(0,0)[l]{$b_6$}}
\put(39,10){\makebox(0,0)[l]{$b_7$}}
\put(39,20){\makebox(0,0)[l]{$b_1$}}
\put(39,29){\makebox(0,0)[l]{$b_8$}}
\put(60,34){\makebox(0,0)[l]{$b_9$}}
\put(68,10){\makebox(0,0)[l]{$C_2$}}
\put(68,17){\makebox(0,0)[l]{$C_3$}}
\put(68,24){\makebox(0,0)[l]{$C_4$}}
\put(68,37){\makebox(0,0)[l]{$C_5$}}
\put(68,49){\makebox(0,0)[l]{$C_6$}}
\put(102,21){\makebox(0,0)[l]{$C_7$}}
\put(90,32){\makebox(0,0)[l]{$C_1$}}
\put(102,33){\makebox(0,0)[l]{$C_8$}}
\put(122,34){\makebox(0,0)[l]{$C_9$}}
\put(30,2){\makebox(0,0)[l]{$\overline{D}_1$}}
\put(90,2){\makebox(0,0)[l]{$\overline{C}_1$}}
\end{picture}
\caption{\label{Fig1} The construction of $K$}
\end{figure}

\par  We finish the proof by constructing 
 a definable homeomorphism
$\Phi:|K|\to A$ such that $p\circ\Phi=p_{|K|}$ and
$\Phi(\overline{D}_\alpha)=\overline{C}_\alpha$ for $\alpha=1,\ldots,\ell$. We
cannot just take $\Phi|_{\overline{D}_\alpha}=\theta_\alpha$, because
$\theta_\alpha|_{\overline{D}_\beta}$ may be different from $\theta_\beta$ for
$\overline{D}_\beta\subset \partial{D}_\alpha$ (see what happens for
$\overline{D}_3\subset
\partial{D}_1$ in Figure \ref{Fig1}). So we modify $\theta_\alpha$ to obtain a
definable homeomorphism
$\Phi_\alpha: \overline{D}_\alpha\to \overline{C}_\alpha$ verifying
$\Phi_\alpha|_{\overline{D}_\beta}=\Phi_\beta$ and $p\circ\Phi_\alpha
=p|_{\overline{D}_\alpha}$. We proceed by induction on
$\dim(\overline{D}_\alpha)$. We take
$\Phi_\alpha=\theta_\alpha$ if $\overline{D}_\alpha$ is a point. If
$\dim(\overline{D}_\alpha)>0$, denote by $\rho_\alpha: \partial D_\alpha\to \partial
D_\alpha$ the homeomorphism defined by
$\rho_\alpha|_{\overline{D}_\beta}=\theta_\alpha^{-1}\circ\Phi_\beta$ for all
$\overline{D}_\beta\subset \partial D_\alpha$. We extend the homeomorphism
$\rho_\alpha$ to a homeomorphism $\eta_\alpha: \overline{D}_\alpha\to
\overline{D}_\alpha$ using the conic structure of $\overline{D}_\alpha=b_\alpha *
\partial D_\alpha$: for every
$x\in
\partial D_\alpha$ and $t\in[0,1]$, we set
$\eta_\alpha(tx+(1-t)b_\alpha)=t\rho_\alpha(x)+(1-t)b_\alpha$. Then we set
$\Phi_\alpha=\theta_\alpha\circ\eta_\alpha$. Observe that $p\circ
\eta=p|_{\overline{D}_\alpha}$ and $p\circ\Phi_\alpha=p|_{\overline{D}_\alpha}$. Now we
can take $\Phi$ defined by $\Phi_{\overline{D}_\alpha}=\Phi_\alpha$ for
$\alpha=1,\ldots,\ell$.
\end{proof}

\section{Triangulation of Definable Functions}

\begin{theorem}\label{triangfunctions} Let $X$ be a closed and bounded definable subset
of $R^n$ and $f:X\to R$
 a continuous definable function. Then there exist a finite simplicial complex $K$ in
$R^{n+1}$ and a definable homeomorphism $\rho : |K|_R\to X$ such that
$f\circ\rho$ is an affine function on each simplex of $K$.\par Moreover, given finitely
many definable subsets $B_1,\ldots,B_k$ of $X$, we may choose the triangulation
$\rho:|K|_R\to X$ so that each $B_i$ is a union of images of open simplices of $K$.
\end{theorem}

We replace $X$ with the set  $A=\{(f(x),x)\in R\times R^n\;;\ x\in X\}$ which is
definably homeomorphic to $X$. This set is a closed and bounded definable subset of
$R^{n+1}$. Let
$\pi=p_{n+1,1}:R^{n+1}\to R$ be the projection on the first coordinate. In order to
prove Theorem
\ref{triangfunctions}, it is sufficient to construct a finite simplicial complex $K$ in
$R^{n+1}$ and a definable homeomorphism $\Phi: |K|\to A$ such that
$\pi\circ\Phi=\pi|_{|K|}$. Indeed, the composition of $\pi$ with the homeomorphism
$X\to A$ is $f$. So Theorem \ref{triangfunctions} follows from the next proposition.\par

\begin{proposition}\label{triang_1coord} Let $A$ be a closed and bounded definable
subset of $R\times R^n$ and let
$B_i$, $i=1,\ldots,k$, be definable subsets of $A$. Let $\pi: R\times R^n\to R$ be
the projection on the first coordinate. Then there exist a finite simplicial complex
$K$ with vertices in
$R\times {\Bbb Q}^n$ and a definable homeomorphism  $\Phi:
|K|\to A$  such that
$\pi\circ\Phi=\pi|_{|K|}$ and each $B_i$ is a union of images by
$\Phi$ of open simplices of $K$.
\end{proposition}

\begin{proof}
  We proceed by induction on $n$. The case of $n=0$ is obvious. We can
subdivide
$R$ with finitely many points $x_1<\ldots<x_p$ such that $A$ and the $B_i$ are
unions of points $x_i$ and intervals $(x_j,x_{j+1})$. We choose for $K$ the
collection of points $x_i$ and closed and bounded intervals $[x_j,x_{j+1}]$
contained in $A$.\par
Assume that $n>0$ and the proposition is proved for $n-1$. As in the proof of Theorem
\ref{triangsets}, we can assume that all $B_i$ are closed.\par
 Let $G_0$ be the boundary of $A$ and $G_i$ the boundary of $B_i$, for $i=1,\ldots,k$.
Set $G=\bigcup_{i=0}^k G_i$. Then
$G$ is a closed and bounded definable set of dimension at most $n$. Denote by
$C$ the finite set of points $c\in R$ such that $\{x\in R^n\;;\  (c,x)\in G\}$ is of
dimension $n$. Let $F_i$, $i=0,\ldots,k$ be the union of the closure of
$G_i\setminus (C\times R^n)$ with the boundary of $G_i\cap (C\times R^n)$ in $C\times
R^n$. Set
$F=\bigcup_{i=0}^k F_i$. Each $F_i$ is a closed and bounded definable set. We claim
that, for every
$t$ in $R$, the dimension of $F_t=\{x\in R^n\;;\ (t,x)\in F\}$ is $<n$. Since $G_i\cap
(C\times R^n)$ has dimension $n$, the dimension of its boundary is at most $n-1$.
Hence, it suffices to check that the dimension of $(\clos(G_i\setminus (C\times
R^n)))_t$ is not greater than $n-1$. This follows from the next lemma applied to
$X=G_i\setminus (C\times R^n)$.  

\begin{lemma}\label{fam1param} Let $X$ be a definable subset of $R\times R^n$ such that,
for every
$t\in R$, the dimension of $X_t$ is $<n$. Then, for every $t\in R$, we have
$\dim((\clos(X))_t)<n$.
\end{lemma}

\begin{proof} It follows
from the assumption that $\dim X\leq n$. Hence, the dimension of the boundary of $X$ is
$< n$. Since $(\clos(X))_t \setminus X_t$ is contained in the boundary of $X$, we have
$\dim((\clos(X))_t)<n$.
\end{proof}
 
We return to the proof of Proposition \ref{triang_1coord}. Let 
$p=p_{n+1,n}:R\times R^n\to R\times R^{n-1}$ be the projection  on the first $n$
coordinates. By Proposition \ref{goodcoord}, since $\dim(F_t)<n$ for every $t\in R$, we
can assume that for all $(t,x') \in R\times R^{n-1}$, $\left(p^{-1}(t,x')\cap F\right)$
is finite.  We choose a cdcd of $R\times R^n$ adapted to   $F_0,\ldots,F_k$ and to
$\{c\}\times R^n$ for every $c\in C$. We get a finite partition of
$p(A)$ into definably connected definable subsets
$X_\lambda$ of $R\times R^{n-1}$, and definable continuous functions
$$\zeta_{\lambda,1}<\ldots<\zeta_{\lambda,m_\lambda} : X_\lambda\longrightarrow R\;,$$
such that every graph of $\zeta_{\lambda,\mu}$ is contained in some
$F_i$ and every
$F_i$ is a union of graphs of  $\zeta_{\lambda,\mu}$. Note that every graph of 
$\zeta_{\lambda,\mu}$ and every band 
$(\zeta_{\lambda,\mu},\zeta_{\lambda,\mu+1})$ is either contained in $G_i$
or disjoint from $G_i$. It follows that $A$ and the
$B_i$ are unions of such cells.\par
 Applying the inductive assumption, we may assume that there is a simplicial complex
$L$ with vertices in
$R\times{\Bbb Q}^{n-1}$ and a definable homeomorphism $\Psi:|L|\to p(A)$ such that
$p_{n,1}\circ\Psi=p_{n,1}|_{|L|}$ and all
$X_\lambda$ are unions of images of open simplices of 
$L$ by $\Psi$. From now on, we can follow the proof of Theorem \ref{triangsets}
and construct a triangulation $\Phi:|K|\to A$ such that all vertices of $K$ are
in $R\times{\Bbb Q}^n$, each $B_i$ is the image by $\Phi$ of a union of open
simplices of $K$, and $p\circ\Phi=\Psi\circ p|_{|K|}$. It follows that $$\pi\circ\Phi=
p_{n,1}\circ p\circ \Phi= p_{n,1}\circ \Psi\circ p|_{|K|}=p_{n,1}\circ
p|_{|K|}=\pi|_{|K|}\;.$$
\end{proof}

The vertices of the simplicial complex constructed in Proposition \ref{triang_1coord}
have all their coordinates rational but the first. We shall need to have also the first
coordinate rational. But this has to be paid with some extra complication. 

\begin{proposition}\label{triang_1coordrat} Let $A$ be a closed and bounded definable
subset of $R\times R^n$ and let
$B_i$, $i=1,\ldots,k$, be definable subsets of $A$. Let $\pi: R\times R^n\to R$ be the
projection on the first factor. Then there exist a finite simplicial complex
$K$ with vertices in
${\Bbb Q}\times {\Bbb Q}^n$ and definable homeomorphisms  $\Phi: |K|\to A$ and
$\tau:R\to R$  such that
$\tau\circ\pi\circ\Phi=\pi|_{|K|}$ and each $B_i$ is a union of images by
$\Phi$ of open simplices of $K$.
\end{proposition}

\begin{proof} We  just modify the step $n=0$ of the proof of the preceding proposition.
We can subdivide
$R$ with finitely many points $x_1<\ldots<x_p$ such that $A$ and the $B_i$ are unions of
points $x_i$ and intervals $(x_j,x_{j+1})$. We choose for $K$ the collection of points
$1,\ldots,p$ and segments $[j,j+1]$ such that $[x_j,x_{j+1}]$ is contained in $A$. We
choose for $\tau$ a piecewise affine homeomorphism sending $x_i$ to $i$. The rest of the
proof is the same, using the appropriate inductive assumption.
\end{proof}

The result of triangulation of continuous definable functions  cannot be
generalized to all continuous definable maps. Consider for instance, the map
$f : [0,1]^2\to {\Bbb R}^2$ defined by $f(x,y)=(x,xy)$. There is no way to choose
triangulations $\Phi: |K|\to [0,1]^2$ and $\Psi : |L|\to f([0,1]^2)$ such that
$\Psi^{-1}\circ f\circ \Phi$ is affine on every simplex of $K$. 

\debexo Prove the preceding assertion. Hint: Setting $a=(0,0)\in R^2$ and
$A=f((0,1)^2)$, remark that $a\in\clos(A)$, $\dim(\clos(f^{-1}(A))\cap f^{-1}(a))=1$
and $\dim(f^{-1}(x))<1$ for all $x\in A$. Show that this cannot happen for a map
 which is defined on a polyhedron and affine on each simplex of a simplicial
decomposition.
\finexo

Note that the crucial point in the proof of Theorem \ref{triangfunctions} which is only
valid for functions with values in $R$ is the claim that $\dim (F_t)$ is $<n$ for every
$t$ in $R$. Indeed, Lemma \ref{fam1param} is no longer true if the parameter $t$
varies in $R^m$ with $m>1$. The counterexample above also gives a counterexample to the
lemma for $m>1$: take 
$$X=\{(x,z,y)\in R^2\times R\;;\ x\not = 0 \et z=xy\}\;.$$
Then, for all $(x,z)\in R^2$, $\dim(X_{(x,z)})<1$, but $\dim(\clos(X)_{(0,0)})=1$.
\par
\medskip
We can deduce from the triangulation of definable functions a result which is usually
proved in other ways (as a consequence of Hardt's theorem in \cite{vdD}). If $a\in R^n$
and $r>0$, we denote by $S(a,r)$ the sphere with center $a$ and radius $r$.

\begin{theorem}[Local Conic Structure]\label{localconic} Let $A\subset R^n$ be a
closed definable set, $a$ a point in $A$. There is $r>0$ such that there exists a
definable homeomorphism $h$ from the cone with vertex $a$ and base $S(a,r)\cap
A$ onto $\overline{B}(a,r)\cap A$, satisfying $h|_{S(a,r)\cap A}={\rm Id}$ and
$\Vert h(x)-a\Vert=\Vert x-a \Vert$ for all $x$ in the cone.
\end{theorem}

\begin{proof} We triangulate the function $f:x\mapsto \Vert
x-a\Vert$ restricted to $A\cap \overline{B}(a,1)$. We obtain a definable homeomorphism
$\Phi : |K|\to A\cap \overline{B}(a,1)$ such that $f\circ \Phi$ is affine on each
simplex of $K$. Since $f$ takes its minimum in $a$, this point $a$ is the image of 
vertex $w$ of $K$. Let $\mu>0$ be the minimum of the values of $f\circ \Phi$ at all
other vertices of $K$, and take $r$ such that $0<r<\mu$. Take a point $y\in |K|$ such
that
$f\circ \Phi(y)=r$, i.e.\ $\Phi(y)\in S(a,r)\cap A$. The point $y$ belongs to a simplex
$[v_0,\ldots,v_d]$ of
$K$. We have $y=\sum_{i=0}^d \lambda_iv_i$ with $\sum_{i=0}^d \lambda_i=1$. Since $f\circ
\Phi(y)=\sum_{i=0}^d \lambda_if\circ\Phi(v_i)$, one of the $v_i$, say $v_0$, must be
$w$. Define $h$ on  the cone with vertex $a$ and base $S(a,r)\cap A$ by
$h(ta+(1-t)\Phi(y))=\Phi(tw+(1-t)y)$. It is easily checked that $h$ has the properties
stated in the theorem.
\end{proof}

The following exercise is a generalization of the local conic structure theorem.

\debexo 
Prove the following result:\\ 
{\it Let $Z\subset S$ be two closed and bounded
definable sets. Let $f$ be a nonnegative continuous definable function on $S$
such that $f^{-1}(0)=Z$. Then there are $\delta>0$ and a continuous definable
map  $h: f^{-1}(\delta)\times [0,\delta]\to f^{-1}([0,\delta])$, such that
$f(h(x,t))=t$ for every $(x,t)\in f^{-1}(\delta)\times [0,\delta]$,
$h(x,\delta)=x$ for every $x\in f^{-1}(\delta)$, and
$h|_{f^{-1}(\delta)\times\;]0,\delta]}$ is a homeomorphism onto $f^{-1}(]0,\delta])$.
}\\
Hints: 1) Triangulate $f$, so that one can assume $S=|K|$ for a finite simplicial
complex, $Z$ is a union of simplices and $f$ is affine on each simplex.\\
2) Choose $\delta>0$ so small that, for every vertex $a$ of
$K$ such that $a\not\in Z$, $\delta<f(a)$.\\
3) Let $x\in f^{-1}(\delta)$. The point $x$ belongs to a simplex $[a_0,\ldots,a_d]$ of
$K$. We can assume that $a_i\in Z$, for
$i=0,\ldots,k$, and $a_i\not\in Z$, for $i=k+1,\ldots,d$. Let
 $x=\sum_{i=0}^d\lambda_i a_i$ with $\sum_{i=0}^d\lambda_i=1$. Show that
$\alpha=\sum_{i=0}^k
\lambda_i$ satisfies $0<\alpha<1$.\\
4) For $t\in[0,\delta]$, define $h(x,t)$ as the point $y$ in the segment joining
$\sum_{i=0}^k (\lambda_i/\alpha) a_i$ to $x$ such that $f(y)=t$. Check that $h$ is
well defined and satisfies the required properties. 
\finexo

%% file: OMINcinq.tex
\chapter{Generic Fibers for Definable Families}

\section{The Program}

We motivate what we are going to do with a very simple example from
algebraic geometry. Consider the family of conics
$x^2-t(1+y^2)=0$ in ${\Bbb C}^2$, parametrized by $t\in{\Bbb C}$. We
can also regard $x^2-t(1+y^2)=0$ as the equation of a conic defined
over the field ${\Bbb C}(t)$. This conic is the ``generic fiber of
the family''. The fact that the generic fiber of the family is a
nondegenerate conic (the discriminant is $t^2$, which is invertible
in ${\Bbb C}(t)$) corresponds to the fact that the conics in the
family are non degenerate for almost all values of $t$ (the only
exception is $t=0$). This is formalized by using the prime spectrum
${\rm Spec}({\Bbb C}[t])$ of the ring ${\Bbb C}[t]$. We have ${\Bbb
C}$ embedded in ${\rm Spec}({\Bbb C}[t])$, each $z\in {\Bbb C}$
corresponding to the maximal ideal $(t-z)$ of ${\Bbb C}[t]$. The
residue field at these maximal ideals is ${\Bbb C}$. There is another
point in ${\rm Spec}({\Bbb C}[t])$: the ideal $(0)$, whose residue
field is ${\Bbb C}(t)$. This point is the ``generic point'' of ${\Bbb
C}$.\par
We are going to embed $R^m$ into a bigger space $\widetilde{R^m}$. We
shall associate to each point $\alpha\in\widetilde{R^m}$ a real
closed field $\kappa(\alpha)$ with an o-minimal structure. Then, given
a definable family $X\subset R^m\times R^n$ and
$\alpha\in\widetilde{R^m}$, we shall define the fiber
$X_\alpha$ which will be a definable subset of $\kappa(\alpha)^n$. If
$\alpha$ is the image by the imbedding of a point $t\in R^m$, then
$\kappa(\alpha)$ will be $R$ and $X_\alpha$ will be the usual fiber
$X_t$. The interesting things will happen for ``generic fibers''
$X_\alpha$ corresponding to points $\alpha\in\widetilde{R^m}\setminus
R^m$. We shall relate the properties of the generic fiber
$X_\alpha$ with properties of the family $X$ over ``large'' subsets
of the parameter space $R^m$. This will be used to deduce results
of triviality for definable families from the triangulation
theorems of the preceding chapter.\par 
The tools that we are going to
introduce are actually a reformulation of classical model-theoretic
notions and results ($m$-types, definable ultrapower,\dots). This
reformulation is modelled upon the theory of the real spectrum (cf.
\cite{BCR}).

\section{The Space of Ultrafilters of Definable Sets}

Let $S_m$ be the Boolean algebra of definable subsets of $R^m$.
The content of this section is nothing but the construction of the
Stone space of the Boolean algebra $S_m$ (cf.\ \cite{BS}).
Nevertheless, we give a more or less self-contained account of this
construction for the reader who is not familiar with ultrafilters
and Stone spaces.\par 

We denote by $\widetilde{R^m}$ the set of ultrafilters of ${\cal
S}_m$. Recall that an {\em ultrafilter of $S_m$}\index{ultrafilter}
is a subset
$\alpha$ of
$S_m$ such that
\begin{enumerate}
\item $R^m\in \alpha$
\item $A\cap B\in \alpha$ if and only if $A\in \alpha$ and $B\in
\alpha$
\item $\emptyset\not\in\alpha$ 
\item $A\cup B\in \alpha$ if and only if $A\in \alpha$ or $B\in
\alpha$
\end{enumerate}
We say that a family ${\cal F}$ of definable subsets of $R^m$
 {\em generates the ultrafilter
$\alpha$} if $\alpha$ is the set of definable sets $A\subset R^m$ such
that there exists $B\in {\cal F}$ with $B\subset A$. If ${\cal F}$
is a nonempty family of definable subsets of $R^m$, closed under
finite intersections, then it generates an ultrafilter if and only
if $\emptyset\not\in{\cal F}$ and, for every $A\in S_m$, either $A$ or
$R^m\setminus A$ contains a $B\in {\cal F}$.\par 
If $t$ is a point of $R^m$, the definable subsets of $R^m$ containing
$t$ form an ultrafilter $\alpha_t$. This is called the {\em principal
ultrafilter}\index{ultrafilter!principal --} generated by $t$. The map
$t\mapsto
\alpha_t$ embeds
$R^m$ as a subset of $\widetilde{R^m}$. The following exercises give
examples of points of $\widetilde{R^m}\setminus R^m$

\debexo Points of $\widetilde R$. Show that that the following
families generate ultrafilters of $S_1$:\\
1) The family of all intervals $(x,+\infty)$ for $x\in R$.\\
2) The family of all intervals $(-\infty,x)$ for $x\in R$.\\
3) For a fixed $a\in R$, the family of all intervals
$(a,a+\varepsilon)$ for $\varepsilon>0$.\\
4) For a fixed $a\in R$, the family of all intervals
$(a-\varepsilon,a)$ for $\varepsilon>0$.\\
If $R={\Bbb R}$, show that all points of $\widetilde{{\Bbb
R}}\setminus {\Bbb R}$ correspond to one of the four cases
above.\finexo

\debexo\label{alphaup} Construction of points of $\widetilde{R^m}$ by
induction.  Let $\alpha$ be a point of $\widetilde{R^m}$. Show that the
family of definable subsets of $R^{m+1}$ of the form 
$$\{(x,y)\in R^m\times R\;;\ x\in A \et 0<y<f(x)\}\;,$$
where $A\in\alpha$ and $f:A\to R$ is a positive definable function,
generates an ultrafilter $\alpha_\uparrow$ of $S_{m+1}$. Hint:
given a definable subset $B$ of $R^{m+1}$, use a cdcd of $R^{m+1}$ to
show that either $B$ or its complement belongs to $\alpha_\uparrow$.
\finexo

\debexo Dimension of a point of $\widetilde{R^m}$. For
$\alpha\in\widetilde{R^m}$, define $\dim \alpha$ as the minimum of
$\dim A$ for $A\in \alpha$. Show that $\dim \alpha=d$ if and only if
$\alpha$ is generated by a family of definable sets of dimension $d$.
Show that $\dim(\alpha_\uparrow)=\dim\alpha +1$ (with the notation
of Exercise \ref{alphaup}). Show that there is an
ultrafilter of dimension $m$ in $\widetilde{R^m}$.  Show that
$\dim A$ is the maximum of
$\dim \alpha$ for $\alpha\ni A$.
\finexo

Now we put a topology on $\widetilde{R^m}$. We consider
$\widetilde{R^m}$ as a subset of the powerset $2^{S_m}$, identifying
an ultrafilter $\alpha$ with its characteristic function ${\bf
1}_\alpha:S_m\to 2=\{0,1\}$. We equip $2$ with the discrete topology
and $2^{S_m}$ with the product topology. By Tychonoff theorem,
$2^{S_m}$ is compact Hausdorf. We take as topology on
$\widetilde{R^m}$ the topology induced by the topology of $2^{S_m}$.

\begin{proposition}\label{spectrumcompact}
The topology of $\widetilde{R^m}$ has a basis of open and closed sets
consisting of all
$\widetilde{A}=\{ \alpha\in \widetilde{R^m}\;;\ \alpha\ni A\}$ for
$A\in S_m$.  The space $\widetilde{R^m}$ is compact Hausdorf.
\end{proposition}

\begin{proof} First note that the operation
$A\mapsto
\widetilde A$ preserves finite unions, finite intersections and
taking complement.\par
The product topology of
$2^{S_m}$ induces on
$\widetilde{R^m}$ the topology which has as basis of open sets the
sets of the form
$$U=\{\alpha\in \widetilde{R^m}\;;\ {\bf 1}_\alpha(A_1)=\ldots={\bf
1}_\alpha(A_k)=1 \et {\bf 1}_\alpha(B_1)=\ldots={\bf
1}_\alpha(B_\ell)=0\}\;,$$
where $(A_1,\ldots,A_k)$ and $(B_1,\ldots,B_\ell)$ are finite
families of elements of $S_m$. Note that $U=\widetilde A$, where
$A=\bigcap_{i=1}^k A_i\cap \bigcap_{j=1}^\ell(R^m\setminus B_j)$.
This proves the first assertion of the proposition.\par
Finally, note that each one of the four conditions of the definiton
of an ultrafilter defines a closed subset of $2^{S_m}$. Hence,
$\widetilde{R^m}$ is closed in $2^{S_m}$. This proves the second
assertion.
\end{proof}

\debexo Prove that the application $A\mapsto \widetilde A$
 is a bijection from $S_m$ to the set of open and
closed subsets of $\widetilde{R^m}$.
\finexo  

\section[The O-minimal Structure $\kappa(\alpha)$]{O-minimal Structure
Associated with an Ultrafilter}

If $A$ is a definable subset of $R^m$, we denote by ${\rm Def}(A,R)$
the ring of definable functions $A\to R$. Given $\alpha\in
\widetilde{R^m}$, we define $\kappa(\alpha)$ as the inductive limit
of ${\rm Def}(A,R)$ for $A\in \alpha$. This means the following. We
form the disjoint union $\bigsqcup_{A\in \alpha}{\rm Def}(A,R)$. We
say that two elements $f:A\to R$ and $g:B\to R$ of this disjoint
union are equivalent if there exists $C\in \alpha$, $C\subset A\cap
B$, such that $f|_C=g|_C$. Then $\kappa(\alpha)$ is the set of
equivalence classes for this equivalence relation. We denote by
$f(\alpha)\in \kappa(\alpha)$ the equivalence class of $f: A\to R$.
By definition, $f(\alpha)=g(\alpha)$ if and only if $f$ and $g$
coincide on a definable set belonging to $\alpha$. The inductive limit
$\kappa(\alpha)$ has a canonical structure of commutative
$R$-algebra:
The sum $f(\alpha)+g(\alpha)$ is $(f+g)(\alpha)$, where $f+g$ is
defined on the intersection of the domains of $f$ and $g$, which
belongs to
$\alpha$. We have a similar definition for the
product $f(\alpha)g(\alpha)$. The images of elements of $R$ are the
classes of the constant functions. The verifications are easy.\par
We define $f(\alpha)\in\kappa(\alpha)$ to be positive if $f$ is
positive on a definable set belonging to $\alpha$. This does not
depend on the choice of the representant $f$.

\begin{proposition} The commutative $R$-algebra $\kappa(\alpha)$, with
positive elements defined as above, is an ordered field which is an
ordered extension of $R$.
\end{proposition}

\begin{proof} We check that, for an element $f(\alpha)$ of
$\kappa(\alpha)$, we have exactly one of the three possibilities
$f(\alpha)>0$, $f(\alpha)=0$ and $-f(\alpha)>0$. This is because the
domain of $f$ is partitioned into the three definable sets where $f$
is positive, zero or negative, respectively. Exactly one of these
three sets belongs to $\alpha$. Moreover, if $f(\alpha)\not=0$, then
$1/f$ is defined on a definable set belonging to $\alpha$. Hence,
$(1/f)(\alpha)$ is the inverse of $f(\alpha)$ in $\kappa(\alpha)$.
The other verifications are also easy.
\end{proof}

Now we shall construct an o-minimal structure on $\kappa(\alpha)$.
For this we need the notion of the fiber of a definable family
$X\subset R^m\times R^n$ at $\alpha\in \widetilde{R^m}$. If
$f=(f_1,\ldots,f_n): A\to R^n$ is a definable map and $A\in \alpha$,
we denote by $f(\alpha)$ the point
$(f_1(\alpha),\ldots,f_n(\alpha))\in
\kappa(\alpha)^n$.

\begin{definition} Let $X$ be a definable subset of $R^m\times R^n$.
The {\em fiber of $X$ at $\alpha\in\widetilde{R^m}$}\index{fiber} is
the set
$X_\alpha$ of those $f(\alpha)$ in
$\kappa(\alpha)^n$ such that there exist $A\in \alpha$ on which $f$
is defined and $(t,f(t))\in X$ for all $t\in A$.
\end{definition}

 In other words,
$f(\alpha)$ belongs to $X_\alpha$ if and only if there is $A\in
\alpha$ such that $f(t)\in X_t$ for all $t\in A$. This definition
makes sense because, if we take another representant $g$ of
$f(\alpha)$ (i.e.\ $g(\alpha)=f(\alpha)$), $g$ and $f$ coincide on
some $B\in \alpha$ and, therefore,  $g(t)\in X_t$ for all $t\in A\cap
B$, with $A\cap B\in\alpha$.\par
\medskip

It is worth looking at what the preceding constructions and
definitions give in the case of a principal ultrafilter $\alpha_t$
corresponding to a point $t\in R^m$. Since the ultrafilter $\alpha_t$
is generated by the singleton $\{t\}$, the inductive limit
$\kappa(\alpha_t)$ is (canonically isomorphic to) ${\rm
Def}(\{t\},R)=R$. Moreover, if $f:A\to R$ is a definable function and
$A\ni t$, then $f(\alpha_t)=f(t)\in R$. Finally, given a definable
set $X\subset R^m\times R^n$, its fiber $X_{\alpha_t}$ is simply the
set of $x\in R^n$ such that $(t,x)\in X$, i.e.\ $X_t$.\par
\medskip

Now we can define the o-minimal structure on $\kappa(\alpha)$. Let
$S_n(\alpha)$ be the family of all fibers $X_\alpha$, for $X$ a
definable subset of $R^m\times R^n$. 

\begin{theorem}\label{residualo-minstruct} The collection
$(S_n(\alpha))_{n\in{\Bbb N}}$ is an o-minimal structure expanding
the ordered field $\kappa(\alpha)$. In particular, $\kappa(\alpha)$
is real closed.
\end{theorem}

\begin{proof} We have to check the five properties in the definition of
an o-minimal structure.\par
\smallskip
1) {\it For every $n\in{\Bbb N}$, $S_n(\alpha)$ is a Boolean subalgebra
of the powerset of $\kappa(\alpha)^n$.}\par
It is sufficient to check the following fact. Let $X$ and $Y$ be
definable subsets of
$R^m\times R^n$. Then we have:
$$\displaylines{\qquad (X\cap Y)_\alpha=X_\alpha\cap Y_\alpha\;,\qquad
 (X\cup Y)_\alpha=X_\alpha\cup Y_\alpha\;,\hspace*{\fill}\cr
\hspace*{\fill}((R^m\times R^n)\setminus X)_\alpha=
\kappa(\alpha)^n\setminus X_\alpha\;.\qquad\cr}$$
Let us check the last equality. Let $f(\alpha)$ be an element of
$\kappa(\alpha)^n$, with $f_i\in {\rm Def}(A,R)$ and $A\in \alpha$.
The two definable sets $B=\{ t\in A\;;\ f(t)\in X_t\}$ and $C=\{ t\in
A\;;\ f(t)\in ((R^m\times R^n)\setminus X)_t\}$ partition $A$.
Therefore, exactly one of $B$ and $C$ belongs to $\alpha$. This
means that $f(\alpha)$ belongs to $((R^m\times R^n)\setminus
X)_\alpha$ if and only if it does not belong to $X_\alpha$.\par
The verification of the other equalities is left to the reader.\par
\smallskip

2) {\it All algebraic subsets of $\kappa(\alpha)^n$ are in
$S_n(\alpha)$.}
It is sufficient to check this for an algebraic set $V\subset
\kappa(\alpha)^n$ given by one equation
$$\postdisplaypenalty=-10000 % page break of the original edition
\sum_{i\in I} a_i(\alpha) x^i=0\;,$$
where $I$ is a finite subset of ${\Bbb N}^n$, $i=(i_1,\ldots,i_n)$ and
$x^i=x_1^{i_1}\cdots x_n^{i_n}$. Let $A\in \alpha$ be such that all
definable functions $a_i$ are defined on $A$. Set 
$$X= \{(t,x) \in R^m\times R^n\;;\ t\in A \et \sum_{i\in I}
a_i(t)x^i=0\}\;.$$
Then $X$ is a definable set. We have $f(\alpha)\in V$ if and only if
there exists $B\in \alpha$ such that $\sum_{i\in I} a_i(t)
f^i(t)=0$ for all $t\in B$, i.e.\ $f(t)\in X_t$ for all $t\in B$.
Hence, $V=X_\alpha$.\par
\smallskip

3) {\it If $X_\alpha\in S_n(\alpha)$ and $Y_\alpha\in S_p(\alpha)$,
then $X_\alpha\times Y_\alpha\in S_{n+p}(\alpha)$.}
Set
$$X\times_{R^m} Y= \{(t,x,y)\in R^m\times R^n \times R^p\;;\ (t,x)\in
X \et (t,y)\in Y\}\;.$$
Note that $(X\times_{R^m} Y)_t= X_t\times Y_t$ for all $t\in R^m$. It
follows that $(X\times_{R^m} Y)_\alpha= X_\alpha\times Y_\alpha$.
\par\smallskip

4) {\it If $p:\kappa(\alpha)^{n+1}\to \kappa(\alpha)^n$ is the
projection on the first
$n$ coordinates and $X_\alpha\in S_{n+1}(\alpha)$, then
$p(X_\alpha)\in S_n(\alpha)$.}
Let $\pi : R^m\times R^{n+1}\to R^m\times R^n$ be the projection on
the first $m+n$ coordinates. It is sufficient to check that
$$p(X_\alpha) = \pi(X)_\alpha\;.$$
Let $f(\alpha)$ be an element of $\kappa(\alpha)^n$. If $f(\alpha)$
belongs to $p(X_\alpha)$, there is $g(\alpha)\in \kappa(\alpha)$ and
$A\in \alpha$ such that $(f(t),g(t))\in X_t$ for all $t\in A$. Hence
$f(t)\in \pi(X)_t$ for all $t\in A$, which shows $f(\alpha)\in
\pi(X)_\alpha$. Conversely, assume $f(\alpha)\in
\pi(X)_\alpha$. Then there is $B\in \alpha$ such that, for all $t\in
B$, there exists $y\in R$ with $(f(t),y)\in X_t$. By Definable Choice
\ref{defchoice}, we can choose such a $y$ as a definable function
$g:B\to R$. It follows that $(f(\alpha),g(\alpha))\in X_\alpha$ and
$f(\alpha)\in p(X_\alpha)$.\par 
\smallskip 

5) {\it The elements of $S_1(\alpha)$ are precisely the finite unions
of points and intervals in $\kappa(\alpha)$.} Take $X_\alpha\in
S_1(\alpha)$ and choose a cdcd of $R^m\times R$ adapted to $X$. We
have a partition of $R^m$ into cells $C_1,\ldots,C_k$ and, for each
cell $C_i$, continuous definable functions
$\zeta_{i,1}<\ldots<\zeta_{i,\ell(i)}: C_i\to R$ such that $X$ is a
union of graphs of $\zeta_{i,j}$ and bands
$(\zeta_{i,j},\zeta_{i,j+1})$. Exactly one of the $C_i$ belongs to
the ultrafilter $\alpha$, say $C_1\in \alpha$. We claim that $X_\alpha$
is the union of the points $\zeta_{1,j}(\alpha)$ such that the graph of
$\zeta_{1,j}$ is contained in $X$ and the intervals
$(\zeta_{1,j}(\alpha),\zeta_{1,j+1}(\alpha))$ such that the band
$(\zeta_{1,j},\zeta_{1,j+1})$ is contained in $X$. Since $C_1\in
\alpha$, we can replace $X$ with $X\cap (C_1\times R)$. Since taking
the fiber at $\alpha$ commutes with union, we can assume that
$X$ is a graph or a band over $C_1$. If $X$ is the band
$(\zeta_{1,j},\zeta_{1,j+1})$, then $f(\alpha)\in X_\alpha$ if and
only if there is $B\in \alpha$ such that
$\zeta_{1,j}(t)<f(t)<\zeta_{1,j+1}(t)$ for all $t\in B$, i.e.\
$f(\alpha)\in (\zeta_{1,j}(\alpha),\zeta_{1,j+1}(\alpha))$.
\end{proof}

Note that there are two crucial points in the proof: the use of
Definable Choice in the proof of the stability by projection (point
4) and the use of the cdcd to prove the o-minimality (point 5). All
the rest is rather formal. \par

The proof of Theorem \ref{residualo-minstruct} shows that taking
fibers at $\alpha$ commutes with the boolean operations and
the projections. We shall formalize this remark in order to have a
useful tool for translating properties of the fiber $X_\alpha$ to
properties of the fibers $X_t$ for all $t$ in a definable set
belonging to $\alpha$. First we introduce some notation.
We consider the {\em definable families of
formulas}\index{definable family!-- of formulas}
$\Phi_t(x)$, where $t$ ranges over $R^m$, which are constructed
according to the following rules:
\begin{enumerate}
\item If $X\subset R^m\times R^n$ is definable, $x\in X_t$ is
a definable family of formulas.
\item A polynomial equation $P(x)=0$ or inequality $P(x)>0$ with
coefficients in $R$ is a definable (constant) family of formulas.
\item If $\Phi_t(x)$ and $\Psi_t(x)$ are definable families of
formulas, then $\Phi_t(x) * \Psi_t(x)$ (where $*$ is one of ``and'',
``or'', $\Rightarrow$, $\Leftrightarrow$) and the negation ``not
$\Phi_t(x)$'' are definable families of formulas.
\item If $\Phi_t(x,y)$ is a definable family of formulas (with
$y$ ranging over $R^p$) and $Y\subset R^m\times R^p$ is definable,
the existential quantification $\exists y\in Y_t\ \Phi_t(x,y)$ and
the universal quantification $\forall y\in Y_t\ \Phi_t(x,y)$ are
definable families of formulas.
\end{enumerate}
Given a definable family of formulas $\Phi_t(x)$  and
$\alpha\in\widetilde{R^m}$, we define the {\em fiber
formula}\index{fiber!-- formula}
$\Phi_\alpha(x)$ in the following way:
\begin{enumerate}
\item If $\Phi_t(x)$ is $x\in X_t$, $\Phi_\alpha(x)$ is $x\in
X_\alpha$
\item If $\Phi_t(x)$ is the constant family $P(x)=0$ or 
$P(x)>0$, $\Phi_\alpha(x)$ is $P(x)=0$ or $P(x)>0$.
\item If $\Phi_t(x)$ is $\Theta_t(x) * \Psi_t(x)$, 
 or ``not
$\Theta_t(x)$'', $\Phi_\alpha(x)$ is $\Theta_\alpha(x) *
\Psi_\alpha(x)$ or ``not
$\Theta_\alpha(x)$'', respectively.
\item If $\Phi_t(x)$ is $\exists y\in Y_t\ \Psi_t(x,y)$ or $\forall
y\in Y_t\ \Psi_t(x,y)$, $\Phi_\alpha(x)$ is $\exists y\in Y_\alpha\
\Psi_\alpha(x,y)$ or $\forall y\in Y_\alpha\ \Psi_\alpha(x,y)$,
respectively.
\end{enumerate}
Note that, if $\Phi_\alpha(x)$ is a fiber formula with
$x=(x_1,\ldots,x_n)$, the set of $x\in \kappa(\alpha)^n$ such
that $\Phi_\alpha(x)$ holds is definable in the o-minimal structure
on $\kappa(\alpha)$.\par
Let us take an example. To the definable family of formulas
$$\forall x\in X_t\quad \exists \varepsilon >0\quad \forall y\in R^n
\quad(\Vert x-y\Vert < \varepsilon \Rightarrow y\in X_t)$$
corresponds the fiber formula
$$\forall x\in X_\alpha\quad \exists \varepsilon >0\quad \forall y\in
\kappa(\alpha)^n\quad (\Vert x-y\Vert < \varepsilon \Rightarrow y\in
X_\alpha)\;.$$
Note that $\Vert x-y\Vert < \varepsilon$ can be expressed as a
polynomial inequality.

\begin{proposition}\label{fibandfam} Let $X$ be a definable subset
of $R^m \times R^n$ and $\Phi_\alpha(x)$ the fiber formula of a
definable family of formulas $\Phi_t(x)$, with
$x=(x_1,\ldots,x_n)$. The equality
$$X_\alpha=\{ x\in \kappa(\alpha)^n\;;\ \Phi_\alpha(x)\}$$
holds if and only if there exists $A\in \alpha$ such that the
equality
$$X_t=\{ x\in R^n\;;\ \Phi_t(x)\}$$
holds for all $t\in A$. In particular, if $n=0$ (i.e. there is
no free variable), the fiber formula $\Phi_\alpha$ holds if and only
if there exists $A\in \alpha$ such that $\Phi_t$ holds for all $t\in
A$.
\end{proposition}

\begin{proof}
We proceed by induction on the construction of the definable family
of formulas according to rules 1-4. \par
Rule 1: $\Phi_t(x)$ is $x\in Y_t$. We have to show that
$X_\alpha=Y_\alpha$ if and only if there is $A\in \alpha$ such that
$X_t=Y_t$ for all $t\in A$. The ``if'' part is easy. To show the
``only if'' part, it suffices to consider the case $Y=\emptyset$,
since taking the fiber at $\alpha$ preserves the boolean
operations. Suppose that the set of $t$ such that $X_t=\emptyset$ does
not belong to $\alpha$. Then the projection of $X$ on the space $R^m$
of the first $m$ coordinates belongs to $\alpha$. By Definable Choice,
there is a definable map $f$ from this projection to $R^n$ whose
graph is contained in $X$. It follows that $f(\alpha)\in X_\alpha$,
which contradicts $X_\alpha=\emptyset$.\par
Rule 2. Note that we can actually omit this rule. Indeed, if 
$X_t=\{ x\in R^n\;;\ P(x) = 0\}$ for all $t\in R^m$, then
$X_\alpha=\{ x\in \kappa(\alpha)^n\;;\ P(x) = 0\}$ and we are reduced
to Rule 1 (the same with $>$ instead of $=$).\par
Rule 3. Here we use the observation that taking fibers at $\alpha$
commutes with the boolean operations.\par
Rule 4. Here we use the observation that taking fibers at $\alpha$
commutes with the projections (for existential quantification). The
universal quantification reduces to negations and existential
quantification.
\end{proof}

\debexo\label{closedfib} Let $X$ and $Y$ be definable subsets of
$R^m\times R^n$. Show that $X_\alpha$ is closed in
$Y_\alpha$ if and only if there exists $A\in\alpha$ such that
$X_t$ is closed in $Y_t$ for every $t\in A$.
\finexo

\section{Extension of Definable Sets}

Let $\alpha$ be an element of $\widetilde{R^m}$. We begin with the
extension of polyhedra from $R$ to $\kappa(\alpha)$. This will be
needed for the proof of Hardt's Theorem \ref{Hardt}.\par
Let
$a_0,\ldots,a_d$ be points of $R^n$ which are affine independent. We
denote by $[a_0,\ldots,a_d]_R$ the simplex they generate in $R^n$, and
by $[a_0,\ldots,a_d]_{\kappa(\alpha)}$ the simplex they generate in
$\kappa(\alpha)^n$. If $K$ is a finite simplicial complex with
vertices in $R^n$, we denote by $|K|_R$ its realization in $R^n$, and
by $|K|_{\kappa(\alpha)}$ its realization in $\kappa(\alpha)^n$, i.e.
the union of the simplices $[a_0,\ldots,a_d]_{\kappa(\alpha)}$ for all
$[a_0,\ldots,a_d]_R$ in $K$. If $V_R$ is a union of open simplices of
$K$, we define $V_{\kappa(\alpha)}$ as the union of the corresponding
open simplices in $\kappa(\alpha)^n$.

\begin{lemma}\label{extpoly} $|K|_{\kappa(\alpha)} = (R^m \times
|K|_R)_\alpha$. $V_{\kappa(\alpha)}= (R^m\times V_R)_\alpha$. 
\end{lemma}

\begin{proof} Since taking fibers at $\alpha$ preserves finite
unions, it is sufficient to prove the lemma for an open simplex
$\sigma$ with vertices $a_0,\ldots,a_d$. The open simplex $\sigma$ is
the set of $x\in R^n$ satisfying the formula
$$ \displaylines{\qquad\exists \lambda_1\in R\ldots\exists
\lambda_d\in R\ \Big(
\lambda_1>0 \et\ldots\et \lambda_d>0\et\hspace*{\fill}\cr
\hspace*{\fill}\sum_{i=1}^d\lambda_i=1\et
\sum_{i=1}^d \lambda_ia_i=x\Big)\;.\qquad\cr}$$ Hence, the fiber
$(R^m\times\sigma)_\alpha$ is described by the fiber formula of the
constant family of formulas:
$$\displaylines{\qquad\exists \lambda_1\in \kappa(\alpha)\ldots\exists
\lambda_d\in
\kappa(\alpha)\ \Big(
\lambda_1>0 \et\ldots\et
\lambda_d>0\et\hspace*{\fill}\cr
\hspace*{\fill}\sum_{i=1}^d\lambda_i=1\et
\sum_{i=1}^d \lambda_ia_i=x\Big)\;.\qquad\cr}$$ This fiber formula
describes
$\sigma_{\kappa(\alpha)}$.
\end{proof}

The notation for polyhedra will also be used for fibers of constant
definable families. Let $S$ be a definable subset of $R^n$. Then $R^m
\times S$ can be regarded as a constant definable family. We denote
by $S_{\kappa(\alpha)}\subset  \kappa(\alpha)^n$ the fiber
$(R^m\times S)_\alpha$, and we call $S_{\kappa(\alpha)}$ the {\em
extension of $S$}\index{extension of a definable set} to
$\kappa(\alpha)$. 

\debexo Show that $S_{\kappa(\alpha)}\cap R^n=S$.\finexo

In model-theoretic terms, the o-minimal structure over
$\kappa(\alpha)$ is an extension of the o-minimal stucture over $R$.
Any first-order formula of the o-minimal structure over $R$ can be
interpreted in the o-minimal structure over $\kappa(\alpha)$, taking
the extension to $\kappa(\alpha)$ of the definable sets appearing in
this formula. Proposition \ref{fibandfam} implies that
$\kappa(\alpha)$ is an elementary extension of $R$: every formula
without free variables holds over $R$ if an only if it holds over
$\kappa(\alpha)$.

\section{Definable Families of Maps}

Let $X$ and $Y$ be definable subsets of $R^m\times R^n$ and
$R^m\times R^p$, respectively. We regard them as definable families
parametrized by $R^m$. A {\em definable family of
maps}\index{definable family!-- of maps} from
$X$ to
$Y$ is a definable map $f:X\to Y$ such that the following diagram
commutes,
\setlength{\unitlength}{1mm}
$$\begin{picture}(40,23)
%je place le diagramme
\put(3,20){\makebox(0,0){
%je fais le diagramme
\begin{picture}(0,0)
\put(0,0){\makebox(0,0){$X$}}
\put(5,0){\vector(1,0){16}} \put(13,2){\makebox(0,0){$f$}}
\put(26,0){\makebox(0,0){$Y$}}
\put(0,-3){\vector(1,-1){11}}
\put(26,-3){\vector(-1,-1){11}}
\put(13,-17){\makebox(0,0){$R^m$}}
\end{picture}
%fin du diagramme 
}}
\end{picture}$$
where the maps to $R^m$ are the projections on the first $m$
coordinates. We obtain a family of maps $f_t:X_t\to Y_t$ for $t\in
R^m$ defined by $f(t,x)=(t,f_t(x))$. \par
Given $\alpha\in\widetilde{R^m}$, we shall define $f_\alpha:
X_\alpha\to Y_\alpha$, which will be a definable map for the
o-minimal structure on $\kappa(\alpha)$. Set 
$$\Gamma=\{ (t,x,y)\in R^m\times R^n\times R^p\;:\ (t,x)\in X \et
f(t,x)=(t,y)\}\;.$$
The set $\Gamma$ is a definable family parametrized by $R^m$ and, for
every $t\in R^m$, $\Gamma_t$ is the graph of $f_t:X_t\to Y_t$. Hence,
the formulas in the definable families of formulas
$$(*)_t\qquad\cases{
\forall x\in R^n\ \big((\exists y\in R^p\ (x,y)\in \Gamma_t)
\Leftrightarrow x\in X_t\big)\cr
\forall x\in R^n\ \forall y\in R^p\ \big((x,y)\in\Gamma_t \Rightarrow
y\in Y_t\big)\cr
\forall x\in R^n\ \forall y\in R^p\ \forall z\in R^p
\Big(\big((x,y)\in\Gamma_t \et (x,z)\in\Gamma_t\big) \Rightarrow y=z
\Big)\cr}
$$
hold true for every $t\in R^m$. By Proposition \ref{fibandfam}, the
fiber formulas for $\alpha$, which are
$$(*)_\alpha\qquad\cases{
\forall x\in \kappa(\alpha)^n\ \big((\exists y\in \kappa(\alpha)^p\ 
(x,y)\in \Gamma_\alpha)\Leftrightarrow x\in X_\alpha\big)\cr
\forall x\in \kappa(\alpha)^n\ \forall y\in \kappa(\alpha)^p\ 
\big((x,y)\in\Gamma_\alpha \Rightarrow y\in Y_\alpha\big)\cr
\forall x\in \kappa(\alpha)^n\ \forall y\in \kappa(\alpha)^p\ 
\forall z\in \kappa(\alpha)^p \Big(\big((x,y)\in\Gamma_\alpha \et
(x,z)\in\Gamma_\alpha\big)\cr
\hspace*{\fill} \Rightarrow y=z \Big)\cr}
$$
 also hold true. These fiber
formulas express the fact that $\Gamma_\alpha$ is the graph of a map
$f_\alpha: X_\alpha\to Y_\alpha$. Since $\Gamma_\alpha$ is definable,
the map $f_\alpha$ is definable for the o-minimal structure on
$\kappa(\alpha)$. \par
The next proposition says that all definable maps are obtained in
this way. We introduce a notation. If $X$ is a definable subset of
$R^m\times R^n$ and $A$ a definable subset of $R^m$, we denote by
$X_A$ the definable subset $X\cap(A\times R^n)$ of $R^m\times R^n$.
In other words, $X_A$ is the definable family $X$ restricted to $A$.

\begin{proposition}\label{fibmap} Let $\varphi: X_\alpha\to Y_\alpha$
be a definable map for the o-minimal structure on $\kappa(\alpha)$.
Then there exist $A\in \alpha$ and a definable family of maps
$f:X_A\to Y_A$ such that $\varphi=f_\alpha$.
\end{proposition}

\begin{proof} Let $\Gamma_\alpha$ be the graph of $\varphi$. Then the
fiber formulas  $(*)_\alpha$ above hold
true. By Proposition \ref{fibandfam}, there exists $A\in \alpha$
such that the formulas $(*)_t$ hold true for every $t\in A$. This means
that, for every $t\in A$, $\Gamma_t$ is the graph of a map $f_t:
X_t\to Y_t$. We obtain in this way a definable family of maps
$f:X_A\to Y_A$, which satisfies $f_\alpha=\varphi$.
\end{proof}

\debexo\label{equalfibmap} Show that $f_\alpha=g_\alpha$ if and only if
there is
$A\in
\alpha$ such that $f_t=g_t$ for every $t\in A$.
\finexo

\debexo\label{composfib} Let $f:X\to Y$ and $g:Y\to Z$ be definable
families of maps parametrized by $R^m$, $X'$ a definable subset of
$X$. For
$\alpha\in
\widetilde{R^m}$, show that
$(g\circ f)_\alpha= g_\alpha \circ f_\alpha$. Show that
$f(X')_\alpha=f_\alpha(X'_\alpha)$.
\finexo

\debexo\label{contfib} Show that $f_\alpha$ is continuous if and only
if there is $A\in \alpha$ such that $f_t$ is continuous for every
$t\in A$. (Use Proposition \ref{fibandfam}.)
\finexo

\section{Fiberwise and Global Properties}

We have seen several examples which show that a property holds for
the fiber at $\alpha$ if and only if there is $A\in \alpha$ such that
the same property holds for the fibers at $t$ for all $t\in A$. In
this section we are interested in topological properties which hold
globally for the family. An example of a  global property is the
property ``$X$ is closed'', which is stronger than the fiberwise
property ``every fiber $X_t$ is closed''. \par
The main tool will be the following. 

\begin{lemma}\label{pwclos}  Let $X$ be a definable subset of
$R^m\times R^n$. There is a partition of $R^m$ into definable sets
$C_1,\ldots,C_k$ such that, for $i=1,\ldots,k$, and for every $t\in
C_i$, we have
$\clos(X_t)=(\clos X_{C_i})_t$ (recall that $X_{C_i}$ denotes $X\cap
(C_i\times R^n)$).
\end{lemma}

\begin{proof} We proceed by induction on $m$. So we assume that the
lemma is proved for all $d$, $0\leq d<m$.  Let
$G$ be the definable subset of those
$t\in R^m$ such that $\clos(X_t)=(\clos X)_t$. 
Lemma \ref{closfibandfibclos} implies that the complement of $G$ in
$R^m$ has dimension $<m$. Choose a cdcd of $R^m$ adapted to $G$.
Every open cell $C$ of this cdcd is contained in $G$, and for every
$t\in C$ we have $\clos(X_t)=(\clos X)_t=(\clos X_C)_t$. Now let $D$
be a cell of dimension $d<m$. There is a definable homeomorphism
$\theta_D: D\to R^d$.  We apply the inductive assumption to the
definable family of $(u,x)\in R^d\times R^n$ such that
$(\theta_D^{-1}(u),x)\in X$. It follows that we can partition $D$ into
finitely many definable subsets $D_j$ such that $\clos(X_t)=(\clos
X_{D_j})_t$ for every $t\in D_j$. The lemma is proved.
\end{proof}

\begin{proposition}\label{closedfibglob} Let $X\subset Y$ be definable
subsets of
$R^m\times R^n$. Let $\alpha\in \widetilde{R^m}$. The fiber
$X_\alpha$ is closed (resp.\ open) in $Y_\alpha$ if and only if there
exists
$A\in
\alpha$ such that $X_A$ is closed (resp.\ open) in $Y_A$. 
\end{proposition}

\begin{proof} It suffices to prove the closed version. The open
version follows by taking the relative complement of $X$ in $Y$. \par
 If
$X_A$ is closed in
$Y_A$, then $X_t$ is closed in $Y_t$ for every $t\in A$. It follows
(cf. Exercise \ref{closedfib}) that $X_\alpha$ is closed in
$Y_\alpha$.\par
Conversely, assume that $X_\alpha$ is closed in $Y_\alpha$. We use
Lemma \ref{pwclos} to obtain a partition of
$R^m$ into finitely many definable sets $C_1,\ldots,C_k$ such that
$\clos(X_t)=(\clos X_{C_i})_t$ for every $t\in C_i$. Since $\alpha$
is an ultrafilter of definable sets, it contains exactly one of the
$C_i$, say $C_1$. From the assumption that $X_\alpha$ is
closed in $Y_\alpha$, it follows that there is $A\in \alpha$ such
that $X_t$ is closed in $Y_t$ for every $t\in A$ (cf. Exercise
\ref{closedfib}). Replacing $A$ with $A\cap C_1$, we can assume
$A\subset C_1$. Then $\clos(X_t)=(\clos X_A)_t$ for every $t\in
A$. We deduce $(\clos X_A)_t \cap Y_t = \clos(X_t)\cap Y_t=X_t$ for
every $t\in A$ (for the second equality we used the fact that $X_t$
is closed in $Y_t$). This implies
$(\clos X_A)\cap Y_A=X_A$, i.e.\
$X_A$ is closed in $Y_A$.
\end{proof}

\begin{theorem}\label{contfibglob} Let $f: X\to Y$ be a definable
family of maps parametrized by $R^m$, and $\alpha\in \widetilde{R^m}$.
Then
$f_\alpha$ is continuous if and only if there exists $A\in \alpha$
such that $f_A:X_A\to Y_A$ is continuous. 
\end{theorem}

\begin{proof} Say $Y$ is a definable subset of $R^m\times R^p$. We can
assume that, for all
$t$,
$Y_t$ is contained in the bounded box $(-1,1)^p$ (then $Y_\alpha$ is
contained in the bounded box $(-1,1)^p_{\kappa(\alpha)}$ in
$\kappa(\alpha)^p$). The reason for this is the following. Let $\mu:
R^p\to (-1,1)^p$ be the semialgebraic homeomorphism defined by
$$\mu(y_1,\ldots,y_p)=
\left({y_1\over \sqrt{1+y_1^2}},\ldots,{y_p\over
\sqrt{1+y_p^2}}\right)\;.$$
We can replace $f$ with the composition 
$$g=({\rm Id}_{R^m}\times
\mu)\circ f : X\to R^m\times (-1,1)^p\;.$$
Note that $g_\alpha$ is $\mu_{\kappa(\alpha)}\circ f_\alpha$, where
$\mu_{\kappa(\alpha)} : \kappa(\alpha)^p\to (-1,1)^p_{\kappa(\alpha)}$
is the semialgebraic homeomorphism defined by the same formula as
$\mu$. The continuity of $g_A$ is equivalent to the
continuity of $f_A$, and the continuity of $g_\alpha$ is equivalent
to the continuity of $f_\alpha$.\par
 Now let 
$$\Gamma=\{(t,x,y)\in X\times (-1,1)^p\;;\ f(t,x)=(t,y)\}\;.$$
The map $f_A: X_A\to A\times (-1,1)^p$ is continuous if and only if
$\Gamma_A$ is closed in $X_A\times R^p$. The map
$f_\alpha:X_\alpha\to (-1,1)^p_{\kappa(\alpha)}$ is continuous if and
only if $\Gamma_\alpha$ is closed in $X_\alpha\times
\kappa(\alpha)^p$. Hence, the theorem follows from Proposition
\ref{closedfibglob}.
\end{proof}

\debexo\label{fibhomeo} Assume that $\varphi: X_\alpha\to Y_\alpha$ is
a definable homeomorphism. Show that there is $A\in \alpha$ and a
definable family of maps $f:X_A\to Y_A$,  such that $f$ is a
homeomorphism and
$f_\alpha=\varphi$.
\finexo

\debexo Let $f$ be a definable family of maps parametrized by $R^m$.
Assume that $f_t$ is continuous for every $t\in R^m$. Show that there
is a finite partition of $R^m$ into definable subsets
$C_1,\ldots,C_k$ such that $f_{C_i}$ is continuous for $i=1,\ldots,k$.
Hint: one can use Lemma \ref{pwclos} and the proof of Theorem
\ref{contfibglob}. Another possibility is as follows:\\
1) Show that for every $\alpha\in
\widetilde{R^m}$, there is $C(\alpha)\in \alpha$ such that
$f_{C(\alpha)}$ is continuous.\\
2) Use the compactness of $\widetilde{R^m}$ to show that $R^m$ is a
finite union $R^m=C(\alpha_1)\cup\ldots\cup C(\alpha_k)$. Modify
$C(\alpha_1)$,\ldots, $C(\alpha_k)$ in order to get a partition of
$R^m$.
\finexo

\section{Triviality Theorems}

Let $X\subset R^m \times R^n$ be a definable family. Let $A$ be a
definable subset of $R^m$. We say that the family $X$ is {\em
definably trivial over $A$}\index{trivial family} if there exist a
definable set
$F$ and a definable homeomorphism $h: A\times F\to X_A$ such that the
folllowing diagram commutes:
\setlength{\unitlength}{1mm}
$$\begin{picture}(40,27)
\put(0,23){\makebox(0,0)[t]{$A \times F$}}
\put(10,22){\vector(1,0){20}}
\put(20,25){\makebox(0,0){$h$}}
\put(38,23){\makebox(0,0)[t]{$X_A$}}
\put(5,19){\vector(1,-1){13}}
\put(0,12){\makebox(0,0){projection}}
\put(35,19){\vector(-1,-1){13}}
\put(40,12){\makebox(0,0){projection}}
\put(25,2){\makebox(0,0){$A\subset{\Bbb R}^m$}}
\end{picture}$$
We say that $h$ is a {\em definable trivialization of $X$ over $A$}.
Now let $Y$ be a definable subset of $X$. We say that the
trivialization $h$ is {\em compatible with $Y$} if there is a definable
subset $G$ of $F$ such that $h(A\times G)=Y_A$. Note that if $h$ is
compatible with $Y$, its restriction to $Y_A$ is a trivialization of
$Y$ over $A$.

\begin{theorem}[Hardt's Theorem for Definable Families]\label{Hardt}\
\\
 Let $X\subset R^m \times R^n$ be a definable family. Let
$Y_1,\ldots,Y_\ell$ be definable subsets of $X$. There exists a finite
partition of $R^m$ into definable sets
$C_1$, \ldots, $C_k$ such that $X$ is definably trivial over each
$C_i$ and, moreover, the trivializations over each $C_i$ are
compatible with $Y_1,\ldots,Y_\ell$.
\end{theorem}

This theorem was proved by R.~Hardt in the semialgebraic case.

\begin{proof} 
We can assume that
$X$ is closed and contained in $R^m\times[-1,1]^n$. The reason for
this is the following.  Let $\mu: R^n\to (-1,1)^n$ be the semialgebraic
homeomorphism defined by
$$\mu(x_1,\ldots,x_n)=
\left({x_1\over \sqrt{1+x_1^2}},\ldots,{x_n\over
\sqrt{1+x_n^2}}\right)\;.$$
First, we can replace $X$ (and each $Y_j$) with its image by ${\rm
Id}_{R^m}\times
\mu$, which is contained in $R^m\times(-1,1)^n$. Then, we can replace
$X$ with $\clos X$, which is closed and contained in
$R^m\times[-1,1]^n$, and add  $X$ to the list of definable
subsets $Y_1$, \ldots, $Y_\ell$. A definable trivialization of
$\clos(X)$ over $C_i$ compatible with $X$ will induce a
definable trivialization of $X$ over $C_i$.\par
Let $\alpha$ be an element of $\widetilde{R^m}$. The assumption above
implies that $X_\alpha$ is bounded and closed in $\kappa(\alpha)^n$.
Hence, by the triangulation theorem \ref{triangsets}, there is a
finite simplicial complex $K$ with vertices in ${\Bbb Q}^n\subset R^n$
and a definable homeomorphism $\Phi: |K|_{\kappa(\alpha)}\to X_\alpha$,
such that each $(Y_j)_\alpha$ is the image by $\Phi$ of a union
$(V_j)_{\kappa(\alpha)}$ of open simplices of $K$. By Lemma
\ref{extpoly}, we have $|K|_{\kappa(\alpha)}=(R^m\times
|K|_R)_\alpha$. Hence (cf.\ Exercise \ref{fibhomeo}), there exists
$C(\alpha)\in \alpha$ and a definable family of maps $h :
C(\alpha)\times |K|_R \to X_A$ such that $h$ is a homeomorphism and
$h_\alpha=\Phi$. Moreover, since 
$$h_\alpha((R^m\times (V_j)_R)_\alpha)
= h_\alpha((V_j)_{\kappa(\alpha)})= (Y_j)_\alpha\;,$$
we can assume, replacing $C(\alpha)$ with a smaller definable set
still in $\alpha$, that $h( C(\alpha)\times
(V_j)_R)=(Y_j)_{C(\alpha)}$. We have proved the following fact: for
every $\alpha\in \widetilde{R^m}$, there exist $C(\alpha)\in \alpha$
and a definable trivialization of $X$ over
$C(\alpha)$ compatible with each $Y_j$.\par
 The open and closed subset
$\widetilde{C(\alpha)}$ cover $\widetilde{R^m}$. Since 
$\widetilde{R^m}$ is compact, we can extract a finite subcover from
this cover. Hence, there are $\alpha_1$, \ldots, $\alpha_k$ such that
$R^m= C(\alpha_1)\cup\ldots\cup C(\alpha_k)$. Replacing $C(\alpha_i)$
with $C_i= C(\alpha_i)\setminus \bigcup_{j<i} C(\alpha_j)$, we obtain
a partition of $R^m$ into definable sets $C_1,\ldots,C_k$, and, for
each $i$, a definable trivialization of $X$ over $C_i$ compatible with
each $Y_j$.
\end{proof}

\debexo Use Hardt's theorem to give another proof of the local conic
structure of definable sets (cf.\ \ref{localconic}). Hint: trivialize
the definable family of 
$$A_\varepsilon=\{x\in A\;;\  \Vert x-a \Vert=\varepsilon\}\;.$$
\finexo

We have a similar trivialization theorem for definable familes of
functions (with values in $R$). First, we have to define the notion
of trivialization for a family of functions.\par
Let 
$f:X\to R^m\times R$ be definable family of functions parametrized
by $R^m$ and $A$ a definable subset of $R^m$. We say that $f$ is {\em
definably trivial over $A$}\index{trivial family} if there exist a
definable function
$g: F\to R$ and definable homeomorphisms
$\rho : A\times F\to X_A$ and
$\lambda : A\times R\to A\times R$ such that the following diagram
commutes,
\setlength{\unitlength}{1mm}
$$\begin{picture}(40,40)
%je place le diagramme
\put(3,36){\makebox(0,0){
%je fais le diagramme
\begin{picture}(0,0)
\put(0,0){\makebox(0,0){$A\times F$}}
\put(7,0){\vector(1,0){14}} \put(13,2){\makebox(0,0){$\rho$}}
\put(26,0){\makebox(0,0){$X_A$}}
\put(0,-3){\vector(0,-1){12}} 
\put(-8,-9){\makebox(0,0){${\rm Id}_A\times g$}}
\put(26,-3){\vector(0,-1){12}}  \put(29,-9){\makebox(0,0){$f$}}
\put(0,-18){\makebox(0,0){$A\times R$}}
\put(20,-18){\vector(-1,0){14}} \put(13,-20){\makebox(0,0){$\lambda$}}
\put(26,-18){\makebox(0,0){$A\times R$}}
\put(0,-21){\vector(1,-1){11}}
\put(26,-21){\vector(-1,-1){11}}
\put(13,-35){\makebox(0,0){$A$}}
\end{picture}
%fin du diagramme 
}}
\end{picture}$$
where the maps to $A$ are the projections on the first factor. The
couple
$(\rho,\lambda)$ is called a definable trivialization of $f$ over
$A$. If $Y$ is a definable subset of $X$, we say that the
trivialization  $(\rho,\lambda)$ is compatible with $Y$ if there is
a definable subset $G$ of $F$ such that $\rho(A\times G)=Y_A$. Note
that, in this case, $(\rho,\lambda)$ induces a definable
trivialization of $f|_Y$ over $A$.

\begin{theorem}[Definable Triviality of Families of
Functions]
\label{Hardtfunctions}\ \\
Let 
$f:X\to R^m\times R$ be a definable family of functions parametrized by
$R^m$. Let $Y_1,\ldots, Y_\ell$ be definable subsets of $X$. There
exists a finite partition of $R^m$ into definable sets
$C_1$, \ldots, $C_k$ such that $f$ is definably trivial over each
$C_i$ and, moreover, the trivializations over each $C_i$ are
compatible with $Y_1,\ldots,Y_\ell$.
\end{theorem}

\begin{proof}
We can assume that $X$ is closed in $R^m\times [-1,1]\times
[-1,1]^n$ and that $f$ is the restriction to $X$ of the projection $p$
on the first
$m+1$ coordinates. The reason for this is the following.  Let $\nu:
R\to (-1,1)$ be the semialgebraic homeomorphism defined by
$\nu(y)=y/\sqrt{1+y^2}$, and $\mu : R^n\to (-1,1)^n$ the
semialgebraic homeomorphism defined by
$\mu(x_1,\ldots,x_n)=(\nu(x_1),\ldots,\nu(x_{n}))$. First we
replace $X$ with its homeomorphic image $\Gamma\subset R^m\times
(-1,1)\times (-1,1)^n$ by the definable map
$$(t,x)\longmapsto (t,\nu(f_t(x)),\mu(x))\;.$$ We
replace also the
$Y_j$ with their images $\Delta_j\subset \Gamma$. We obtain the
following commutative diagram, where the horizontal maps are
definable homeomorphisms and the maps to $R^m$ are the
projections on the $m$ first coordinates.
\setlength{\unitlength}{1mm}
$$\begin{picture}(44,44)
%je place le diagramme
\put(3,40){\makebox(0,0){
%je fais le diagramme
\begin{picture}(0,0)
\put(0,0){\makebox(0,0){$\Gamma$}}
\put(5,0){\vector(1,0){26}}
\put(36,0){\makebox(0,0){$X$}}
\put(0,-3){\vector(0,-1){12}} 
\put(-3,-9){\makebox(0,0){$p$}}
\put(36,-3){\vector(0,-1){12}}  \put(39,-9){\makebox(0,0){$f$}}
\put(-4,-18){\makebox(0,0){$R^m\times (-1,1)$}}
\put(29,-18){\vector(-1,0){20}} 
\put(18,-21){\makebox(0,0){${\rm Id}_{R^m}\times\nu$}}
\put(38,-18){\makebox(0,0){$R^m\times R$}}
\put(0,-21){\vector(1,-1){15}}
\put(36,-21){\vector(-1,-1){15}}
\put(18,-39){\makebox(0,0){$R^m$}}
\end{picture}
%fin du diagramme 
}}
\end{picture}$$
Now we replace
$\Gamma$ with
$\clos(\Gamma)\subset R^m\times [-1,1]\times [-1,1]^n$, and we add
$\Gamma$ to the list of definable subsets
$\Delta_1,\ldots,\Delta_\ell$. A definable trivialization of
$p|_{{\rm clos}(\Gamma)}$ compatible with
$\Gamma,\Delta_1,\ldots,\Delta_\ell$ will induce a
definable trivialization of 
$p|_\Gamma$ compatible with
$\Delta_1,\ldots,\Delta_\ell$. Composing with the
definable homeomorphisms represented by the horizontal arrows in the
diagram above, we obtain a definable trivialization of
$f$ compatible with $Y_1,\ldots,Y_\ell$.\par
Let $\alpha\in \widetilde{R^m}$. The assumption above implies that
$X_\alpha$ is closed and bounded in
$\kappa(\alpha)\times\kappa(\alpha)^n$ and $f_\alpha$ is the
restriction to
$X_\alpha$ of the projection $\pi$ on the first coordinate. By
Proposition \ref{triang_1coordrat}, there exist a finite simplicial
complex $K$ with vertices in
${\Bbb Q}\times {\Bbb Q}^n$ and definable homeomorphisms  $\Phi:
|K|_{\kappa(\alpha)}\to X_\alpha$ and
$\tau:\kappa(\alpha)\to \kappa(\alpha)$  such that
the following left-hand side diagram commutes and each
$Y_j$ is a union of images by $\Phi$ of open simplices. 
\setlength{\unitlength}{1mm}
$$\begin{picture}(100,40)
%je place le diagramme
\put(3,36){\makebox(0,0){
%je fais le diagramme
\begin{picture}(0,0)
\put(0,0){\makebox(0,0){$|K|_{\kappa(\alpha)}$}}
\put(8,0){\vector(1,0){13}} \put(13,2){\makebox(0,0){$\Phi$}}
\put(26,0){\makebox(0,0){$X_\alpha$}}
\put(0,-3){\vector(0,-1){12}} 
\put(-3,-9){\makebox(0,0){$\pi$}}
\put(26,-3){\vector(0,-1){12}}  \put(29,-9){\makebox(0,0){$f_\alpha$}}
\put(0,-18){\makebox(0,0){$\kappa(\alpha)$}}
\put(21,-18){\vector(-1,0){16}} 
\put(13,-21){\makebox(0,0){$\tau$}}
\put(26,-18){\makebox(0,0){$\kappa(\alpha)$}}
\end{picture}
%fin du diagramme 
}}
%je place le deuxieme diagramme
\put(63,36){\makebox(0,0){
%je fais le deuxieme diagramme
\begin{picture}(0,0)
\put(-4,0){\makebox(0,0){$C(\alpha)\times |K|_R$}}
\put(8,0){\vector(1,0){12}} \put(13,2){\makebox(0,0){$\rho$}}
\put(26,0){\makebox(0,0){$X_{C(\alpha)}$}}
\put(0,-3){\vector(0,-1){12}} 
\put(-3,-9){\makebox(0,0){$p$}}
\put(26,-3){\vector(0,-1){12}}  \put(29,-9){\makebox(0,0){$f$}}
\put(-2,-18){\makebox(0,0){$C(\alpha)\times R$}}
\put(19,-18){\vector(-1,0){12}} 
\put(13,-21){\makebox(0,0){$\lambda$}}
\put(28,-18){\makebox(0,0){$C(\alpha)\times R$}}
\put(0,-21){\vector(1,-1){11}}
\put(26,-21){\vector(-1,-1){11}}
\put(13,-35){\makebox(0,0){$C(\alpha)$}}
\end{picture}
%fin du deuxieme diagramme 
}}
\end{picture}$$ 
We obtain
$C(\alpha)\in \alpha$ and  definable families of maps $\rho :
C(\alpha)\times |K|_R\to X_{C(\alpha)}$ and $\lambda :
C(\alpha)\times R\to C(\alpha)
\times R$ which are homeomorphisms and such that $\rho_\alpha=\Phi$
and $\lambda_\alpha=\tau$. After possibly shrinking $C(\alpha)$,
$(\rho,\lambda)$ is a definable trivialization of
$f$ over $C(\alpha)$, compatible with the $Y_j$. We conclude as in the
proof of Hardt's theorem.
\end{proof}

\section{Topological Types of Sets and Functions}

Two embeddings $Y\subset X$ and $Y'\subset X'$, where all sets are
definable, are said to have the same {\em (definable) topological type}
if there is a (definable) homeomorphism $h: X\to X'$ such that
$h(Y)=Y'$. Two definable maps $f:X\to Y$ and $f': X'\to Y'$ are said
to have the same (definable) topological type if there are (definable)
homeomorphisms $\rho : X\to X'$ and $\lambda : Y'\to Y$ such that
$\lambda\circ f'\circ \rho=f$. \par
Assume that the definable family $X\subset R^m \times R^n$ is
definably trivial over $A\subset R^m$, and that the trivialization is
compatible with $Y\subset X$. This means that there is a definable
homeomorphism $h: A\times F\to X_A$ of the form $h(t,x)=
(t,h_t(x))$, and a definable subset $G$ of $F$ such that $h(A\times
G)=Y_A$. Now let $t$ and $t'$ be two points of $A$. Then
$h_{t'}\circ h_t^{-1}$ is a definable homeomorphism from $X_t$ onto
$X'_t$, and the image of $Y_t$ by this homeomorphism is $Y_{t'}$.
Hence, all embeddings $Y_t\subset X_t$ for $t\in A$ have the same
topological type. From this discussion and Hardt's theorem we obtain
the following.

\begin{corollary}\label{fintoptypes} Let $X\subset R^m \times R^n$ be a
definable  family, $Y$ a definable subset of $X$. There are finitely
many definable topological types of embeddings $Y_t\subset X_t$ for
$t\in R^m$.
\end{corollary}

We give an example of application.

\begin{theorem} Given $n$ and $d$ positive integers, there are
finitely many topological types of embedding $V\subset R^n$, where
$V$ is an algebraic subset defined by equations of degrees at most
$d$.
\end{theorem}

\begin{proof} If $V$ is defined by equations $P_1=\ldots=P_\ell=0$ of
degrees at most $d$, it is also defined by the single equation
$P_1^2+\cdots+P_\ell^2=0$ of degree at most $2d$. Hence, it suffices
to show that there are finitely many topological types of embedding
$V\subset R^n$, where $V$ is defined by one equation of degree at
most $2d$. If $i=(i_1,\ldots,i_n)\in {\Bbb N}^n$, we set
$|i|=i_1+\cdots+i_n$ and $x^i= x_1^{i_1}\cdots x_n^{i_n}$. Let
$I_{2d}$ be the finite set of $i\in {\Bbb N}^n$ such that $|i|\leq
2d$. We denote by $a=(a_i)_{i\in I_{2d}}$ an element of the affine
space $R^{I_{2d}}$. Set
$$Y= \{ (a,x)\in R^{I_{2d}}\times R^n\;;\ \sum_{i\in I_{2d}}
a_ix^i=0\}\subset R^{I_{2d}}\times R^n\;.$$
For $a \in R^{I_{2d}}$, $Y_a$ is an algebraic subset of $R^n$ defined
by one equation of degree $\leq 2d$, and we obtain in this way all
such algebraic subsets. Hence, the theorem follows from Corollary
\ref{fintoptypes} applied to the semialgebraic families $Y\subset 
R^{I_{2d}}\times R^n$ parametrized by $R^{I_{2d}}$.
\end{proof}

We now turn to the case of families of functions. Assume that the
definable family of functions $f_t: X_t\to R$ parametrized by $t\in
R^m$ is trivial over $A$. This means that there are definable
homeomorphism $\rho : A\times F\to X_A$ of the form
$\rho(t,x)=(t,\rho_t(x))$ and $\lambda: A\times R\to A\times R$ of
the form $\lambda(t,y)=(t,\lambda_t(y))$ and a definable function
$g: F\to R$, such that $\lambda\circ f\circ \rho= {\rm Id}_A\times
g$. Let $t$ and $t'$ be two points of $A$. Then  
$\rho_{t'}\circ \rho_t^{-1} : X_t\to
X_{t'}$ and
$\lambda_t^{-1}\circ \lambda_{t'}: R\to R$ are definable homeomorphisms
which satisfy
$$\lambda_t^{-1}\circ \lambda_{t'}\circ f_{t'}\circ
\rho_{t'}\circ \rho_t^{-1} = f_t\;.$$
This shows that all $f_t$ for $t\in A$ have the same definable
topological type. From this discussion and Theorem
\ref{Hardtfunctions} we obtain the following.

\begin{corollary}\label{fintoptypesfunctions} In a 
definable  family of functions $f_t: X_t\to R$ parametrized by
$t\in R^m$, there are finitely many definable topological types.
\end{corollary}

The similar statement for definable families of maps is not true. For
instance, there is a semialgebraic family of maps ${\Bbb R}^2\to
{\Bbb R}^2$ in which there are infinitely many topological types
(cf.\ \cite{CR}).\par
\medskip

We give an example of application for the o-minimal structure ${\Bbb
R}_{\exp}$ expanding ${\Bbb R}$ and which is generated by the
exponential function. This means that ${\Bbb R}_{\exp}$ is
a structure $(S_n)_{n\in {\Bbb N}}$, where
$S_n\subset {\cal P}({\Bbb R}^n)$ and $S_2$ contains the graph of the
exponential function, and $S_n\subset T_n$ for every structure
$(T_n)_{n\in {\Bbb N}}$ with these properties. The fact that ${\Bbb
R}_{\exp}$ is o-minimal (i.e.\ $S_1$ consists of finite unions of
points and intervals) is a fundamental result of A.~Wilkie \cite{Wi}.
In
${\Bbb R}_{\exp}$ the function $(x,\lambda)\mapsto x^\lambda$ on $\{x
>0\}\times {\Bbb R}$ is definable, since
$$y=x^\lambda \Leftrightarrow \exists z\in {\Bbb R} (x=\exp(z)
\et y=\exp(\lambda z))\;.$$

\begin{theorem}\label{fewnomials} Given $n$ and $k$ two positive
integers, there are finitely many topological types in the family of
polynomials
${\Bbb R}^n\to {\Bbb R}$ with at most $k$ monomials (and no bound on
the degree).
\end{theorem}

\begin{proof} 
We extend the power function
$(x,\lambda)\mapsto x^\lambda$  to two
definable functions
$M_\epsilon : {\Bbb R}^2\to {\Bbb R}$ for $\epsilon=0,1$, defined by
$$M_\epsilon(x,\lambda)=\cases{ x^\lambda & if $x>0\;,$\cr 0 & if
$x=0$ and $\lambda\not=0\;,$\cr 1 & if $x=0$ and $\lambda=0\;,$\cr
(-1)^\epsilon \,|x|^\lambda & if $x<0\;.$\cr}$$

Now consider the family of all functions ${\Bbb R}^n\to {\Bbb R}$ 
$$(x_1,\ldots,x_n)\mapsto \sum_{i=1}^k  \left(a_i\prod_{j=1}^n
M_{\epsilon_{i,j}}(x_j,\lambda_{i,j})\right)\;.$$ This is a definable
family of functions parametrized by the
$((a_i),(\lambda_{i,j}),(\epsilon_{i,j}))\in {\Bbb R}^k\times
{\Bbb R}^{kn}\times
\{0,1\}^{kn}$, which is a definable subset of ${\Bbb R}^{k+2kn}$. In
this family we have all polynomials in
$n$ variables with at most
$k$ monomials. Hence, the theorem follows from Corollary
\ref{fintoptypesfunctions}.\end{proof}

Note that the finiteness of topological
types for polynomials ${\Bbb R}^n\to {\Bbb R}$ with bounded degree (a
result of Fukuda \cite{Fu} obtained by completely different methods)
can be obtained by working in the semialgebraic structure (it is not
difficult to put all polynomials ${\Bbb R}^n\to {\Bbb R}$ with degree
$\leq d$ in a semialgebraic family). Theorem \ref{fewnomials} is
stronger. Its statement does not refer to o-minimal structures, but
its proof uses this theory in an essential way.

%% file: OMINsix.tex
\chapter{Smoothness}

{\em In this chapter, we assume for simplicity that $R={\Bbb R}$,
i.e.\ we consider an o-minimal expansion of the field ${\Bbb R}$ of
real numbers}. The results of this chapter still hold for o-minimal
structures expanding an arbitrary real closed field. In order to
prove these results, one would first have to establish the basic
facts of analysis (such as the local inversion theorem) for
definable mappings. This is done in \cite{vdD}.

\section[Definable Functions in One Variable]{Derivability of Definable
Functions in One Variable}

If $A$ is a definable subset of ${\Bbb R}^n$, we say that a function
$g: A\to {\Bbb R}\cup\{-\infty,+\infty\}$ is {\em definable} if the
inverse images of $-\infty$ and of $+\infty$ are definable subsets of
$A$ and the restriction of $g$ to $g^{-1}({\Bbb R})$ is a definable
function with values in ${\Bbb R}$. 

\begin{lemma}\label{lrderivative} Let $f: I\to {\Bbb R}$ be a definable
continuous function on an open interval $I$ of ${\Bbb R}$. Then $f$ has
left and right derivatives in ${\Bbb R}\cup\{-\infty,+\infty\}$ at
every point
$x$ in $I$. These derivatives are denoted respectively by
$f'_\ell(x)$ and $f'_r(x)$. The functions $f'_\ell$ and $f'_r$ from
$I$ to ${\Bbb R}\cup\{-\infty,+\infty\}$ are definable.
\end{lemma}

\begin{proof} Apply the monotonicity theorem to the function
$y\mapsto (f(y)-f(x))/(y-x)$ to show that it has a limit in
${\Bbb R}\cup \{-\infty,+\infty\}$ as $y$ tends to
$x$ from the left or from the right (cf.\ Exercise
\ref{leftrightlimit}). The definability of the left and right
derivatives is left as an exercise.
\end{proof}

\begin{lemma}\label{variation}  Let $f: I\to {\Bbb R}$ be a definable
continuous function on an open interval $I$ of ${\Bbb R}$. If $f'_\ell
>0$ (resp.\ $f'_r>0$) on $I$, $f$ is strictly increasing on $I$.
\end{lemma}

\begin{proof} Apply the monotonicity theorem to the function $f$,
plus the fact that there is no subinterval of $I$ on which $f$ is
constant or $f$ is strictly decreasing (otherwise, one would have
$f'_\ell\leq 0$ and $f'_r\leq0$ on such a subinterval).
\end{proof}

Remark that the preceding lemma can be proved without the assumption
that $f$ is definable. 

\debexo Assume that $f$ is any function having a
positive (or $+\infty$) left derivative at every point of $I$. Take
$a\in I$. Show that 
$$\sup\{x\in I\;;\ f\geq f(a)\ {\rm on}\ [a,x]\}$$
is the right endpoint of $I$ (or $+\infty$ if $I$ is not bounded
from above). Deduce that $f$ is non decreasing on $I$ and, hence,
that $f$ is strictly increasing.
\finexo

\begin{theorem}\label{pwderiv1variable}
Let $f: I\to {\Bbb R}$ be a definable continuous function on an open
interval $I$ of ${\Bbb R}$. Then, for all but finitely many points in
$I$, we have $f'_\ell(x)=f'_r(x)\in {\Bbb R}$. Hence, $f$ is
derivable outside a finite subset of $I$.
\end{theorem}

\begin{proof}
First we prove that there is no subinterval of $I$ on which
$f'_\ell = +\infty$, or $f'_\ell = -\infty$, or $f'_r = +\infty$, or
$f'_r = -\infty$. Suppose for instance that $f'_\ell = +\infty$ on a
subinterval $J$ of $I$. Take $a<b$, both in $J$ and set
$$g(x) = f(x)- {f(b)-f(a)\over b-a}\, x\qquad\mbox{\rm for}\
x\in J\;.$$
We have $g'_\ell=+\infty$ on $J$. By Lemma \ref{variation}, $g$ is
strictly increasing on $J$. We deduce $g(b)>g(a)$,
which is clearly impossible. We can show in a similar way that the
other cases lead to a contradiction.\par
Suppose now that there is a subinterval $J$ of $I$ on which
$f'_\ell$ and $f'_r$ have values in ${\Bbb R}$ and $f'_\ell<f'_r$.
Replacing $J$ with a smaller subinterval, we can assume that
$f'_\ell$ and $f'_r$ are continuous. Taking a still smaller
subinterval, we can assume that there is $c\in {\Bbb R}$ such that
$f'_\ell < c < f'_r$. Lemma \ref{variation} implies that $x\mapsto
f(x)-cx$ is at the same time strictly increasing and strictly
decreasing on this subinterval, which is impossible. We can show in a
similar way that there is no subinterval on which $f'_\ell$ and $f'_r$
have values in ${\Bbb R}$ and $f'_\ell>f'_r$.\par
The definability of $f'_\ell$ and $f'_r$, together with the facts
just proved, imply the conclusion of the theorem.
\end{proof}

\begin{corollary}\label{pwCk1variable}
Let $f:I\to {\Bbb R}$ be a definable function. For all $k\in {\Bbb
N}$, there exists a finite subset $M(k)$ of $I$ such that $f$ is of
class ${\cal C}^k$ on $I\setminus M(k)$.
\end{corollary}

\begin{proof} By induction on $k$, using Theorem
\ref{pwderiv1variable}, the definability of the derivative, and
piecewise continuity.
\end{proof}

\section{${\cal C}^k$ Cell Decomposition}

A ${\cal C}^k$ cylindrical definable cell decomposition of ${\Bbb
R}^n$ is a cdcd satisfying extra smoothness conditions which imply, in
particular, that each cell is a ${\cal C}^k$ submanifold of ${\Bbb
R}^n$.\par
\begin{itemize}
\item A ${\cal C}^k$ cdcd of ${\Bbb R}$ is any cdcd of ${\Bbb R}$
(i.e. a finite subdivision of ${\Bbb R}$).
\item If $n>1$, a ${\cal C}^k$ cdcd of ${\Bbb R}^n$ is given by a
${\cal C}^k$ cdcd of ${\Bbb R}^{n-1}$ and, for each cell $D$ of
${\Bbb R}^{n-1}$, definable functions of class ${\cal C}^k$
$$\zeta_{D,1}<\ldots<\zeta_{D,\ell(D)} : D\to {\Bbb R}\;.$$
The cells of ${\Bbb R}^n$ are, of course, the graphs of the
$\zeta_{D,i}$ and the bands delimited by these graphs.
\end{itemize}
It is clear from the description by induction that the cells of a
${\cal C}^k$ cdcd are ${\cal C}^k$ submanifolds of ${\Bbb R}^n$.
Moreover, for each cell $C$, there is a definable ${\cal
C}^k$ diffeomorphism $\theta_C: C\to {\Bbb R}^{\dim C}$. The
diffeomorphism $\theta_C$ can be defined by induction on $n$, by
the same formulas as in Proposition \ref{cell=affine}.\par

We now state the two theorems of this section.

\begin{theorem}[${\cal C}^k$ Cell Decomposition: ${\cal C}^k{\rm
CDCD}_n$]\label{CkCDCD}
Given finitely many\break definable subsets $X_1,\ldots,X_\ell$ of
${\Bbb R}^n$, there is a ${\cal C}^k$ cdcd of ${\Bbb R}^n$ adapted to
$X_1,\ldots,X_\ell$ (i.e.\ each $X_i$ is a union of cells).
\end{theorem}

\begin{theorem}[Piecewise ${\cal C}^k$: ${\rm P}{\cal
C}^k_n$]\label{PCk}
Given a definable function $f: A\to {\Bbb R}$, where $A$ is a
definable subset of ${\Bbb R}^n$, there
is a finite partition of $A$ into definable ${\cal C}^k$
submanifolds  $C_1,\ldots,C_\ell$, such that each restriction
$f|_{C_i}$ is ${\cal C}^k$.
\end{theorem}

\begin{proof} We prove the two theorems simultaneously, by induction
on $n$. The case $n=1$ is already done: ${\cal C}^k{\rm CDCD}_1$ is
obvious, and ${\rm P}{\cal C}^k_1$ follows from Theorem
\ref{pwderiv1variable}. So we assume $n>1$ and the theorems proved
for smaller dimensions.\par
1) ${\cal C}^k{\rm CDCD}_n$. Start with a cdcd of ${\Bbb R}^n$
adapted to $X_1,\ldots,X_\ell$. Using ${\rm P}{\cal
C}^k_{n-1}$, partition the cells of the induced cdcd of ${\Bbb
R}^{n-1}$ in order to have all $\zeta_{C,j}$ of class ${\cal C}^k$.
Then, using ${\cal C}^k{\rm CDCD}_{n-1}$, refine to a ${\cal C}^k$
cdcd of ${\Bbb R}^{n-1}$.\par
2) ${\rm P}{\cal C}^k_n$. We use the following lemma.
\begin{lemma}\label{almosteverywhereCk} Let $g:U\to {\Bbb R}$ be a
definable function, with
$U$ definable open subset of ${\Bbb R}^n$. Then there is a definable
open subset $V$ of $U$ such that $f|_V$ is ${\cal C}^k$ and
$\dim(U\setminus V)<n$.
\end{lemma}
\begin{proof} By induction on $k$, it is sufficient to show that the
set
$G_i$ where the partial derivative $\partial g/\partial x_i$
exists is definable, that its complement in $U$ has empty interior,
and that $\partial g/\partial x_i$ is definable on $G_i$. We leave the
checking of definability to the reader and prove that $U\setminus
G_i$ has empty interior. Otherwise, there would exist a nonempty
open box where $\partial g/\partial x_i$ does not exists. Considering
the restriction of $g$ to an interval of a line parallel to the $x_i$
axis contained in this box, we obtain a contradiction with Theorem
\ref{pwderiv1variable}.
\end{proof}
 We return to
the proof of
${\rm P}{\cal C}^k_n$. We choose a cdcd of
${\Bbb R}^n$ adapted to $A$. By Lemma \ref{almosteverywhereCk}, for
each open cell $C_i$ of ${\Bbb R}^n$ contained in $A$, there is a
definable open subset $V_i$ such that $f|_{V_i}$ is ${\cal C}^k$ and
$\dim(C_i\setminus V_i)<n$. By ${\cal C}^k{\rm CDCD}_n$, we can refine
the cdcd to a ${\cal C}^k$ cdcd of ${\Bbb R}^n$ adapted to $A$ and the
$V_i$. On each open cell of this new cdcd contained in $A$, $f$ is
${\cal C}^k$. Let $D$ be a cell of dimension $<n$ contained in $A$.
Using a ${\cal C}^k$ definable diffeomorphism from $D$ to ${\Bbb
R}^{\dim D}$ and the inductive assumption, we can partition $D$ into
finitely many definable ${\cal C}^k$ submanifolds on which $f$ is
${\cal C}^k$. This completes the proof of the two theorems.
\end{proof}

\section{Definable Manifolds and Tubular Neighborhoods}

Let $M\subset {\Bbb R}^n$ be a definable ${\cal C}^k$ submanifold (we
always assume $1\leq k< \infty$). We introduce the tangent and normal
bundles of
$M$.\par The {\em tangent
bundle}\index{tangent bundle} ${\rm T}M$ is the set of $(x,v) \in M
\times {\Bbb R}^n$ such that $v$ is a tangent vector to $M$ at $x$. We
denote by $p: TM\to M$  the projection defined by $p(x,v)=x$ and by
${\rm T}_x M$ the tangent space to $M$ at $x$, i.e.\ ${\rm T}_x M =
p^{-1}(x)$ for
$x\in M$. We can argue as follows in order to prove that
$TM$ is a definable subset of ${\Bbb R}^n\times{\Bbb R}^n$. We
consider the set
$S$ of triples $(x,y,v)$, where $v$ is a vector parallel to the line
joining two distinct points $x$ and $y$ of $M$, that is
$$S=\{ (x,y,\lambda(y-x))\in M\times M\times{\Bbb R}^n\;;\ x\not=y
\ {\rm and}\ \lambda\in{\Bbb R}\}\;.$$
This set $S$ is obviously definable. Its closure ${\rm clos}(S)$ in
$M\times M\times{\Bbb R}^n$ is also definable. We claim that
$${\rm T}M=\{(x,v)\in M\times {\Bbb R}^n\;;\ (x,x,v)\in {\rm
clos}(S)\}\;,$$ which shows that it is definable. 

\debexo Prove the claim and show in the same time that ${\rm T}M$ is
a ${\cal C}^{k-1}$ submanifold (this is a classical fact of
differential geometry). Hint: One can assume that the origin $0$
belongs to $M$ and work in a neighborhood of $0$. One can also
assume that the first $d=\dim(M)$ cordinates $x_1,\ldots,x_d$ are the
coordinates of a chart of $M$ in a neighborhood of $0$. Then the other
coordinates are ${\cal C}^k$ functions $x_j=\xi_j(x_1,\ldots,x_d)$ on
a neighborhood of $0$ in $M$ (for $j=d+1,\ldots,n$). Then show that
$(0,0,v)$ belongs to the closure of $S$ if and only if $v$ is a
linear combination of the vectors
$$ e_i + \sum_{j=d+1}^n {\partial \xi_j\over x_i}(0) e_j \qquad
\mbox{for } i=1,\ldots,d\;,$$ where $e_i$ denotes the $i$-th vector of
the canonical basis of
${\Bbb R}^n$. These vectors form a basis of ${\rm T}_0M$.
\finexo

\debexo Let $f : M\to {\Bbb R}$ be a definable ${\cal C}^k$ function.
Prove that ${\rm d}f : {\rm T}M\to {\rm T}{\Bbb R}={\Bbb
R}\times{\Bbb R}$ is a definable map (of class ${\cal C}^{k-1}$).
Hint: construct the graph of $df$ replacing the set $S$ above with
the set of 
$$(x,y,\lambda(y-x),f(x),\lambda(f(y)-f(x)))\;,$$
 where $x$ and $y$ are
distinct point of $M$ and $\lambda\in {\Bbb R}$.
\finexo

The normal bundle\index{normal bundle} ${\rm N}M$ is the set of $(x,v)$
in
$M\times {\Bbb R}^n$ such that $v$ is orthogonal to ${\rm T}_x M$.
This is a
${\cal C}^{k-1}$ submanifold of ${\Bbb R}^n\times{\Bbb R}^n$, and it is
definable since ${\rm T}M$ is definable.\par
Now we introduce a definable ${\cal C}^{k-1}$ map
\begin{eqnarray*}
\varphi : {\rm N} M &\longrightarrow & {\Bbb R}^n\\
 (x,v)& \longmapsto & (x+v)
\end{eqnarray*}
The map $\varphi$ induces the canonical diffeomorphism from the zero
section
$M\times\{0\}$ of the normal bundle onto $M$, and it is a local
diffeomomorphism at each point $(x,0)$. Indeed,  the derivative
${\rm d}_{(x,0)} \varphi : {\rm T}_xM\times {\rm N}_x M \to {\Bbb
R}^n$ is the isomorphism which maps $(\xi,v)$ to $\xi+v$. \par
Our aim in this section is to prove the following result, which is
the definable version of the tubular neighborhood theorem.

\begin{theorem}[Definable Tubular Neighborhood]\label{tubnbhd}
Let $M$ be a definable ${\cal C}^{k}$ submanifold of ${\Bbb R}^n$.
There exists a definable open neighborhood $U$ of the zero-section
$M\times\{0\}$ in the normal bundle ${\rm N}M$ such that the
restriction $\varphi|_U$ is a ${\cal C}^{k-1}$ diffeomorphism onto an
open neighborhood $\Omega$ of $M$ in ${\Bbb R}^n$. Moreover, we can
take $U$ of the form 
$$U=\{ (x,v)\in {\rm N}M\;;\ \Vert v\Vert < \varepsilon(x)\}\;,$$
where $\varepsilon$ is a positive definable ${\cal C}^k$ function on
$M$.
\end{theorem}

A neighborhood $\Omega$ as in the theorem above is called a {\em
definable ${\cal C}^{k-1}$ tubular neighborhood of $M$}\index{tubular
neighborhood}. We have on
$\Omega$ a definable ${\cal C}^{k-1}$ retraction $\pi : \Omega\to M$
and a definable ${\cal C}^{k-1}$ ``square of distance function'' $\rho
:
\Omega\to {\Bbb R}$, which are defined by
\begin{eqnarray*}
\pi(\varphi(x,v))&=&x\\
\rho(\varphi(x,v))&=&\Vert v\Vert^2\;.
\end{eqnarray*}
\par

For the proof of Theorem \ref{tubnbhd}, we need to consider first the
case of closed definable submanifolds.

\begin{lemma}\label{Ckminorclosed}
Let $M$ be a definable ${\cal C}^k$ submanifold of ${\Bbb R}^n$,
closed in ${\Bbb R}^n$. Let $\psi : M\to {\Bbb R}$ be a positive
definable function, which is locally bounded from below by positive
constants (for every $x$ in $M$, there exist $c>0$ and a neighborhood
$V$ of $x$ in $M$ such that $\psi>c$ on $V$). Then there exists a
positive definable ${\cal C}^k$ function $\varepsilon : M\to {\Bbb R}$
such that $\varepsilon<\psi$ on $M$.
\end{lemma}

\begin{proof} For $r\in {\Bbb R}$, set 
$$M_r=\{\ x\in M \;;\  \Vert x\Vert^2\leq r\}\;.$$
We can assume $M_{r_0}\not=\emptyset$ for some $r_0$. Observe that
$M_r$ is compact for $r\geq r_0$. For such $r$, we define
$\mu(r)$ by $\mu(r)=\inf\{\psi(x)\;;\ x\in
M_r\}$. The function $\mu:[r_0,+\infty)\to {\Bbb R}$ is
definable and nonincreasing. It is positive, since $M_r$ is covered by
finitely many subsets where  $\psi$ is bounded from below by positive
constants. By Theorem \ref{pwderiv1variable}, there is $a>r_0$ such
that
$\mu$ is ${\cal C}^k$ on an open interval containing $[a,+\infty)$.
Choose a definable
${\cal C}^k$ function $\theta: {\Bbb R}\to {\Bbb R}$ such that
$\theta=0$ on
$(-\infty,a]$, $\theta$ increases from $0$ to $1$ on $[a,a+1]$ and
$\theta=1$ on $[a+1,+\infty)$. Define the function $\mu_1: {\Bbb R}\to
{\Bbb R}$ by
$$\mu_1(r)= \theta(r)\mu(r) + (1-\theta(r))\mu(a+1)\;.$$
Observe that $\mu_1$ is ${\cal C}^k$, definable, and satisfies
$\mu_1\leq\mu$ on $[r_0,+\infty)$. We set $\varepsilon(x)=
{1\over2}\,\mu_1(\Vert x\Vert^2))$ for $x\in M$. By construction, the
positive, definable, ${\cal C}^k$ function $\varepsilon$ satisfies
$\varepsilon<\psi$ on $M$.
\end{proof}

Note that Lemma \ref{Ckminorclosed} also holds if we assume only that
there is a definable ${\cal C}^k$ diffeomorphism from $M$ onto a
closed submanifold of some ${\Bbb R}^m$. Later we shall show that it
holds for every definable submanifold.

\debexo 
Give an explicit formula for $\theta$ in the proof above (you
can take for $\theta$ on $[a,a+1]$ a primitive of a well-chosen
polynomial).
\finexo

\debexo
Having in view the generalization of the lemma to an o-minimal
structure expanding an arbitrary real closed field, replace the
argument for the positivity of $\mu$ in the proof of the preceding
lemma by an argument valid in this general situation. (Hint: if
$\psi(r)=0$, show that there is $b>0$ and a definable path
$\gamma:(0,b)\to M_r$ such that $\psi(\gamma(t))=t$.)
\finexo

\begin{lemma}\label{tubnbhdclosed} The definable tubular neighborhood
\ref{tubnbhd} holds if $M$ is closed in ${\Bbb R}^n$, or definably
${\cal C}^k$ diffeomorphic to a closed submanifold in some ${\Bbb
R}^m$ (e.g.\ if
$M$ is a cell of a ${\cal C}^k$ cdcd).
\end{lemma}

\begin{proof}
We want to avoid the following two ``bad'' situations:
\begin{enumerate}
\item $\varphi$ is not a local diffeomorphism at $(x,v)$.
\item there are $(x,v)$ and $(y,w)$ such that
$\varphi(x,v)=\varphi(y,w)$.
\end{enumerate}
Let $Z$ be the subset of $(x,v)$ in ${\rm N}M$ such that 
$${\rm d}_{(x,v)}\varphi: {\rm T}_{(x,v)}({\rm N}M)\longrightarrow
{\Bbb R}^n$$
is not an isomorphism. The set $Z$ is definable, closed in ${\rm N}M$
and disjoint from the zero section $M\times\{0\}$. For $x$ in $M$,
let $\psi(x)$ be the minimum of $1$, ${\rm dist}((x,0),Z)$ and the
infimum of $r\in {\Bbb R}$ such that 
$$\exists (y,w)\in{\rm N}M\quad \exists v\in {\rm N}_xM\quad
\Vert w\Vert \leq\Vert v \Vert = r \et y+w=x+v\;.$$
The function $\psi: M\to {\Bbb R}$ is definable.\par We claim that it
is
locally bounded from below by positive constants. Take $x\in M$. The
inverse function theorem implies that there is $c_x>0$ such that the
restriction of $\varphi$ to the intersection of ${\rm N}M$ with
the open ball $B((x,0),c_x)\subset {\Bbb R}^{2n}$ is a diffeomorphism
onto a neighborhood of $x$ in ${\Bbb R}^n$. We show that $\psi\geq
c_x/4$ on the intersection of $M$ with the open ball
$B(x,c_x/4)\subset R^n$. We can, of course, take $c_x/4\leq 1$.  Since
$B((x,0),c_x)$ is disjoint from
$Z$ we have $c_x/4\leq {\rm dist}((y,0),Z)$ for every $y\in
M\cap B(x,c_x/4)$. Suppose now that there are $(y,v)$ and $(z,w)$ in
${\rm N}M$ such that $y\in B(x,c_x/4)$, $\Vert w\Vert\leq\Vert
v\Vert<c_x/4$ and $y+v=z+w$. Then $(y,v)$ and $(z,w)$ both belong to
$B((x,0),c_x)$, and we arrive to a contradiction with the fact that
$\varphi$ is injective in restriction to the intersection of this
ball with ${\rm N}M$. Hence, the claim is proved.
We can apply Lemma \ref{Ckminorclosed} to the function $\psi$. We
obtain a definable, positive ${\cal C}^k$ function $\varepsilon: M\to
{\Bbb R}$ such that, in restriction to the definable open subset
$$U=\{(x,v)\in {\rm N}M\;;\ \Vert v\Vert < \varepsilon(x)\}\;,$$
the map $\varphi$ is an injective local diffeomorphism. We conclude
that $\varphi|_U$ is a diffeomorphism onto an open definable
neighborhood of $M$ in ${\Bbb R}^n$.
\end{proof}

\begin{proposition}\label{eqnforsbmnfld}
Let $M$ be a definable ${\cal C}^k$ submanifold of ${\Bbb R}^n$,
closed in ${\Bbb R}^n$ or definably ${\cal C}^k$ diffeomorphic to a
closed submanifold of some ${\Bbb R}^m$. Then $M$ has a nonnegative
definable ${\cal C}^{k-1}$ equation in the complement of $\partial M=
\clos(M)\setminus M$. This means that there exists a
nonnegative definable
${\cal C}^{k-1}$ function $f: {\Bbb R}^n\setminus \partial M\to{\Bbb
R}$ such that $M=f^{-1}(0)$.
Moreover, one can choose $f\leq 1$.
\end{proposition}

\begin{proof} By Lemma \ref{tubnbhdclosed}, we can choose a
definable tubular neighborhood $\Omega$ of $M$ of radius
$\varepsilon: M\to (0,+\infty)$. By Lemma
\ref{Ckminorclosed}, we can assume that $\varepsilon(x) < \min(1,{\rm
dist}(x, \partial M))$ for all $x$ in $M$. Observe that the
square of distance function
$\rho:\Omega\to [0,+\infty)$ is a definable nonnegative ${\cal
C}^{k-1}$ equation of $M$ in $\Omega$. We are going to extend it to
${\Bbb R}^n\setminus \partial M$ by using a partition of unity. We
choose a ${\cal C}^{k-1}$ definable function $\theta:{\Bbb R}\to
[0,1]$ which is $0$ on $(-\infty, 0]$ , $1$ on $[1,+\infty)$ and
increases from $0$ to $1$ on $[0,1]$. Then we set
$$f(x)=\cases{
\theta(2\rho(x)/\varepsilon^2(x))
& on $\Omega$\cr
1&outside $\partial M \cup \{x\in \Omega\;;\ \rho(x)\leq
{1\over2}\varepsilon^2(x)\}$\cr}\;.$$ 
The function $f:{\Bbb
R}^n\setminus\partial M\to {\Bbb R}$ is nonnegative, definable, and
$f^{-1}(0)=M$. The assumption $\varepsilon(x) < \min(1,{\rm dist}(x,
\partial M))$ implies that $\partial M \cup \{x\in \Omega\;;\
\rho(x)\leq {1\over2}\varepsilon^2(x)\}$ is closed (we leave the
proof as an exercise). This shows that
$f$ is indeed ${\cal C}^{k-1}$, and completes the proof of the
proposition.
\end{proof}

\begin{theorem}\label{eqnforclosed} Let $F$ be a closed definable
subset of
${\Bbb R}^n$. For every positive integer $k$, there exists a definable,
continuous, nonnegative function $f: {\Bbb R}^n\to {\Bbb R}$ such
that $f^{-1}(0)= F$ and the restriction of $f$ to ${\Bbb
R}^n\setminus F$ is of class ${\cal C}^k$.
\end{theorem}

\begin{proof} We proceed by induction on the dimension of $F$. If
$F$ is empty, we can take for $f$ any positive constant function.
Now let
$d=\dim F\geq 0$ and  assume that the theorem holds for all closed
definable subsets of ${\Bbb R}^n$ of dimension $<d$. Take a ${\cal
C}^{k+1}$ cdcd of ${\Bbb R}^n$ adapted to $F$. Let ${\cal D}$
be the finite set of cells of dimension $d$ contained in $F$. By
Proposition
\ref{eqnforsbmnfld}, every $d$-cell $C\in {\cal D}$ has a
nonnegative, definable ${\cal C}^k$ equation $g_C$ in ${\Bbb
R}^n\setminus \partial C$. Let $Z$ be the union of all cells of
dimension $<d$ contained in $F$ and all $\partial C$, for $C\in
{\cal D}$. The set $Z$ is definable and closed in ${\Bbb R}^n$.
The product of  all functions $g_C$ for
$C\in{\cal D}$ defines a function $g: U={\Bbb R}^n\setminus Z\to {\Bbb
R}$. The function $g$ is of class ${\cal C}^k$, and it is a definable,
nonnegative equation of
$F\cap U$ in $U$. Moreover, we can take all $g_C$ bounded by $1$ and,
hence, we have $g\leq 1$. By the inductive assumption, there is a
definable, continuous, nonnegative function $h: {\Bbb R}^n\to {\Bbb
R}$ such that $h^{-1}(0)= Z$ and the restriction of $h$ to $U$ is of
class ${\cal C}^k$. The product $f=gh$ is defined and ${\cal C}^k$ on
$U$. It can be continuously extended to ${\Bbb R}^n$ by setting
$f=0$ on $Z$. The zeroset of $F$ is then $Z\cup(F\cap U)=F$.  The
function $f$ satisfies the properties of the theorem.
\end{proof}

\begin{corollary}\label{diffeotoclosed} Let $M$ be a definable ${\cal
C}^k$ submanifold of
${\Bbb R}^n$. Then there is a definable ${\cal C}^k$ submanifold $N$ of
${\Bbb R}^{n+1}$, closed in ${\Bbb R}^{n+1}$, such that the
projection $\pi:{\Bbb R}^{n+1}\to {\Bbb R}^n$ on the first $n$
coordinates induces a diffeomorphism from $N$ onto $M$.
Hence, every  definable ${\cal C}^k$ submanifold of
${\Bbb R}^n$ is definably ${\cal C}^k$ diffeomorphic to a closed
submanifold of ${\Bbb R}^{n+1}$, and Lemma
\ref{Ckminorclosed} and Proposition \ref{eqnforsbmnfld} hold for every
definable ${\cal C}^k$ submanifold.
\end{corollary}

\begin{proof} Let $M$ be a definable ${\cal C}^k$ submanifold of ${\Bbb
R}^n$. Since $M$ is locally closed, $\partial M=\clos(M)\setminus M$
is closed in ${\Bbb R}^n$. By Theorem \ref{eqnforclosed}, there is a
continuous, nonnegative definable function $f:{\Bbb R}^n\to {\Bbb R}$
such that $f^{-1}(0)=\partial M$ and $f$ is ${\cal C}^k$ on ${\Bbb
R}^n\setminus \partial M$. Let
$$N=\{(x,t)\in {\Bbb R}^{n+1}\;;\ x\in M \et t f(x)=1\}\;.$$
Then $N$ is a definable ${\cal C}^k$ submanifold, closed in ${\Bbb
R}^{n+1}$, and the mapping $x\mapsto (x,1/f(x))$ is a definable
${\cal C}^k$ diffeomorphism from $M$ onto $N$, which is the inverse of
$\pi|_N$. 
\end{proof}

The preceding corollary, together with Lemma \ref{tubnbhdclosed},
completes the proof of the definable tubular neighborhood theorem
\ref{tubnbhd}. We conclude with a property which plays an important
role in the course of Macintyre on constructing o-minimal structures.

\begin{proposition}[${\rm D}{\cal C}^k$ - all $k$] Let $A$ be a
definable subset of ${\Bbb R}^n$. For every positive integer $k$,
there exists a definable ${\cal C}^k$ function $f: {\Bbb R}^{n+1}\to
{\Bbb R}$ such that $A=\pi(f^{-1}(0))$, where $\pi={\Bbb R}^{n+1}\to
{\Bbb R}^n$ is the projection on the first $n$ coordinates.
\end{proposition}

\begin{proof} Fix $k$. Take a
${\cal C}^{k+1}$ cdcd of ${\Bbb R}^n$ adapted to $A$. By
corollary \ref{diffeotoclosed}, for every cell $C$ there is a definable
${\cal C}^{k+1}$ submanifold $D$ of
${\Bbb R}^{n+1}$, closed in ${\Bbb R}^{n+1}$, such that $\pi|_D$ is a
diffeomorphism onto $C$. By Proposition \ref{eqnforsbmnfld}, there is
a definable ${\cal C}^{k}$ function $f_C:{\Bbb R}^{n+1}\to {\Bbb R}$
such that $f_C^{-1}(0)=D$. Let $f$ be the product of the functions
$f_C$ for all cells $C$ contained in $A$. Then $f$ is a definable
${\cal C}^k$ function on ${\Bbb R}^{n+1}$ and $A=\pi(f^{-1}(0))$.
\end{proof} 

%% file: OMINroot.bbl
\begin{thebibliography}{[GLRR]}
\addcontentsline{toc}{chapter}{Bibliography}





\bibitem[BS]{BS}
J.L.~Bell, A.R.~Slomson: {\em Models and ultraproducts: An introduction}.
North Holland (1969)


\bibitem[BCR]{BCR} J.~Bochnak, M.~Coste, M-F.~Roy: {\em Real Algebraic Geometry}.
Springer 1998

\bibitem[Co1]{CoTT} M.~Coste: {\em Topological Types of Fewnomials}. In {\em
Singularities Symposium -- {\L}ojasiewicz 70}. Banach Center Pub.\ 44 (1998), 81--92

\bibitem[Co2]{Co} M.~Coste: {\em An Introduction to Semialgebraic Geometry}.  Dottorato
di Ricerca in Matematica, Dip. Mat. Univ. Pisa. Istituti Editoriali e Poligrafici
Internazionali 2000

\bibitem[CR]{CR} M. Coste, M. Reguiat: {\it Trivialit\'es en famille}, dans
{\it Real algebraic geometry}, Lect. Notes Math. 1524, Springer 1992, 193--204

\bibitem[vD]{vdD} L.~van den Dries: {\em Tame Topology and O-minimal
Structures}. London Math.\ Soc.\ Lecture Note 248. Cambridge Univ. Press
1998

\bibitem[vDMi]{vDMi} L. van den Dries, C. Miller: {\it Geometric categories and o-minimal
structures}, Duke Math.\ J. 84 (1996) 497--540

\bibitem[Fu]{Fu} T. Fukuda: {\it Types topologiques des polyn\^omes}, Publ. Math.
I.H.E.S. 46 (1976) 87--106

\bibitem[KPS]{KPS} J. Knight, A. Pillay, C. Steinhorn: {\it Definable sets in o-minimal
structures II}, Trans. Amer. Math. Soc. 295 (1986) 593--605

\bibitem[PS]{PS} A. Pillay, C. Steinhorn: {\it Definable sets in ordered structures I},
Trans. Amer. Math. Soc. 295 (1986) 565--592 

\bibitem[Pr]{Pr} A.~Prestel: {\em Model Theory for the Real Algebraic Geometer}. Dottorato di
Ricerca in Matematica, Dip. Mat. Univ. Pisa. Istituti Editoriali e Poligrafici Internazionali 1998
 
\bibitem[Wi]{Wi} A. Wilkie: {\it Model completeness results for expansion of the real
field by restricted Pfaffian functions and the exponential function}, J. of the Amer.
Math. Soc. 9 (1996) 1051--1094

\end{thebibliography}
